\documentclass[smallextended]{svjour3}       
\smartqed  

\usepackage{graphicx}
\usepackage{amssymb}
\usepackage{amsmath}
\usepackage{float}
\usepackage{subfigure}
\usepackage{hyperref}
\usepackage{prettyref}
\usepackage{setspace}
\usepackage{mathtools}
\usepackage{amsfonts}
\usepackage{enumerate}
\usepackage{algorithm}
\usepackage{algorithmic}
\DeclareMathOperator*{\esssup}{ess\ sup}

\usepackage{amsmath}

\begin{document}

\title{A non-local gradient based approach of infinity Laplacian with  $\operatorname{\mathbf{\Gamma}}$-convergence
\thanks{This work was supported by the National Natural Science Foundation of China (12071244).}
}

\titlerunning{Infinity Laplacian with $\mathrm{\Gamma}$-convergence}        

\author{Weiye Gan \and Xintong Liu \and Yicheng Li \and Zuoqiang Shi}


\institute{Weiye Gan \at
              Department of Mathematical Sciences, Tsinghua University, 100084 Beijing, China \\
              \email{gwy18@mails.tsinghua.edu.cn}           
          \and
          Xintong Liu \at
              Department of Mathematical Sciences, Tsinghua University, 100084 Beijing, China\\
              \email{liu-xt18@mails.tsinghua.edu.cn}
            \and
            Yicheng Li \at
            Department of Mathematical Sciences, Tsinghua University, 100084 Beijing, China\\
            \email{liyc18@mails.tsinghua.edu.cn}
            \and
            Zuoqiang Shi \at
            Department of Mathematical Sciences, Tsinghua University, 100084 Beijing, China\\
            \email{zqshi@tsinghua.edu.cn}
        }

\date{}

\maketitle

\begin{abstract}
We propose an infinity Laplacian method to address the problem of interpolation on an unstructured point cloud. In doing so, we find the labeling function with the smallest infinity norm of its gradient. By introducing the non-local gradient, the continuous functional is approximated with a discrete form. The discrete problem is convex and can be solved efficiently with the split Bregman method. Experimental results indicate that our approach provides consistent interpolations and the labeling functions obtained are globally smooth, even in the case of extreme low sampling rate. More importantly, convergence of the discrete minimizer to the optimal continuous labeling function is proved using $\mathrm{\Gamma}$-convergence and compactness, which guarantees the reliability of the infinity Laplacian method in various potential applications.
\end{abstract}

\section{Introduction}

Consider a point cloud $P = S \cup T$ in $\mathbb{R}^d$, in which $S = \{x_i\}_{i=1}^m$ is the set of labeled points and $T = \{x_{m+j}\}_{j=1}^n$ is the set of unlabeled points. We denote by $u_i$ the label value of the point $x_i$ for $i=1,2,\dots,m$. The task of predicting the label values of points in $T$, widely known as interpolation or semi-supervised learning, has received much attention in recent years due to its various applications in the fields of data science \cite{NPDETOWGIP,OPIG}, image processing \cite{NPIP,LDMM}, and machine learning \cite{NILEGAIM,WNLL}. The interpolation task can be completed by finding some labeling function $u: \mathcal{D}\rightarrow \mathbb{R}$, in which $P \subseteq \mathcal{D}$ and the values of $u$ on the set $S$ is given by the label values. By considering the regularization of $u$, the labeling function can be determined as the minimizer of certain functional. 
It was observed in many tasks that points close to each other share similar label values \cite{LDMM,WNLL,CURE,LDMMcolorization}. Thus a natural idea is to minimize the $L_2$ norm of the gradient of the labeling function, which leads to the popular Graph Laplacian (GL) method in machine learning \cite{GL}. However, it was found that when the number of unlabeled points goes to infinity, harmonic extension may fail and the resulting labeling function is discontinuous \cite{INFLIMGL,HOR}. To deal with this issue, the authors of \cite{WNLL} propose a Weighted Non-local Laplacian (WNLL) method that minimizes a weighted $L_2$ norm and obtain much better results in handwritten digit classification and image inpainting tasks. Intuitively, if $u \in W^{k,p}(\mathbb{R}^n$), the Soblev embedding theorem implies that $u$ will be H\"older continuous with positive component $k-n/p$ when $kp > n$. In high dimensions, this is not the case for the GL method, where $k = 1$ and $p =2$. For the case of $k = 2$, curvature can be used as regularization of the labeling function \cite{CURE}. 

The condition $kp > n$ can be fulfilled by increasing either $k$ or $p$. In this paper, we study the case of $k = 1$, $p = \infty$ and propose a new method that minimizes the $L_\infty$ norm of the gradient for the interpolation task. We assume that the point cloud $P$ samples a smooth manifold $\mathcal{M} \subset \mathbb{R}^d$. We also find a labeling function $u:\mathcal{M} \rightarrow \mathbb{R}$ and determine it by minimizing the infinity norm of its gradient. Besides, we also add the $L_2$ norm of the gradient of $u$ as regularization for a unique solution. The functional we minimize writes
\begin{equation}
\label{con_functional}
\mathcal{J}(u) = \Vert \nabla _\mathcal{M} u \Vert_\infty^p  + \lambda \Vert \nabla _\mathcal{M} u \Vert_{L_2}^2
\end{equation}
in which $p > 1$ is a constant, $\mathcal{M}$ is the manifold that the point cloud $P$ samples. $\lambda$ is a fixed parameter that balances the infinity norm and the $L_2$ norm. Note that the values of $u$ at the label points are fixed with the label values, we find the labeling function by solving the following constrained optimization problem.
\begin{equation}
\label{opt_problem}
\begin{aligned}
\mathop{\min}_u \quad &\Vert \nabla _\mathcal{M} u \Vert_\infty^p  + \lambda \Vert \nabla _\mathcal{M} u \Vert_{L_2}^2\\
s.t.\quad & u(x_i) = u_i, \quad i = 1,\dots,m
\end{aligned}
\end{equation}

When $\mathcal{M}$ is an open set in $\mathbb{R}^n$, $S = \partial \mathcal{M}$ and $\lambda = 0$, the absolute minimizer among all Lipschitz continuous functions of this problem is the solution of the infinity Laplacian equation with Dirichlet boundary condition \cite{ILAEG}.
\begin{equation*}
\mathrm{\Delta}_{\infty} u:=\left\langle D u, D^{2} u D u\right\rangle= \sum_{i,j=1} ^n u_{x_{i}} u_{x_{j}} u_{x_{i} x_{j}}=0
\end{equation*}
In real applications, however, the manifold $\mathcal{M}$ is not known, making it impossible to determine the function $u$ explicitly. Instead of finding an analytic solution, we aim at determining values of $u$ on the unlabel set $T$. In doing so, we approximate the gradient with the non-local gradient \cite{NPIP} defined as
\begin{equation}
\label{non-local gradient}
    Du(x,y) = \sqrt{w(x,y)}\big(u(x)-u(y)\big)
\end{equation}
for all pairs of points $(x,y) \in \mathcal{M} \times \mathcal{M}$. In this definition, $w(x,y) = \eta(\Vert x-y \Vert)$ and $\eta$ is a monotonously decreasing kernel function satisfying some regularity conditions. The infinity norm and $L_2$ norm of $\nabla_{\mathcal{M}}u$ are then approximated as
\begin{equation*}
    \Vert \nabla_\mathcal{M}u \Vert_{\infty} = \mathop{\max}_{x\in \mathcal{M}} \bigg( \int_\mathcal{M}\Big \vert Du(x,y)\Big \vert^p d y \bigg)^{1/p}
\end{equation*}
and
\begin{equation*}
    \Vert \nabla_\mathcal{M}u \Vert_{L_2} =  \Bigg(  \int_\mathcal{M} \int_\mathcal{M}\Big \vert Du(x,y)\Big \vert^2 d x d y \Bigg) ^{1/2}
\end{equation*}
Using these approximations, the constrained optimization problem \eqref{opt_problem} is reformulated as

\begin{equation}
\label{final opt}
\begin{aligned}
\mathop{\min}_u \quad & \mathop{\max}_{x\in \mathcal{M}} \int_\mathcal{M}\Big \vert Du(x,y)\Big \vert^p d y + \lambda  \int_\mathcal{M} \int_\mathcal{M}\Big \vert Du(x,y)\Big \vert^2 d x d y\\
s.t.\quad & u(x_i) = u_i,\quad i = 1,\dots,m
\end{aligned}
\end{equation}
in which $Du$ is the non-local gradient defined in \eqref{non-local gradient}. Assuming that the point cloud $P$ samples the manifold $\mathcal{M}$ uniformly, we obtain the discrete version of problem \eqref{final opt} as follows
\begin{equation}
\label{prob_cons}
\mathop{\min}_{\{u_{m+j}\}_{j=1}^n} \Big\{\mathop{\max}_{i\in I} \big\{\sum_{j \in I} w_{ij}\vert u_i-u_j \vert ^p\big\}  + \alpha \sum_{i\in I} \sum_{j \in I} w_{ij}(u_i-u_j)^2\Big\}\\
\end{equation}
in which $I = \{1,\dots, m+n\}$. The constant $\alpha$ is the multiplication of $\lambda$ and the measure at each point in $P$. $u_{m+j} = u(x_{m+j})$ is the label value of $x_{m+j} \in T$. $w_{ij}$ is the weight of points $x_i$ and $x_j$. In this way, only the label values and the weights are needed to fulfill the interpolation task. For an efficient numerical algorithm, we choose $p = 2$, in which case the objective function is convex and can be minimized efficiently using the split Bregman method \cite{SBI}. In the following we call the proposed approach the infinity Laplacian method and abbreviate it to IL.

A crucial problem we need to consider is whether the problem \eqref{prob_cons} is a good discrete approximation of \eqref{opt_problem}. For this, we will show that the minimizers of \eqref{prob_cons} converges to the minimizer of \eqref{opt_problem} with respect to $L^\infty$ norm when the number of unlabeled points goes to infinity. The main tool we use is called $\mathrm{\Gamma}$-convergence for functionals which is popular in studying the consistency of graph-based methods for semi-supervised learning, see~\cite{ref1,ref2,trillos2016continuum}. In~\cite{ref2}, the authors consider the convergence of the following functional
\begin{equation*}
    F_{p,p}(u)=\frac{1}{|I|^2}\sum_{i,j\in I}w_{ij}|u_i-u_j|^p
\end{equation*}
for $1<p<\infty$ and remark that there are no essential difficulties to extend the results to the manifolds which have already been done for related problems about the graph Laplacian~\cite{trillos2020error}. When $p$ is smaller than the ambient space dimension $d$, however, such discrete method called $p$-Laplacian regularization comes with the drawback that the label information will be ignored if the number of unlabel data tends to infinity. For this motivation,~\cite{ref1} consider the situation for $p=\infty$, that is, the functional
\begin{equation*}
    F_{\infty,\infty}(u)=\mathop{\max}_{i\in I}\mathop{\max}_{j\in I}w_{ij}|u_i-u_j|
\end{equation*}
While in this work, we will break the symmetry of $i,j$ to consider the functional like
\begin{equation*}
    F_{\infty,p}(u)=\mathop{\max}_{i\in I} \big\{\frac{1}{|I|}\sum_{j \in I} w_{ij}|u_i-u_j|^p\big\}
\end{equation*}
and show that $F_{\infty,p}$ converges to a functional whose value is the $L^\infty$ norm of a function multiplying a constant similarly to $F_{\infty,\infty}$.

The rest of the paper is organized as follows. In \autoref{sec:Preliminaries}, we will give a brief introduction to the main mathematical tools we use. In \autoref{sec:main}, we firstly consider $\mathcal{M}$ as a domain in $\mathbb{R}^d$ and prove that the minimizers of the discrete problem \eqref{prob_cons} converges to the minimizer of the continuous version \eqref{con_functional} when the number of unlabeled points goes to infinity. In \autoref{sec:manifold}, we will generalize our convergence results to the manifolds. In \autoref{sec:alg}, we propose an algorithm to solve the discrete problem \eqref{prob_cons}. In \autoref{sec:experiments}, two toy examples are used to test the proposed algorithm. Conclusions are made in \autoref{sec:conclusion}.

\section{Preliminaries}
\label{sec:Preliminaries}
This section reviews two mathematical tools used in this article. The first one is the concept of $\mathrm{\Gamma}$-convergence which we have mentioned in the introduction. And the second one is the transportation map between probability measures which we employ in order to turn functions on samples into continuum ones. Such connection between discrete function and continuum ones is necessary for comparing them. Although our conclusion is more similar to~\cite{ref1} considering the case $p=\infty$, we adopt the way used in~\cite{ref2} not~\cite{ref1} which use the closest point projection. From this point of view, the problem we studied is like a combination of the problems studied in the previous two articles. 
\subsection{$\mathrm{\Gamma}$-Convergence}
	In this part, we will introduce the concept of $\mathrm{\Gamma}$-convergence and one of its important properties-implying the convergence of minimizers under compactness assumptions. More detailed introduction about $\mathrm{\Gamma}$-convergence can be found in~\cite{braides2002gamma}.
	\begin{definition}[$\mathrm{\Gamma}$-convergence]
		\label{def:gamma}
		Let $X$ be a metric space and $F_n:X\rightarrow\mathbb{R}\cup\{-\infty,\infty\}$ be a sequence of functionals on $X$. We say that $F_n$ $\mathrm{\Gamma}$-converges to $F:X\rightarrow\mathbb{R}\cup\{-\infty,\infty\}$ if\\
		\begin{enumerate}
			\item[(1)](liminf inequality) for any sequence $\{x_n\}_{n\in\mathbb{N}}\subset X$ converging to $x\in X$,
			\[\liminf_{n\rightarrow\infty}F_{n}(x_n)\geq F(x)\]
			\item[(2)](limsup inequality) for any $x\in X$, there exists a sequence $\{x_n\}_{n\in\mathbb{N}}\subset X$ converging to $x$ such that
			\[\limsup_{n\rightarrow\infty}F_n(x_n)\leq F(x)\]  
		\end{enumerate}
	which also denoted by $F_n\stackrel{\Gamma}{\rightarrow}F(n\rightarrow\infty)$.
	\end{definition}
	\begin{lemma}[Convergence of minimizers]
		\label{lem:minimizer}
		Let $X$ be a metric space and $F_n:X\rightarrow[0,\infty]$ $\mathrm{\Gamma}$-converges to $F:X\rightarrow[0,\infty]$ which is not identically $\infty$. If there exists a relatively compact sequence $\{x_n\}_{n\in\mathbb{N}}\subset X$ such that
		\[\lim_{n\rightarrow\infty}(F_n(x_n)-\inf_{x\in X}F_n(x))=0\]
		then we have
		\[\lim_{n\rightarrow\infty}\inf_{x\in X}F_n(x)=\min_{x\in X}F(x)\]
		and any cluster point of $\{x_n\}_{n\in\mathbb{N}}$ is a minimizer of $F$.
	\end{lemma}
	\begin{proof}
		For any $y\in X$, we know that there exists a sequence $\{y_n\}_{n\in\mathbb{N}}\subset X$ satisfying the limsup inequality. So we have
		\[F(y)\geq\limsup_{n\rightarrow\infty}F_n(y_n)\geq\limsup_{n\rightarrow\infty}\inf_{x\in X}F_n(x)\]
		which yields
		\[\min_{x\in X}F(x)\geq\limsup_{n\rightarrow\infty}\inf_{x\in X}F_n(x)\]
		On the other hand, consider the sequence $\{x_n\}_{n\in\mathbb{N}}\subset X$ mentioned in the assumption, let $\tilde{x}$ be one of the cluster points of $\{x_n\}_{n\in\mathbb{N}}$, using the liminf ineuqality, we get
		\[\liminf_{n\rightarrow\infty}\inf_{x\in X}F_n(x)=\liminf_{n\rightarrow\infty}F_n(x_n)\geq F(\tilde{x})\geq\min_{x\in X}F(x)\]
		Therefore,
		\[\limsup_{n\rightarrow\infty}\inf_{x\in X}F_n(x)\leq\min_{x\in X}F(x)\leq F(\tilde{x})\leq\liminf_{n\rightarrow\infty}\inf_{x\in X}F_n(x)\]
		and we can get the conclusion.
	\end{proof}
\subsection{Transprotation Map}
	Let $\Omega\subset\mathbb{R}^d$ and $\mu,\nu$ be two probability measures on $\Omega$, a map $T:\Omega\rightarrow\Omega$ called a transportation map from $\nu$ to $\mu$ if it satisfies the push-forward condition
	\[\mu=T_{\#}\nu=\nu\circ T^{-1},\]
	In this arcticle, we always let $\mu$ be a empirical measure of $\nu$. That is, if we have $\{x_i\}_{i=1}^n$ independently sampled from $\Omega$ with distribution $\nu$, then
	\[\mu(A)=\frac{1}{n}\sum_{i=1}^n\delta_{x_i}(A),\quad A\in\mathcal{B}(\Omega)\]
	where $\delta_{x_i}(A):=\mathbf{1}_A(x_i)$ is the Dirac measure of $x_i$. In this case, the push-forward condition implies the following transformation:
	\begin{equation}
	\label{eq:trans}
		\frac{1}{n}\sum_{i=1}^nu(x_i)=\int_\Omega u(x)d\mu(x)=\int_\Omega u(T(x))d\nu(x)
	\end{equation}
	In our proof, we will also use a simple property of transportation map about the essential supremum:
	\begin{lemma}
	\label{lem:esssup}
		 Let $\nu,\mu$ be two probability measures on $\Omega$ and $T:\Omega\rightarrow\Omega$ be a transportation map from $\nu$ to $\mu$ then for any measurable function $f:\Omega\rightarrow\mathbb{R}$, we have that
		 \[\mu\text{-}\esssup_{x\in\Omega}f(x)=\nu\text{-}\esssup_{y\in\Omega}f(T(y))\]		 
	\end{lemma}
	\begin{proof}
		Note that $A$ is a null set with respect to $\mu$ if and only if $T^{-1}(A)$ is a null set with respect to $\nu$, from which we can get the conclusion.
	\end{proof}
	\indent In general, the transportation map $T$ from $\nu$ to $\mu$ is not unique. The infimum of the $L^\infty$ distance between different $T$ and the identity actually gives a distance between $\nu$ and $\mu$, so-called the $\infty$-Wasserstein distance, or $\infty$-optimal transportation distance.
	\[d_\infty(\nu,\mu)=\inf_{T_\#\nu=\mu}\nu\text{-}\esssup_{x\in\Omega}|x-T(x)|\]
	When $\mu$ is the empirical measure of $\nu$, we have the following estimates on the $\infty$-Wasserstein distance, which is also used as Thm 3.3 in~\cite{ref2}. Such estimates are established in~\cite{trillos2015rate} for $d\geq 2$, which extend the results in~\cite{ajtai1984optimal,leighton1989tight,shor1991minimax,talagrand2014upper}. 
	\begin{theorem}
		\label{thm:transportation}
		Let $\Omega\subset\mathbb{R}^d$ be a bounded domain with Lipschitz boundary. Let $\nu$ be a probability measure on $\Omega$ with density(with respect to Lebesgue measure) $\rho$ which is bounded above and below by positive constants. Let $x_1,x_2,\dots$ be a sequence of independent random samples with distribution $\nu$ and let $\nu_n$ be the empirical measure. Then, there exists constant $C\geq c\geq0$ such that almost surely there exists a sequence of transprotation maps $\{T_n\}_{n=1}^\infty$ from $\nu$ to $\nu_n$ with the property
		\begin{equation}
		c\leq
		\liminf_{n\rightarrow\infty}\frac{\lVert T_n-Id\rVert_{L^\infty(\Omega)}}{\delta_n}
		\leq
		\limsup_{n\rightarrow\infty}\frac{\lVert T_n-Id\rVert_{L^\infty(\Omega)}}{\delta_n}\leq C
		\label{eq:e2}
		\end{equation}
		where
		\begin{equation*}
		\begin{aligned}
		\delta_n=\left\{\begin{array}{cc}
		\sqrt{\frac{\ln\ln(n)}{n}}                  & \text{if }d=1,    \\
		\frac{(\ln n)^\frac{3}{4}}{\sqrt{n}}        & \text{if }d=2,    \\
		\frac{(\ln n)^\frac{1}{d}}{n^{\frac{1}{d}}} & \text{if }d\geq 3
		\end{array}\right.
		\end{aligned}
		\end{equation*}
	\end{theorem}

\section{Theoretical Analysis for Domains}
\label{sec:main}
    In this section, we will prove the convergence results from the discrete functionals to the continuum functional when $\mathcal{M}$ is actually a domain in $\mathbb{R}^d$. In order to distinguish from the situation of manifolds, we will rewrite $\mathcal{M}$ as $\Omega$. We consider the convergence results on domains at first because this is the main difficulty of the problems and in the next section, we will show that the results for smooth manifolds can be derived from the results for domains and their local coordinate representations. The convergence of functionals we consider is called $\mathrm{\Gamma}$-convergence defined in \autoref{def:gamma}. The compactness of minimizing sequence of the discrete functionals will also be proved. Then, because of \autoref{lem:minimizer}, we will get the convergence of the minimizers(\autoref{thm:2}). We firstly list the assumptions we need as follows.

    Let $p\in(1,\infty)$ and $\Omega\subset\mathbb{R}^d$ be a bounded  $C^2$  domain.
    Denote by $\nu$ the $\mathbb{R}^d$-Lebesgue measure restricted on $\Omega$.
    Let $\{x_i\}_{i=1}^N$ be arbitrary $N$ points in $\Omega$ and $\{x_i\}_{i=N+1}^\infty$ be an i.i.d. sequence sampled from $\frac{\nu}{\nu(\Omega)}$ which is the probability measure induced by $\nu$.
    Denote by $\mathcal{O}\coloneqq\{x_i\big|1\leq i\leq N\}$ the label set and $\Omega_n\coloneqq\{x_i\big|1\leq i\leq |\Omega_n|\}$
    satisfying
    $N<|\Omega_1|<|\Omega_2|<\dots<\infty$.
    We have a labeling function $g:\bar{\Omega}\longrightarrow\mathbb{R}$ which is Lipschitz continuous.
    For $x\in\mathcal{O}$, we have known that it is labeled by $g(x)$.

    Denote by $\nu_n$ the empirical measure with respect to $\Omega_n$, which means
    \[\nu_n(A)=\frac{1}{|\Omega_n|}\sum_{x\in\Omega_n}\delta_x(A),\quad A\in\mathcal{B}(\Omega)\]
    Here, $\delta_x(A):=\mathbf{1}_A(x)$ is the Dirac measure of x. From \autoref{thm:transportation}, we know that with probability one there exists a
    sequence of transportation maps $\{T_n\}_{n=1}^\infty$ from $\frac{\nu}{\nu(\Omega)}$ to
    $\nu_n$ satisfying~\eqref{eq:e2}. And from~\eqref{eq:trans}, we have that
    \[\frac{1}{|\Omega_n|}\sum_{x\in\Omega_n}u(x)=\frac{1}{\nu(\Omega)}\int_\Omega u(T_n(x))d\nu(x)=\frac{1}{\nu(\Omega)}\int_\Omega u(T_n(x))dx\]
    for a function $u:\Omega\longrightarrow\mathbb{R}$. From this equation, we have $T_n(x)\in\Omega_n$ almost surely. And $\Omega_n$ is a finite set so it is also a $\nu$-null set. So we can let $T_n(\Omega)\subset\Omega_n$ and $T_n=Id$ on $\Omega_n$ without changing the properties $\{T_n\}$ already have.

    We consider $w(x,y)=\eta_{s_n}(|x-y|)$ where $\eta$ is called the kernel function that satisfies the following standard assumptions:
    \begin{enumerate}
        \item[(K1)] $\eta:[0,\infty)\longrightarrow[0,\infty)$ is continuous,
        \item[(K2)] $\eta$ is monotonously decreasing,
        \item[(K3)] $ \mathrm{supp}(\eta)\subset[0,r_{\eta}]$ for some $r_{\eta}>0$.
    \end{enumerate}
    and $\eta_{s_n}(x)\coloneqq\frac{1}{s_n^d}\eta(\frac{x}{s_n})$ is the scaled kernel with the scaling parameter $s_n>0$.
    
    We define a constant that is only related to the kernel $\eta$ and $p$ which will be used in our results, that is:
    \begin{equation}
        \label{eq:kernel-radius}
        \sigma_{\eta}\coloneqq(\int_{B(r_{\eta})}\eta(|z|)|z\cdot e_1|^p dz)^\frac{1}{p}
    \end{equation}
    where $e_1\coloneqq(1,0,\dots,0) \in \mathbb{R}^d$ and $B(r_{\eta}):=\{x\in\mathbb{R}^d\big||x|<r_\eta\}$, note that $\eta(|z|)$ is a radially symmetric function of $z$, we have for any unit vector $h\in\mathbb{R}^d$,
    \[(\int_{B(r_{\eta})}\eta(|z|)|z\cdot h|^p dz)^\frac{1}{p}=\sigma_{\eta}.\]

    Finally, we define the discrete functional $E_n:L^\infty(\Omega_n)\longrightarrow\mathbb{R}$,
\begin{equation}
        \label{eq:discrete-func}
        E_n(u)\coloneqq
        \begin{cases}
            \frac{1}{s_n}\mathop{\max}\limits_{x\in\Omega_n}(\frac{1}{|\Omega_n|}\sum_{y\in\Omega_n}\eta_{s_n}(|x-y|)|u(x)-u(y)|^p)^\frac{1}{p}
            & \text{if } u=g \text{ on } \mathcal{O},\\
            \infty
            & \text{otherwise.}
        \end{cases}
\end{equation}
    Using the transportation maps $\{T_n\}$, we can extend the domain of $E_n$ to $L^\infty(\Omega)$,
    \begin{equation}
        E_{n,\mathrm{cons}}(u)\coloneqq
        \begin{cases}
            E_{n}(\tilde{u}) & \text{if } u=\tilde{u}\circ T_n \text{ for some } \tilde{u}:\Omega_n\longrightarrow\mathbb{R}, \\
            \infty           & \text{otherwise}.
        \end{cases}
        \label{eq:e3}
    \end{equation}
    We also define a continuum functional $\mathcal{E}(u):L^\infty(\Omega)\longrightarrow\mathbb{R}$,
    \begin{equation*}
        \label{eq:continuum-func}
        \mathcal{E}(u)\coloneqq
        \begin{cases}
            \esssup_{x\in\Omega}|\nabla u| & \text{if } u\in W^{1,\infty}(\Omega), \\
            \infty                         & \text{otherwise}.
        \end{cases}
    \end{equation*}
    and its constrained version
    \begin{equation}
        \mathcal{E}_{\mathrm{cons}}(u)\coloneqq\left\{\begin{array}{ll}
                                                          \mathcal{E}(u) \ \  & \text{if } u=g\text{ on } \mathcal{O}, \\
                                                          \infty         & \text{otherwise}.
        \end{array}\right.
        \label{eq:e4}
    \end{equation}

\sloppy Now we can state our main results as the following two theorems:

    \begin{theorem}
        [Discrete to continuum $\mathrm{\Gamma}$-convergence]
        \label{thm:1}
        Let $\Omega\subset\mathbb{R}^d$ be a bounded domain of $C^2$ class,
        and let the kernel $\eta$ fulfil (K1)-(K3).
        Then for any null sequence $\{s_n\}_{n\in\mathbb{N}}\subset\mathbb{R}^+$ which satisfies $\lim_{n\rightarrow\infty}\frac{\delta_n}{s_n}=0$,
        with probability one we have
        \begin{enumerate}
            \item[(1)] for any sequence $\{u_n\}_{n\in\mathbb{N}}\subset L^\infty(\Omega)$ converging to $u\in L^\infty(\Omega)$ in $L^\infty(\Omega)$,
            \[\liminf_{n\rightarrow\infty}E_{n,\mathrm{cons}}(u_n)\geq\sigma_\eta\nu(\Omega)^{-\frac{1}{p}} \mathcal{E}_{\mathrm{cons}}(u)\]
            \item[(2)] for any $u\in L^\infty(\Omega)$, there exists a sequence $\{u_n\}_{n\in\mathbb{N}}\subset L^\infty(\Omega)$ converging to $u$ in $L^\infty(\Omega)$ such that
            \[\limsup_{n\rightarrow\infty}E_{n,\mathrm{cons}}(u_n)\leq\sigma_\eta\nu(\Omega)^{-\frac{1}{p}} \mathcal{E}_{\mathrm{cons}}(u)\]
        \end{enumerate}
With \autoref{def:gamma}, another statement of properties (i) and (ii) is that the discrete functional sequence $\{E_{n,\mathrm{cons}}\}_{n\in\mathbb{N}}$ $\mathrm{\Gamma}$-converges to the continuum functional $\sigma_\eta\nu(\Omega)^{-\frac{1}{p}}\mathcal{E}_{\mathrm{cons}}$.
    \end{theorem}

    \begin{theorem}[Convergence of minimizers]
        \label{thm:2}
        Let domain $\Omega$, kernel $\eta$ and null sequence\\ $\{s_n\}_{n\in\mathbb{N}}\subset\mathbb{R}^+$
        satisfies the assumptions in \autoref{thm:1}.
        Then with probability one, any bounded sequence $(u_n)_{n\in \mathbb{N}}$ such that
        \begin{equation*}
            \lim_{n\rightarrow \infty} \left( E_{n,\mathrm{cons}}(u_n) - \inf_{u\in L^\infty(\Omega)} E_{n,\mathrm{cons}}(u)
            \right) = 0
        \end{equation*}
        is relatively compact in $ L^\infty(\Omega)$ and
        \begin{equation*}
            \lim_{n\rightarrow \infty} E_{n,\mathrm{cons}}(u_n) = \sigma_\eta\nu(\Omega)^{-\frac{1}{p}}\min_{u\in L^\infty(\Omega)} \mathcal{E}_{\mathrm{cons}}(u).
        \end{equation*}
        Furthermore, every cluster point of $(u_n)_{n\in \mathbb{N}}$ is a minimizer of $ \mathcal{E}_{\mathrm{cons}}$.
    \end{theorem}
    \begin{remark}
        In the above theorems, we only consider the convergence of the first term in~\eqref{prob_cons} since the convergence of the second term is just a special case in~\cite{ref2}(choose $p=2$). It can be easily verified that the $\mathrm{\Gamma}$-convergence of the two parts are sufficient to derive the $\mathrm{\Gamma}$-convergence of the whole functional. For liminf inequality, this is because we consider a bounded domain so $\{u_n\}$ which converges to $u$ in $L^{\infty}(\Omega)$ also converges to $u$ in $L^2(\Omega)$. For limsup inequality, we can choose the same $\{u_n\}$ for both terms.
    \end{remark}

    
    \paragraph{Proof outline}
    In \autoref{sec:gamma-convergence} we will give the proof of \autoref{thm:1}.
    The proof will be broken down into two parts based on the idea in~\cite{ref1}.
    First, in \autoref{subsec:non-local-to-local-convergence},
    we will prove the $\mathrm{\Gamma}$-convergence of the non-local functionals
    \begin{equation*}
        \mathcal{E}_{s}(u)\coloneqq\frac{1}{s}\esssup_{x\in\Omega}(\int_{\Omega}\eta_{s}(|x-y|)|u(x)-u(y)|^p dy)^\frac{1}{p}
    \end{equation*}
    to the continuum functional $\mathcal{E}$ defined in~\eqref{eq:continuum-func}.
    Second, in \autoref{subsec:discrete-to-continuum-convergence},
    we will utilize the transportation map to bridge the gap between the non-local functional
    $\mathcal{E}_{s,\mathrm{cons}}$ and the discrete functional $E_{s,\mathrm{cons}}$, and establish
    discrete to continuum $\mathrm{\Gamma}$-convergence as our first main result.

    In \autoref{sec:compactness} we will prove \autoref{thm:2}.
    We first prove the compactness result for the non-local functionals $\mathcal{E}_{s,\mathrm{cons}}$,
    and then with the same estimation derived in \autoref{subsec:discrete-to-continuum-convergence}, we
    get the compactness result and convergence of minimizers as a corollary.

    \subsection{$\mathrm{\Gamma}$-Convergence of Functionals}\label{sec:gamma-convergence}

    \subsubsection{Non-Local to Local Convergence}\label{subsec:non-local-to-local-convergence}
    We first define functional $\mathcal{E}_{s}:L^\infty(\Omega)\longrightarrow\mathbb{R}$,
    \begin{equation}
        \mathcal{E}_{s}(u)\coloneqq\frac{1}{s}\esssup_{x\in\Omega}(\int_{\Omega}\eta_{s}(|x-y|)|u(x)-u(y)|^p dy)^\frac{1}{p}
        \label{eq:e5}
    \end{equation}
    Then for any $u_1,u_2\in L^\infty(\Omega)$ satisfying $\lVert u_1\rVert_{L^\infty},\lVert u_2\rVert_{L^\infty}\leq M$ where $M$ is a positive constant, we have
    \begin{equation*}
        \begin{aligned}
        |\mathcal{E}_{s}(u_1)-\mathcal{E}_{s}(u_2)|&\leq \frac{C_{M,p,\eta}}{s}\esssup_{x\in\Omega}(\int_{\Omega}\eta_{s}(|x-y|)||u_1(x)-u_1(y)|-|u_2(x)-u_2(y)||^pdy)^\frac{1}{p}\\
        &\leq\frac{2C_{M,p,\eta}}{s}\esssup_{x\in\Omega}(\int_{\Omega}\eta_{s}(|x-y|)dy)^\frac{1}{p}\lVert u_1-u_2\rVert_{L^\infty}\\
        &\leq\frac{2C_{M,p,\eta}}{s}(\int_{B(sr_\eta)}\eta_{s}(|y|)dy)^\frac{1}{p}\lVert u_1-u_2\rVert_{L^\infty}\\
        &=\frac{2C_{M,p,\eta}}{s}(\int_{B(r_\eta)}\eta(|y|)dy)^\frac{1}{p}\lVert u_1-u_2\rVert_{L^\infty}
    \end{aligned}
    \end{equation*}
    which means
    \begin{equation}
        |\mathcal{E}_s(u_1)-\mathcal{E}_{s}(u_2)|\leq \tilde{C}_{M,p,\eta}\frac{1}{s}\lVert u_1-u_2\rVert_{L^\infty}
        \label{eq:e11}
    \end{equation}
    This estimation will be used in the proof of Lemma~\ref{lem:2.5}. And in the rest of this part we will prove $\mathcal{E}_s$ $\mathrm{\Gamma}$-converges to $\sigma_\eta\mathcal{E}$ as $s\rightarrow 0$.

    \begin{lemma}
    	\label{lem:2.2}
        Let $\Omega$ be a domain in $\mathbb{R}^d$ and kernel $\eta$ satisfy (K1)-(K3), then
        \[\liminf_{n\rightarrow\infty}\mathcal{E}_{s_n}(u_n)\geq \sigma_{\eta}\mathcal{E}(u)\]
        for all $\mathbb{R}^+\ni s_n\longrightarrow 0,\ L^\infty(\Omega)\ni u_n\longrightarrow u$ in $L^\infty(\Omega)$
    \end{lemma}
    \begin{proof}
        We assume w.l.o.g that
        \[\liminf_{n\rightarrow\infty}\mathcal{E}_{s_n}(u_n)<\infty\]
        Denote by $h=h_x$ a vector in $\mathbb{R}^d$ only related to $x\in \mathbb{R}^d$, which satisfies that 
        \[\lVert h_x\rVert_{L^\infty(\mathbb{R}^d)}<\infty\]
        Define two functions of $x$
        \[D_n^\pm(x)\coloneqq\int_{\Omega}\eta_{s_n}(|x-y|)^\frac{p-1}{p}\left|\frac{(x-y)\cdot h}{s_n}\right|^{p-1}\eta_{h,n}^\pm(x-y)\frac{u_n(x)-u_n(y)}{s_n}dy\]
        where
        \[\eta_{h,n}^\pm(x-y)\coloneqq\left\{\begin{array}{cc}
                                                 \eta_{s_n}(|x-y|)^\frac{1}{p} & \text{if }\pm(x-y)\cdot h>0, \\
                                                 0                             & \text{else}
        \end{array}\right.\]
        Note that with such definitions,
        \[|\eta_{h,n}^+(x-y)|+|\eta_{h,n}^-(x-y)|=\eta_{s_n}(x)^\frac{1}{p}\ a.s.\]
So, combining with the H\"older inequality, we have
\begin{equation*}
    \begin{aligned}
    &|D_n^+(x)|+|D_n^-(x)|\\
    \leq&\int_{\Omega}\eta_{s_n}(|x-y|)^\frac{p-1}{p}\left|\frac{(x-y)\cdot h}{s_n}\right|^{p-1}\\
    &\cdot(|\eta_{h,n}^+(x-y)|+|\eta_{h,n}^-(x-y)|)\left|\frac{u_n(x)-u_n(y)}{s_n}\right|dy\\
    =&\int_{\Omega}\eta_{s_n}(|x-y|)^\frac{p-1}{p}\left|\frac{(x-y)\cdot h}{s_n}\right|^{p-1}\eta_{s_n}(|x-y|)^\frac{1}{p}\left|\frac{u_n(x)-u_n(y)}{s_n}\right|dy\\
    \leq &\lVert\eta_{s_n}(|x-y|)^\frac{p-1}{p}\left|\frac{(x-y)\cdot h}{s_n}\right|^{p-1}\rVert_{L^\frac{p}{p-1}}\\
    &\cdot\lVert\eta_{s_n}(|x-y|)^\frac{1}{p}\left|\frac{u_n(x)-u_n(y)}{s_n}\right|\rVert_{L^p}\\
    \leq & (\int_{\Omega}\eta_{s_n}(|x-y|)\left|\frac{(x-y)\cdot h}{s_n}\right|^{p}dy)^\frac{p-1}{p}\mathcal{E}_{s_n}(u_n)\\
    =&\sigma_{\eta}^{p-1}|h|^{p-1}\mathcal{E}_{s_n}(u_n) \quad \text{(when $n$ is large enough)}\\
    \leq & \sigma_{\eta}^{p-1}\lVert h\rVert_{L^\infty}^{p-1}\mathcal{E}_{s_n}(u_n)
\end{aligned}
\end{equation*}
        Therefore,
        \begin{equation}
            \liminf_{n\rightarrow\infty}\lVert|D_n^+(x)|+|D_n^-(x)|\rVert_{L^\infty}\leq\sigma_{\eta}^{p-1}\lVert h\rVert_{L^\infty}^{p-1}\liminf_{n\rightarrow\infty}\mathcal{E}_{s_n}(u_n)<\infty
            \label{eq:e6}
        \end{equation}
        which yields that, by the sequential Banach-Alaoglu theorem, there exists a subsequence $D_{n_k}^\pm$ and $D^\pm\in L^\infty(\Omega)$ such that
        \[D_{n_k}^+\rightharpoonup^* D^+,\ D_{n_k}^-\rightharpoonup^* D^-\]
        in $L^\infty(\Omega)$, i.e., for every $w\in L^1(\Omega)$, we have
        \begin{equation}
            \int_\Omega D_{n_k}^\pm wdx\longrightarrow\int_\Omega D^\pm wdx
            \label{eq:e7}
        \end{equation}
        For the next part, we only consider $D_{n_k}^+$ and $D^+$, and we denote $n_k$ by $n$ to simplify the notation. 
        Let $\phi\in\mathcal{C}_c^\infty(\Omega)\subset L^1(\Omega)$, then when $n$ is large enough so that $B(x,s_nr_\eta)\subset\Omega$ 
        for all $x\in \mathrm{supp}(\phi)$, we have
        \begin{equation*}
            \begin{aligned}
            \int_{\mathbb{R}^d}D_n^+\phi dx&=\int_{\mathbb{R}^d}\int_{\mathbb{R}^d}\eta_{s_n}(|x-y|)^\frac{p-1}{p}\left|\frac{(x-y)\cdot h}{s_n}\right|^{p-1}\eta_{h,n}^+(x-y)\frac{u_n(x)-u_n(y)}{s_n}dy\phi(x)dx\\
            (y=x+z)&=\int_{\mathbb{R}^d}\int_{\mathbb{R}^d}\eta_{s_n}(|z|)^\frac{p-1}{p}\left|\frac{z\cdot h}{s_n}\right|^{p-1}\eta_{h,n}^+(-z)\frac{u_n(x)-u_n(x+z)}{s_n}dz\phi(x)dx\\
            &=\int_{\mathbb{R}^d}\int_{\mathbb{R}^d}\eta_{s_n}(|z|)^\frac{p-1}{p}\left|\frac{z\cdot h}{s_n}\right|^{p-1}\eta_{h,n}^+(-z)\frac{\phi(x)-\phi(x-z)}{s_n}dzu_n(x)dx\\
            &=\int_{\mathbb{R}^d}\int_{\mathbb{R}^d}\eta(|z|)^\frac{p-1}{p}\left|z\cdot h\right|^{p-1}\eta_{h}^+(-z)\frac{\phi(x)-\phi(x-s_nz)}{s_n}dzu_n(x)dx
        \end{aligned}
        \end{equation*}
        So
        \begin{equation*}
            \begin{aligned}
            &|\int_{\mathbb{R}^d}D_n^+\phi dx-\int_{\mathbb{R}^d}\int_{\mathbb{R}^d}\eta(|z|)^\frac{p-1}{p}\left|z\cdot h\right|^{p-1}\eta_{h}^+(-z)z\cdot \nabla\phi(x) dzu(x)dx|\\
            &\leq\int_\Omega \{\lVert u_n-u\rVert_{L^\infty}|\int_{\mathbb{R}^d}\eta(|z|)^\frac{p-1}{p}\left|z\cdot h\right|^{p-1}\eta_{h}^+(-z)\frac{\phi(x)-\phi(x-s_nz)}{s_n}dz|\\
            &+\lVert u\rVert_{L^\infty}\int_{\mathbb{R}^d}\eta(|z|)^\frac{p-1}{p}\left|z\cdot h\right|^{p-1}\eta_{h}^+(-z)\left|\frac{\phi(x)-\phi(x-s_nz)-s_nz\cdot\nabla\phi(x)}{s_n}\right|dz\}dx\\
            &\leq\int_\Omega\{\lVert u_n-u\rVert_{L^\infty}\lVert\nabla\phi\rVert_{L^\infty}|\int_{\mathbb{R}^d}\eta(|z|)^\frac{p-1}{p}\left|z\cdot h\right|^{p-1}\eta_{h}^+(-z)|z|dz|\\
            &+    \frac{s_n}{2}\lVert\mathcal{H}(\phi)\rVert\lVert u\rVert_{L^\infty}\int_{\mathbb{R}^d}\eta(|z|)^\frac{p-1}{p}\left|z\cdot h\right|^{p-1}\eta_{h}^+(-z)|z|^2dz\}dx\\
            &\longrightarrow 0(n\longrightarrow\infty)
        \end{aligned}
        \end{equation*}
So using~\eqref{eq:e7}, we can prove that
\begin{equation}
        \begin{aligned}
            &\int_{\mathbb{R}^d}D^+\phi dx=\lim_{n\rightarrow\infty}\int_{\Omega}D_n^+\phi dx\\=&\int_{\mathbb{R}^d}\int_{\mathbb{R}^d}\eta(|z|)^\frac{p-1}{p}\left|z\cdot h\right|^{p-1}\eta_{h}^+(-z)z\cdot \nabla\phi(x) dzu(x)dx
            \label{E8}
        \end{aligned}
\end{equation}
For $D^-$, we have a similar conclusion.\\
        \indent Now we show $u\in W^{1,\infty}(\Omega)$ by~\eqref{E8}. For all $i\in\{1,2,\dots,d\}$, let
        \[h_x=\left\{\begin{array}{cc}
                         C_i(x)\mathrm{sgn}(\partial_i\phi(x))\nabla\phi(x) & |\nabla\phi(x)|\neq 0, \\
                         0                                         & |\nabla\phi(x)|=0
        \end{array}\right.\]
where
        \[C_i(x)\coloneqq(\frac{|\partial_i\phi(x)|}{|\nabla\phi(x)|^p})^\frac{1}{p-1}\]
        \[sgn(x)=\left\{\begin{array}{cc}
                         1 & x>0, \\
                         0 & x=0,\\
                         -1 & x<0
        \end{array}\right.\]
Since $|C_i(x)|\leq 1$, we have $\lVert h_x\rVert_{L^\infty}\leq\lVert \nabla\phi\rVert_{L^\infty}<\infty$ so that we can use~\eqref{E8} to get
\begin{equation*}
    \begin{aligned}
    \int_{\mathbb{R}^d}D^+\phi dx&=\int_{\mathbb{R}^d}\int_{\mathbb{R}^d}\eta(|z|)^\frac{p-1}{p}\left|z\cdot h\right|^{p-1}\eta_{h}^+(-z)z\cdot \nabla\phi(x) dzu(x)dx\\
    &=-\int_{\mathbb{R}^d}\int_{\mathbb{R}^d}\eta(|z|)^\frac{p-1}{p}\left|z\cdot \nabla\phi(x)\right|^{p}\eta_{h}^+(-z) dz\mathrm{sgn}(\partial_i\phi(x))C_i(x)^{p-1}u(x)dx\\
    &=-\frac{1}{2}\int_{\mathbb{R}^d}\int_{\mathbb{R}^d}\eta(|z|)\left|z\cdot \nabla\phi(x)\right|^{p} dz\mathrm{sgn}(\partial_i\phi(x))C_i(x)^{p-1}u(x)dx\\
    &=-\frac{\sigma_{\eta}^p}{2}\int_{\Omega}|\nabla\phi(x)|^p\mathrm{sgn}(\partial_i\phi(x))C_i(x)^{p-1}u(x)dx\\
    &=-\int_{\mathbb{R}^d}D^-\phi dx
\end{aligned}
\end{equation*}
and
\begin{equation}
\begin{aligned}
&\int_{\mathbb{R}^d}(D^--D^+)\phi dx\\
=&\sigma_{\eta}^p\int_{\Omega}|\nabla\phi(x)|^p\mathrm{sgn}(\partial_i\phi(x))C_i(x)^{p-1}u(x)dx\\
=&\sigma_{\eta}^p\int_{\Omega}\partial_i\phi(x)u(x)dx
\end{aligned}
\label{eq:e9}
\end{equation}
        Because $\frac{D^--D^+}{\sigma_{\eta}^p}\in L^\infty(\Omega)$ and $\phi$ can be any function in $\mathcal{C}_c^\infty(\Omega)$,
        we prove that $u\in W^{1,\infty}(\Omega)$.
        So now we can change~\eqref{E8} into
        \[\int_{\mathbb{R}^d}D^+\phi dx=-\int_{\mathbb{R}^d}\int_{\mathbb{R}^d}\eta(|z|)^\frac{p-1}{p}\left|z\cdot h\right|^{p-1}\eta_{h}^+(-z)z\cdot \nabla u(x) dz\phi(x)dx\]
        by the integration by parts with respect to $x$.
        Let $h_x=\frac{\nabla u(x)}{|\nabla u(x)|}$, then $\lVert h_x\rVert_{L^\infty}=1$ and we can get
        \begin{equation}
        \int_{\mathbb{R}^d}(D^+-D^-)\phi dx=\sigma_{\eta}^p\int_{\Omega}|\nabla u(x)|\phi(x)dx
        \label{eq:example}
        \end{equation}
        with the similar methods about getting~\eqref{eq:e9}.
        Because of the denseness of $\mathcal{C}_c^\infty(\Omega)$ in $L^1(\Omega)$, we actually have that
        \[\int_{\mathbb{R}^d}(D^+-D^-)w dx=\sigma_{\eta}^p\int_{\Omega}|\nabla u(x)|w(x)dx\]
        for any $w\in L^1(\Omega)$, which yields that
        \[\lVert D^+-D^-\rVert_{L^\infty}=\sigma_{\eta}^p\lVert \nabla u(x)\rVert_{L^\infty}\]
        Using~\eqref{eq:e6} again, now we have
\begin{equation*}
\begin{aligned}
        \liminf_{n\rightarrow\infty}\mathcal{E}_{s_n}(u_n)&\geq\frac{1}{\sigma_{\eta}^{p-1}}\liminf_{n\rightarrow\infty}\lVert|D_n^+|+|D_n^-|\rVert_{L^\infty}\\
        &\geq\frac{1}{\sigma_{\eta}^{p-1}}\liminf_{n\rightarrow\infty}\lVert D_n^+-D_n^-\rVert_{L^\infty} \\
        &\geq\frac{1}{\sigma_{\eta}^{p-1}}\lVert D^+-D^-\rVert_{L^\infty}\\
        &=\sigma_{\eta}\lVert \nabla u(x)\rVert_{L^\infty}\\
        &=\sigma_{\eta}\mathcal{E}(u)
\end{aligned}   
\end{equation*}
which completes the proof.
\end{proof}

    In the previous lemma we do not require additional properties on domain $\Omega$.
    However, the next lemma requires the following stronger boundary property on $\Omega$.
    For any $n\in\mathbb{N}^+$, define
    \[\tilde{\Omega}_n=\{x\in\Omega\big|B(x,\frac{1}{n})\subset\Omega\}\]
    where $B(x,r)$ represents the ball in $\mathbb{R}^d$ with center x and radius $r$.
    We need the following condition for $\Omega$:
    \begin{equation}
        d_n\coloneqq\sup_{x\in\Omega \backslash \tilde{\Omega}_n}\mathrm{dist}(x,\tilde{\Omega}_n)\longrightarrow 0(n\longrightarrow\infty)
        \label{eq:prop-boundary}
    \end{equation}
    This property says that the boundary of $\Omega$ can be approximated by the "smoothened interior".
    Being a bounded $C^2$ domain is sufficient for this property, see Proposition \ref{prop:a1} in the appendix.

    \begin{lemma}
    	\label{lem:2.3}
        Let $\Omega$ be a domain of class $C^1$ satisfying~\eqref{eq:prop-boundary} and kernel $\eta$ satisfy (K1)-(K3),
        then for all $u\in L^\infty(\Omega),\ \mathbb{R}^+\ni s_n\longrightarrow 0$, there exists a sequence $\{u_n\}_{n=1}^\infty\subset L^\infty(\Omega)$ such that
        \[\limsup_{n\rightarrow\infty}\mathcal{E}_{s_n}(u_n)\leq \sigma_{\eta}\mathcal{E}(u)\]
        and $u_n\longrightarrow u$ in $L^\infty(\Omega)$
    \end{lemma}

    \begin{proof}
        W.l.o.g, we can assume $u\in W^{1,\infty}(\Omega)$.
        Thanks to the Extension Theorem for functions in $W^{1,\infty}(\Omega)$ that introduced as Theorem 9.7 in~\cite{ref3} and Theorem 5.4.1 in~\cite{Evans}, we can extend $u$ to $\hat{u}\in W^{1,\infty}(\mathbb{R}^d)$ such that
        \[\hat{u}|_{\Omega}=u,\ \lVert\hat{u}\rVert_{L^\infty}\leq C_{\Omega}\lVert u\rVert_{L^\infty}\]
        Now we choose a non-negative function $\phi\in \mathcal{C}_c^\infty(B(1))$ such that $\int_{\mathbb{R}^d}\phi dx=1$ and do convolution between $\phi_\varepsilon\coloneqq\frac{1}{\varepsilon^d}\phi(\frac{\cdot}{\varepsilon})$ and $\hat{u}$ to get
        \[\hat{u}_{\varepsilon}(x)=\phi_\varepsilon*\hat{u}(x)=\int_{\mathbb{R}^d}\phi_\varepsilon(y)\hat{u}(x-y)dy\]
        Because of the property of convolution, we know $\hat{u}_{\varepsilon}$ is also smooth and for any $x\in\mathbb{R}^d$,
\begin{equation*}
    \begin{aligned}
        |\hat{u}_{\varepsilon}(x)-\hat{u}(x)|&\leq\int_{\mathbb{R}^d}\phi_\varepsilon(y)|\hat{u}(x-y)-\hat{u}(x)|dy\\
        &=\int_{\mathbb{R}^d}\phi(y)|\hat{u}(x-\varepsilon y)-\hat{u}(x)|dy\\
        &\leq \varepsilon\lVert\nabla\hat{u}\rVert_{L^\infty}\int_{\mathbb{R}^d}\phi(y)|y|dy
    \end{aligned}
\end{equation*}
which yields
        \[\lVert\hat{u}_\varepsilon-\hat{u}\rVert_{L^\infty}\longrightarrow0(\varepsilon\longrightarrow0)\]
        or more precisely,
\begin{equation}
    \lVert\hat{u}_\varepsilon-\hat{u}\rVert_{L^\infty}=O(\varepsilon)
    \label{eq:e12}
\end{equation}
Moreover, for any $i,j\in\{1,2,\dots,d\}$, we have
\begin{equation*}
    \begin{aligned}
        \frac{\partial^2}{\partial x_i\partial x_j}\hat{u}_{\varepsilon}(x)&\leq\int_{\mathbb{R}^d}|\frac{\partial^2}{\partial x_i\partial x_j}\phi_\varepsilon(x-y)||\hat{u}(y)|dy\\
        &\leq \lVert\hat{u}\rVert_{L^\infty}\int_{\mathbb{R}^d}|\frac{\partial^2}{\partial x_i\partial x_j}\phi_\varepsilon(x-y)|dy\\
        &=\lVert\hat{u}\rVert_{L^\infty}\int_{\mathbb{R}^d}|\frac{\partial^2}{\partial y_i\partial y_j}\phi(y)|dy\\
        &\leq C_{\hat{u},\phi}
    \end{aligned}
\end{equation*}
It means that there exists a constant $C$ which only with respect to $u$ and $\phi$ such that
\[\lVert\mathcal{H}(\hat{u}_\varepsilon)(x)\rVert_2\leq C,\ \forall x\in\mathbb{R}^d,\varepsilon>0\]
where $\mathcal{H}(\hat{u}_\varepsilon)$ represents the Hessian matrix of $\hat{u}_\varepsilon$ and $\lVert A\rVert_2\coloneqq\sup_{x\in\mathbb{R}^d}\frac{\lVert Ax\rVert_2}{\lVert x\rVert_2}$.

Now we choose $u_n=\hat{u}_{\frac{1}{n}}|_\Omega$, then
\[\lVert u_n-u\rVert_{L^\infty(\Omega)}\leq\lVert \hat{u}_{\frac{1}{n}}-\hat{u}\rVert_{L^\infty(\mathbb{R}^d)}\longrightarrow0(n\longrightarrow\infty)\]
and
\begin{equation*}
\begin{aligned}
    \mathcal{E}_{s_n}(u_n)&=\frac{1}{s_n}\esssup_{x\in\Omega}(\int_{\Omega}\eta_{s_n}(|x-y|)|u_n(x)-u_n(y)|^pdy)^\frac{1}{p}\\
    &=\frac{1}{s_n}\esssup_{x\in\Omega}(\int_{\Omega}\eta_{s_n}(|x-y|)|\hat{u}_{\frac{1}{n}}(x)-\hat{u}_{\frac{1}{n}}(y)|^pdy)^\frac{1}{p}\\
    &=\esssup_{x\in\Omega}(\int_{B(r_\eta),x-s_nz\in\Omega}\eta(|z|)\left|\frac{\hat{u}_{\frac{1}{n}}(x)-\hat{u}_{\frac{1}{n}}(x-s_nz)}{s_n}\right|^pdz)^\frac{1}{p}
\end{aligned}
\end{equation*}
With Taylor expansion and $\{\hat{u}_{\frac{1}{n}}\},\  \{\nabla\hat{u}_{\frac{1}{n}}\}$ are uniformly bounded, we have
\begin{equation*}
        \begin{aligned}
            &\left|\left|\frac{\hat{u}_{\frac{1}{n}}(x)-\hat{u}_{\frac{1}{n}}(x-s_nz)}{s_n}\right|^p-\left|\nabla\hat{u}_{\frac{1}{n}}(x)\cdot z\right|^p\right|\\
            \leq & C_p\left|\frac{\hat{u}_{\frac{1}{n}}(x)-\hat{u}_{\frac{1}{n}}(x-s_nz)}{s_n}-\nabla\hat{u}_{\frac{1}{n}}(x)\cdot z\right|\\
            \leq &\frac{C_p}{2}\lVert\mathcal{H}(\nabla\hat{u}_{\frac{1}{n}})\rVert_2s_n|z|^2\\
            \leq &\frac{CC_p}{2}s_n|z|^2
        \end{aligned}   
\end{equation*}
So
\begin{equation*}
        \begin{aligned}
            \mathcal{E}_{s_n}(u_n)&\leq(\esssup_{x\in\Omega}\int_{B(r_\eta),x-s_nz\in\Omega}\eta(|z|)\left|\nabla\hat{u}_{\frac{1}{n}}(x)\cdot z\right|^pdz+o(s_n))^\frac{1}{p}\\
            &=(\sigma_{\eta}^p\lVert\nabla\hat{u}_{\frac{1}{n}}\rVert_{L^\infty(\Omega)}^p+O(s_n))^\frac{1}{p}
        \end{aligned}   
\end{equation*}

Note that
       \begin{equation*}  \begin{aligned}
            \frac{\partial}{\partial x_i}\hat{u}_{\frac{1}{n}}(x)&=\int_{\mathbb{R}^d}\frac{\partial}{\partial x_i}\phi_\varepsilon(x-y)\hat{u}(y)dy\\
            &=\int_{\mathbb{R}^d}\phi_\varepsilon(x-y)\frac{\partial}{\partial y_i}\hat{u}(y)dy\\
            &=\int_{\mathbb{R}^d}\phi(y)\frac{\partial}{\partial y_i}\hat{u}(x-\varepsilon y)dy
       \end{aligned} \end{equation*}  
        which yields
       \begin{equation*}  \begin{aligned}
            |\nabla\hat{u}_{\frac{1}{n}}(x)|&=|\int_{\mathbb{R}^d}\phi(y)\nabla\hat{u}(x-\varepsilon y)dy|\\
            &\leq \int_{\mathbb{R}^d}\phi(y)|\nabla\hat{u}(x-\varepsilon y)|dy\\
            &\leq \lVert\nabla u\rVert_{L^\infty(\Omega)},\ \forall x\in\tilde{\Omega}_n
       \end{aligned} \end{equation*}  
        and for any $x\in\Omega/\tilde{\Omega}_n$, using Taylor expansion for $\nabla\hat{u}_{\frac{1}{n}}$, we have
        \[|\nabla\hat{u}_{\frac{1}{n}}(x)|\leq \lVert\nabla \hat{u}_{\frac{1}{n}}\rVert_{L^\infty(\tilde{\Omega}_n)}+\lVert\mathcal{H}(\nabla\hat{u}_{\frac{1}{n}})(x)\rVert_2d_n\leq\lVert\nabla\hat{u}_{\frac{1}{n}}\rVert_{L^\infty(\tilde{\Omega}_n)}+Cd_n\]
        Combine two parts of estimation, we get
        \[\lVert\nabla\hat{u}_{\frac{1}{n}}\rVert_{L^\infty(\Omega)}\leq\lVert\nabla u\rVert_{L^\infty(\Omega)}+Cd_n\]
        Thanks to~\eqref{eq:prop-boundary} and $s_n\rightarrow 0$,
        \[\mathcal{E}_{s_n}(u_n)\leq\sigma_{\eta}\lVert\nabla u\rVert_{L^\infty(\Omega)}+o(1)=\sigma_{\eta}\mathcal{E}(u)+o(1)\]
        from which we can get
        \[\limsup_{n\rightarrow\infty}\mathcal{E}_{s_n}(u_n)\leq \sigma_{\eta}\mathcal{E}(u)\]
    \end{proof}
    The previous two lemmas directly imply the $\mathrm{\Gamma}$-convergence of the respective functionals, which we
    can restate as a theorem below.

    \begin{theorem}[Non-Local to Local Convergence]
        Let $\Omega$ be a $C^1$ domain in $\mathbb{R}^d$ satisfying \eqref{eq:prop-boundary} (or $\Omega$ be a bounded $C^2$ domain)
        and kernel $\eta$ satisfy (K1)-(K3), then for any null sequence $\{s_n\}_{n=1}^\infty\subset\mathbb{R}^+$, we have
        \[\mathcal{E}_{s_n}\stackrel{\Gamma}{\longrightarrow}\sigma_{\eta}\mathcal{E}\ (n\rightarrow\infty)\]
    \end{theorem}

    \subsubsection{Discrete to Continuum Convergence}\label{subsec:discrete-to-continuum-convergence}
    In this part, we will prove \autoref{thm:1}. Before our proof, we need the following transform of $E_{n,\mathrm{cons}}(u)$ when $E_{n,\mathrm{cons}}(u)<\infty$.
   \begin{equation*}  \begin{aligned}
        E_{n,\mathrm{cons}}(u)&=E_{n}(\tilde{u})\\
        &=\frac{1}{s_n}\max_{x\in\Omega_n}(\frac{1}{|\Omega_n|}\sum_{y\in\Omega_n}\eta_{s_n}(|x-y|)|\tilde{u}(x)-\tilde{u}(y)|^p)^\frac{1}{p}\\
        &=\nu(\Omega)^{-\frac{1}{p}}\frac{1}{s_n}\max_{x\in\Omega_n}(\int_{\Omega}\eta_{s_n}(|x-T_n(y)|)|\tilde{u}(x)-\tilde{u}\circ T_n(y)|^pdy)^\frac{1}{p}\\
        &=\nu(\Omega)^{-\frac{1}{p}}\frac{1}{s_n}\nu_n\text{-}\esssup_{x\in\Omega}(\int_{\Omega}\eta_{s_n}(|x-T_n(y)|)|\tilde{u}(x)-\tilde{u}\circ T_n(y)|^pdy)^\frac{1}{p}\\
        &=\nu(\Omega)^{-\frac{1}{p}}\frac{1}{s_n}\nu\text{-}\esssup_{x\in\Omega}(\int_{\Omega}\eta_{s_n}(|T_n(x)-T_n(y)|)|\tilde{u}\circ T_n(x)-\tilde{u}\circ T_n(y)|^pdy)^\frac{1}{p}\\
        &=\nu(\Omega)^{-\frac{1}{p}}\frac{1}{s_n}\esssup_{x\in\Omega}(\int_{\Omega}\eta_{s_n}(|T_n(x)-T_n(y)|)|u(x)-u(y)|^pdy)^\frac{1}{p}
   \end{aligned} \end{equation*}  
    The last but two equality holds because of~\autoref{lem:esssup}.
    \begin{lemma}
        \label{lem:2.4}
        Let $\Omega\subset\mathbb{R}^d$ be a bounded $C^1$ domain and kernel $\eta$ satisfy (K1)-(K3), then
        with probability one,
        \[\liminf_{n\rightarrow\infty}E_{n,\mathrm{cons}}(u_n)\geq \sigma_{\eta}\nu(\Omega)^{-\frac{1}{p}}\mathcal{E}_{\mathrm{cons}}(u)\]
        for all $\mathbb{R}^+\ni s_n\longrightarrow 0,\ \frac{\delta_n}{s_n}\longrightarrow 0,\ L^\infty(\Omega)\ni u_n\longrightarrow u$ in $L^\infty(\Omega)$
    \end{lemma}
    \begin{proof}
        We assume $\liminf_{n\rightarrow\infty}E_{n,\mathrm{cons}}(u_n)<\infty,\ u_n=g$ on $\mathcal{O}$ w.l.o.g. Then for all $x\in\mathcal{O}$
\begin{equation*}
         \begin{aligned}
            |u(x)-g(x)|&\leq |u(x)-u_n(x)|+|u_n(x)-g(x)|\\
            &\leq \lVert u-u_n\rVert_{L^\infty}+0\\
            &\longrightarrow 0(n\longrightarrow\infty)
        \end{aligned}   
\end{equation*}
        So we prove $u=g$ on $\mathcal{O}$ which means $\mathcal{E}_{\mathrm{cons}}(u)=\mathcal{E}(u)$. We firstly prove our conclusion for a special kind of kernel, we assume there exists $t>0$ such that $\eta$ is a constant on $B(t)$. Note that
        \[|T_n(x)-T_n(y)|>s_nt\Longrightarrow|x-y|>s_nt-2\lVert T_n-Id\rVert_{L^\infty}\]
        Denote $s_n-\frac{2}{t}\lVert T_n-Id\rVert_{L^\infty}$ by $\tilde{s_n}$, then
       \begin{equation*}  \begin{aligned}
            \frac{|x-y|}{\tilde{s_n}}&\geq\frac{|T_n(x)-T_n(y)|-2\lVert T_n-Id\rVert_{L^\infty}}{\tilde{s_n}}\\
            &=\frac{|T_n(x)-T_n(y)|}{s_n}\frac{1-2\lVert T_n-Id\rVert_{L^\infty}/|T_n(x)-T_n(y)|}{1-2\lVert T_n-Id\rVert_{L^\infty}/(s_nt)}\\
            &\geq\frac{|T_n(x)-T_n(y)|}{s_n}
       \end{aligned} \end{equation*}  
        when $|T_n(x)-T_n(y)|>s_nt$. So combine with $\eta$ is a constant on $B(t)$, we have
        \[\eta_{s_n}(T_n(x)-T_n(y))\geq\frac{\tilde{s_n}^d}{s_n^d}\eta_{\tilde{s_n}}(|x-y|)\]
        Now using the transfrom at the beginning of this part,
       \begin{equation*}  \begin{aligned}
            E_{n,\mathrm{cons}}(u_n)&=\nu(\Omega)^{-\frac{1}{p}}\frac{1}{s_n}\esssup_{x\in\Omega}(\int_{\Omega}\eta_{s_n}(|T_n(x)-T_n(y)|)|u_n(x)-u_n(y)|^pdy)^\frac{1}{p}\\
            &\geq\nu(\Omega)^{-\frac{1}{p}}\frac{\tilde{s_n}^d}{s_n^{d+1}}\esssup_{x\in\Omega}(\int_{\Omega}\eta_{\tilde{s_n}}(|x-y|)|u_n(x)-u_n(y)|^pdy)^\frac{1}{p}\\
            &=\nu(\Omega)^{-\frac{1}{p}}\frac{\tilde{s_n}^{d+1}}{s_n^{d+1}}\mathcal{E}_{\tilde{s_n}}(u_n)
       \end{aligned} \end{equation*}  
        Note that we require $\lim_{n\rightarrow\infty}\frac{\delta_n}{s_n}=0$. With~\eqref{eq:e2}, we have $\lim_{n\rightarrow\infty}\tilde{s_n}=0$ and $\lim_{n\rightarrow\infty}\frac{\tilde{s_n}}{s_n}=1$. So
       \begin{equation*}  \begin{aligned}
            \liminf_{n\rightarrow\infty} E_{n,\mathrm{cons}}(u_n)&\geq\nu(\Omega)^{-\frac{1}{p}}\liminf_{n\rightarrow\infty}\frac{\tilde{s_n}^{d+1}}{s_n^{d+1}}\mathcal{E}_{\tilde{s_n}}(u_n)\\
            &=\nu(\Omega)^{-\frac{1}{p}}\liminf_{n\rightarrow\infty}\mathcal{E}_{\tilde{s_n}}(u_n)\\
            (using\ Lemma~\ref{lem:2.2})&\geq\sigma_{\eta}\nu(\Omega)^{-\frac{1}{p}}\mathcal{E}(u)\\
            &=\sigma_{\eta}\nu(\Omega)^{-\frac{1}{p}}\mathcal{E}_{\mathrm{cons}}(u)
       \end{aligned} \end{equation*}  
        We have proved the conclusion for a special kind of $\eta$, for the general situation, we define
        \[\eta_k(t)\coloneqq\left\{\begin{array}{cc}
                                       \eta(\frac{1}{k}) & \text{if }t\in[0,\frac{1}{k}] \\
                                       \eta(t)           & \text{if }t>\frac{1}{k}
        \end{array}
        \right.\quad \forall k\in\mathbb{N}^+\]
        Since $\eta$ is monotonously decreasing, we have $\eta_k\leq\eta$. So for any $k\in\mathbb{N}^+$
       \begin{equation*}  \begin{aligned}
            \liminf_{n\rightarrow\infty} E_{n,\mathrm{cons}}(u_n)&\coloneqq\liminf_{n\rightarrow\infty} E_{n,\mathrm{cons}}^\eta(u_n)\\
            &\geq \liminf_{n\rightarrow\infty} E_{n,\mathrm{cons}}^{\eta_k}(u_n)\\
            &\geq\sigma_{\eta_k}\nu(\Omega)^{-\frac{1}{p}}\mathcal{E}_{\mathrm{cons}}(u)
       \end{aligned} \end{equation*}  
        Because $\eta$ is continuous and non-negative, we have $0\leq\eta_k\uparrow\eta\ (k\longrightarrow\infty)$. So by monotone convergence theorem, let $k$ tends to infinity, we can get
        \[\liminf_{n\rightarrow\infty} E_{n,\mathrm{cons}}(u_n)\geq\sigma_{\eta}\nu(\Omega)^{-\frac{1}{p}}\mathcal{E}_{\mathrm{cons}}(u)\]
    \end{proof}

    \begin{lemma}
        \label{lem:2.5}
        Let $\Omega$ be a bounded $C^1$ domain satisfying~\eqref{eq:prop-boundary}, and let kernel $\eta$ satisfy (K1)-(K3).
        Then with probability one, for all $u\in L^\infty(\Omega),\ \mathbb{R}^+\ni s_n\longrightarrow 0,\ \frac{\delta_n}{s_n}\longrightarrow 0$,
        there exists a sequence $\{u_n\}_{n=1}^\infty\subset L^\infty(\Omega)$ such that
        \[\limsup_{n\rightarrow\infty}E_{n,\mathrm{cons}}(u_n)\leq \sigma_{\eta}\nu(\Omega)^{-\frac{1}{p}}\mathcal{E}_{\mathrm{cons}}(u)\]
        and $u_n\longrightarrow u$ in $L^\infty(\Omega)$.
    \end{lemma}
    \begin{proof}
        Assume $u\in W^{1,\infty}(\Omega)$ and $u=g$ on $\mathcal{O}$ w.l.o.g. Now let
        \[u_n=u|_{\Omega_n}\circ T_n=u\circ T_n\]
        then $E_{n,\mathrm{cons}}(u_n)=E_{n}(u_n)$ and $u_n=u\circ T_n=u=g$ on $\mathcal{O}$, so $E_{n}(u_n)<\infty$.\\
        \indent We firstly prove $u_n\longrightarrow u$ in $L^\infty(\Omega)$, actually using the extension $\hat{u}$ of $u$ defined in Lemma~\ref{lem:2.3},
        \[\lVert u_n-u\rVert_{L^\infty(\Omega)}=\lVert \hat{u}\circ T_n-\hat{u}\rVert_{L^\infty(\Omega)}\leq\lVert \nabla \hat{u}\rVert_{L^\infty}\lVert T_n-Id\rVert_{L^\infty}\]
        With~\eqref{eq:e2}, we have
        \[\lVert u_n-u\rVert_{L^\infty}=O(\delta_n)\longrightarrow 0(n\longrightarrow\infty)\]
        Be similar to Lemma~\ref{lem:2.4}, we firstly assume that there exists $t>0$ such that $\eta$ is a constant on $B(t)$. Then it is similar to prove,
        \[\eta_{s_n}(T_n(x)-T_n(y))\leq\frac{\tilde{s_n}^d}{s_n^d}\eta_{\tilde{s_n}}(|x-y|)\]
        for $\tilde{s_n}\coloneqq s_n+\frac{2}{t}\lVert T_n-Id\rVert_{L^\infty}$\\
        So
       \begin{equation*}  \begin{aligned}
            E_{n,\mathrm{cons}}(u_n)&=\nu(\Omega)^{-\frac{1}{p}}\frac{1}{s_n}\esssup_{x\in\Omega}(\int_{\Omega}\eta_{s_n}(|T_n(x)-T_n(y)|)|u_n(x)-u_n(y)|^pdy)^\frac{1}{p}\\
            &\leq\nu(\Omega)^{-\frac{1}{p}}\frac{\tilde{s_n}^d}{s_n^{d+1}}\esssup_{x\in\Omega}(\int_{\Omega}\eta_{\tilde{s_n}}(|x-y|)|u_n(x)-u_n(y)|^pdy)^\frac{1}{p}\\
            &=\nu(\Omega)^{-\frac{1}{p}}\frac{\tilde{s_n}^{d+1}}{s_n^{d+1}}\mathcal{E}_{\tilde{s_n}}(u_n)
       \end{aligned} \end{equation*}  
        which yields
        \[\limsup_{n\rightarrow\infty}E_{n,\mathrm{cons}}(u_n)\leq \nu(\Omega)^{-\frac{1}{p}}\limsup_{n\rightarrow\infty} \mathcal{E}_{\tilde{s_n}}(u_n)\]
        Now consider the difference between $\{u_n\}$ and the sequence $\{\hat{u}_{\delta_n}\}$ defined in Lemma~\ref{lem:2.3}, we have
       \begin{equation*}  \begin{aligned}
            \lVert u_n-\hat{u}_{\delta_n}\rVert_{L^\infty(\Omega)}&\leq \lVert u_n-u\rVert_{L^\infty(\Omega)}+\lVert\hat{u}_{\delta_n}-\hat{u}\rVert_{L^\infty(\mathbb{R}^d)}\\
            (thanks\ to\ \eqref{eq:e12})&=\lVert u_n-u\rVert_{L^\infty(\Omega)}+O(\delta_n)\\
            &= O(\delta_n)
       \end{aligned} \end{equation*}  
        So combine with~\eqref{eq:e11} and $\frac{\delta_n}{s_n}\longrightarrow 0$, we get
       \begin{equation*}  \begin{aligned}
            \limsup_{n\rightarrow\infty}E_{n,\mathrm{cons}}(u_n)\leq \nu(\Omega)^{-\frac{1}{p}}\limsup_{n\rightarrow\infty}\mathcal{E}_{\tilde{s_n}}(\hat{u}_{\delta_n})\leq\sigma_{\eta}\nu(\Omega)^{-\frac{1}{p}}\mathcal{E}(u)=\sigma_{\eta}\nu(\Omega)^{-\frac{1}{p}}\mathcal{E}_{\mathrm{cons}}(u)
       \end{aligned} \end{equation*}  
        where the second inequality holds because of Lemma~\ref{lem:2.3} (Note that we can change $\{\frac{1}{n}\}$
        in the proof into any positive sequence tending to 0).\\
        \indent For the general situation, it is also similar to Lemma~\ref{lem:2.4}. We define
        \[\eta_k(t)\coloneqq\left\{\begin{array}{cc}
                                       \eta(0) & \text{if }t\in[0,\frac{1}{k}] \\
                                       \eta(t) & \text{if }t>\frac{1}{k}
        \end{array}
        \right.\quad \forall k\in\mathbb{N}^+\]
        We have $\eta_k\geq\eta$.
        So for any $k\in\mathbb{N}^+$
       \begin{equation*}  \begin{aligned}
            \limsup_{n\rightarrow\infty} E_{n,\mathrm{cons}}(u_n)&\coloneqq \limsup_{n\rightarrow\infty} E_{n,\mathrm{cons}}^\eta(u_n)\\
            &\leq \limsup_{n\rightarrow\infty} E_{n,\mathrm{cons}}^{\eta_k}(u_n)\\
            &\leq\sigma_{\eta_k}\nu(\Omega)^{-\frac{1}{p}}\mathcal{E}_{\mathrm{cons}}(u)
       \end{aligned} \end{equation*}  
        Let $k$ tend to infinity, we finally get our conclusion.
    \end{proof}
    Combining Lemma~\ref{lem:2.4} and Lemma~\ref{lem:2.5} we immediately obtain the $\mathrm{\Gamma}$-convergence of
    the discrete functionals to those defined in the continuum, which is the statement of
    \autoref{thm:1}.

    \subsection{Compactness}\label{sec:compactness}

    In this section we prove the following compactness result:

    \begin{theorem}[Compactness]
        \label{thm:comp}
        Under the same assumption in \autoref{thm:2}, with probability one
        the functionals $ \{E_{n,\mathrm{cons}}\}_{n\in \mathbb{N}}$
        satisfies the compactness property, that is, any bounded sequence $ \{u_n\}_{n\in \mathbb{N}} \subseteq L^\infty(\Omega)$
        for which
       \begin{equation*}  \begin{aligned}
            \sup_{k \in \mathbb{N}} E_{n,\mathrm{cons}}(u_n) < +\infty
       \end{aligned} \end{equation*}  
        is relatively compact.
    \end{theorem}

    As soon as we prove \autoref{thm:comp}, with \autoref{lem:minimizer}, we
    can easily get \autoref{thm:2} as a corollary. In order to prove the compactness property, we need the following lemma, which is Lemma 4.1 in~\cite{ref1} using classical ideas from Lem.IV.5.4 in~\cite{dunford1988linear}.

    \begin{lemma}
        Let $\Omega$ be a set with finite measure and $K \subseteq L^\infty(\Omega)$ be a bounded set w.r.t.
        $\|\cdot\|_{L^\infty}$ such that for every $\varepsilon > 0$ there exists a finite partition $ \{V_i\}_{i=1}^n$
        of $\Omega$ with
        \begin{align}
            \esssup_{x,y\in V_i} |u(x) - u(y)| \leq \varepsilon, \; \forall u \in K, i=1,\dots,n,
        \end{align}
        then $K$ is relatively compact.
    \end{lemma}

    With the previous lemma, we now begin to prove a compactness result for the non-local functionals defined in~\eqref{eq:e5}.


    \begin{lemma}
        Let $\Omega\subset\mathbb{R}^d$ be a bounded $C^2$ domain satisfying~\eqref{eq:prop-boundary},
        let the kernel $\eta$ fulfil (K1)-(K3) and $\{s_n\}_{n\in\mathbb{N}}\subset\mathbb{R}^+$ a null sequence.
        Then for any bounded sequence $ \{u_n\}_{n\in \mathbb{N}} \subseteq L^\infty(\Omega)$
        for which
       \begin{equation*}  \begin{aligned}
            \sup_{n \in \mathbb{N}} \mathcal{E}_{s_n}(u_n) < +\infty
       \end{aligned} \end{equation*}  
        is relatively compact.
        \label{lem:nonlocal_comp}
    \end{lemma}
    \begin{proof}
        W.l.o.g. we rescale the kernel such that
       \begin{equation*}  \begin{aligned}
            \eta(t) \geq 1,\; \text{for} \; t \leq 1.
       \end{aligned} \end{equation*}  

        In the following proof we denote $C,C_1,C_2,\dots$ positive constants that is only
        dependent of $d,p$ and $\Omega$.

        We claim that there exists a constant $r_0 > 0$ such that for
        all $0 < r < r_0$,
        \begin{align}
            \nu(B(x,r)\cap \Omega) \geq C r^d,\; \forall x \in \Omega
            \label{eq:vol}
        \end{align}
        here $B(x,r) \subset \Omega$ is an open ball centered at $x$ with radius $r$.
        This property says that the Lebesgue density of $\Omega$ is uniformly bounded below, which
        is also a boundary property of $\Omega$.
        A proof of this property is given as Proposition \ref{prop:a2} in the appendix.

        Since $\{s_n\}_{n\in\mathbb{N}}\subset\mathbb{R}^+$ is a null sequence,
        w.l.o.g. we can assume $s_n < r_0,\; \forall n\in \mathbb{N}$.

        The key idea is the following estimation: Let $u \in L^{\infty}(\Omega)$,
        choose $0 < r < r_0$.
        For almost all $x \in \Omega$, let $U = B(x,r) \cap \Omega$,
       \begin{equation*}  
       \begin{aligned}
            & \frac{1}{\nu(U)} \int_{U} |u(x) - u(y)| dy \\
            \leq & \left( \frac{1}{\nu(U)} \int_{U} |u(x) - u(y)|^p dy  \right)^{\frac{1}{p}} \quad \text{(H\"older inequality)} \\
            \leq&   \left( \frac{1}{Cr^d} \int_{U}|u(x) - u(y)|^p dy  \right)^{\frac{1}{p}}  
            \quad \text{(by condition~\eqref{eq:vol})} \\
            \leq&   \left( \frac{1}{C} \int_{U} \eta_r(|x-y|) |u(x) - u(y)|^p dy  \right)^{\frac{1}{p}} \quad \text{(kernel rescaling)} \\
            \leq&  C_1 \left( \int_{\Omega} \eta_r(|x-y|) |u(x) - u(y)|^p dy  \right)^{\frac{1}{p}} \\
            \leq & C_1 r \mathcal{E}_r(u)
       \end{aligned} 
       \end{equation*}  
        or equivalently
        \begin{align}
            \esssup_{x\in \Omega}\int_{B(x,r)\cap \Omega} |u(x) - u(y)| dy \leq C_2 r^{d+1}  \mathcal{E}_r(u)
            \label{eq:int_bound}
        \end{align}
        For almost all $x,z$ such that $|x-z| < \frac{r}{2}$, let $U = B(x,r) \cap \Omega$,
        then $V = B(z,\frac{r}{2}) \cap \Omega \subset U$, so that
       \begin{equation*}  \begin{aligned}
            |u(x) - u(z)| &=
            \frac{1}{\nu(V)}|\int_{V}u(x)dy - \int_V u(z) dy| \\
            &\leq  \frac{1}{Cr^d} \left(\int_V |u(x)-u(y)| dy + \int_V |u(y)-u(z)| dy\right) \\
            &\leq \frac{1}{Cr^d} \left(\int_U |u(x)-u(y)| dy + \int_V |u(y)-u(z)| dy\right) \\
            &\leq C_3 r  \mathcal{E}_r(u) \quad (\text{by\eqref{eq:int_bound}}) \\
       \end{aligned} \end{equation*}  
        We restate the obtained result as follows.
        Assume $u\in L^{\infty}(\Omega),0 < r < r_0$ and $\mathcal{E}_r(u) < \infty$, then
        the following inequality holds:
        \begin{align}
            \esssup_{z_1,z_2\in\Omega, |z_1-z_2| < \frac{r}{2}} |u(z_1) - u(z_2)| \leq C_3 r\mathcal{E}_r(u).
            \label{eq:bound}
        \end{align}

        Now we follow the same route as the proof of Lemma 4.3 in~\cite{ref1}, applying the previous lemma.
        Let $\varepsilon > 0$ be given.
        With~\eqref{eq:convex} in Proposition \ref{prop:a3},
        there exists $\delta > 0$ such that for every $x,y\in \Omega$ with $|x-y| \leq \delta$,
        there is a path $\gamma : [0,1]\rightarrow \Omega$ connecting $x,y$, and
       \begin{equation*}  \begin{aligned}
            \mathrm{len}(\gamma) \leq (1+\varepsilon)|x-y|.
       \end{aligned} \end{equation*}  
        For every $n$, we can divide the path by points $0=t_0 < \dots < t_{k_n+1} = 1$ such that
        for $z_i \coloneqq \gamma(t_i)$ we have that
       \begin{equation*}  \begin{aligned}
            |z_i - z_{i+1}| \leq \frac{s_n}{2}, \; i = 0,\dots,k_n
       \end{aligned} \end{equation*}  
        where
       \begin{equation*}  \begin{aligned}
            k_n \leq \lfloor (1+\varepsilon)\frac{|x-y|}{s_n} \rfloor
       \end{aligned} \end{equation*}  
        Then we have that for almost all $x,y$ with $|x-y| \leq \delta$,
       \begin{equation*}  \begin{aligned}
            |u_n(x) - u_n(y)| & \leq \sum_{i=0}^{k_n}|u_n(z_i) - u_n(z_{i+1})| \\
            &\leq \sum_{i=0}^{k_n} C_3 s_n \mathcal{E}_{s_n}(u_n) \quad \text{(by\eqref{eq:bound})}\\
            & \leq C_3 s_n k_n \sup_{n\in \mathbb{N}} \mathcal{E}_{s_n}(u_n) \\
            & \leq C_4(1+\varepsilon) |x-y|.
       \end{aligned} \end{equation*}  

        Choosing a partition $ \{V_i\}_{i=1}^N$ of $\Omega$ with sufficiently small diameter such that
       \begin{equation*}  \begin{aligned}
            \mathrm{diam}(V_i) < \min \left( \delta, \frac{\varepsilon}{C_4(1+\varepsilon)} \right)
       \end{aligned} \end{equation*}  
        then it is satisfied that
       \begin{equation*}  \begin{aligned}
            \esssup_{x,y\in V_i} |u(x) - u(y)| \leq \varepsilon, \; \forall u \in K, i=1,\dots,n,
       \end{aligned} \end{equation*}  
        By applying the previous lemma we prove the compactness property.
    \end{proof}

    We can prove the following compactness result for discrete functionals:

    \begin{lemma}
        Let $\Omega\subset\mathbb{R}^d$ be a bounded $C^2$ domain,
        let the kernel $\eta$ fulfil (K1)-(K3) and $\{s_n\}_{n\in\mathbb{N}}\subset\mathbb{R}^+$ a null sequence satisfying
        $\frac{\delta_n}{s_n}\longrightarrow 0$.
        Then with probability one, for any bounded sequence $ \{u_n\}_{n\in \mathbb{N}} \subseteq X$
        for which
       \begin{equation*}  \begin{aligned}
            \sup_{n \in \mathbb{N}} E_{n,\mathrm{cons}}(u_n) < +\infty
       \end{aligned} \end{equation*}  
        is relatively compact.
        \label{lem:compactness}
    \end{lemma}

    \begin{proof}
        With the same argument in the proof of Lemma~\ref{lem:2.4}, we have that there exist a sequence $\{\tilde{s}_n\}_{n\in\mathbb{N}}\subset\mathbb{R}^+$
        such that
       \begin{equation*}  \begin{aligned}
            E_{n,\mathrm{cons}}(u_n) \geq \nu(\Omega)^{-\frac{1}{p}}\left(\frac{\tilde{s}_n}{s_n}\right)^{d+1} \mathcal{E}_{\tilde{s}_n}(u_n),
       \end{aligned} \end{equation*}  
        and
       \begin{equation*}  \begin{aligned}
            \lim_{n\rightarrow\infty} \frac{\tilde{s}_n}{s_n} = 1.
       \end{aligned} \end{equation*}  
        Therefore,
       \begin{equation*}  \begin{aligned}
            \sup_n \mathcal{E}_{\tilde{s}_n}(u_n) \leq
            \nu(\Omega)^{\frac{1}{p}}\sup_n \left(\frac{s_n}{\tilde{s}_n}\right)^{d+1} \sup_n E_{n,\mathrm{cons}}(u_n) < +\infty,
       \end{aligned} \end{equation*}  
        and Lemma~\ref{lem:nonlocal_comp} guarantees the relative compactness of  $\{u_n\}_{n\in \mathbb{N}}$.
    \end{proof}

	\section{Theoretical Analysis for Manifolds}
	\label{sec:manifold}
	In this part, we will generalize our conclusions to manifolds.
	\paragraph{Transportation map $T_n$} In the case of domains, transportation map $T_n$ is an important bridge for us to establish the connection between functionals on $L^\infty(\Omega_n)$ and functionals on $L^\infty(\Omega)$. So for the manifolds, we also hope to define a family of such transportation maps. However, as shown below, $T_n$ we define for manifolds is no longer a transportation map between measures in general, but it can still be considered as transportation map from a manifold to its finite random sample. Now consider $\mathcal{M}$ as a $k$-dimensional manifold in $\mathbb{R}^d$ which is smooth and bounded. Denote by $\nu$ the $k$-dimensional Hausdorff measure restricted on $\mathcal{M}$. We define $\mathcal{O},\ \nu_n$ and $\mathcal{M}_n$ using the same way as $\mathcal{O},\ \nu$ and $\Omega_n$ for a domain $\Omega$ in \autoref{sec:main}.\\
	\indent With Heine-Borel theorem, we can get a finite open covering $\mathcal{M}=\bigcup_{i=1}^m\mathcal{M}^i$ such that each $\mathcal{M}^i$ has a smooth global coordinate representation  $\phi^i:D^i\rightarrow \mathcal{M}^i$ where $D^i$ is a bounded $C^2$ domain in $\mathbb{R}^k$(such as the unit ball with center at origin). We now divide $\mathcal{M}$ into disjoint parts:
	\[\mathcal{M}=\bigcup_{i=1}^{m}\tilde{\mathcal{M}}^i\cup\varGamma\]
	where $\tilde{\mathcal{M}}^i$ is the set of interior points (with respect to $\mathcal{M}$) of $\mathcal{M}^i\backslash\bigcup_{j=1}^{i-1}\mathcal{M}^j$ and $\varGamma:=\mathcal{M}\backslash\bigcup_{i=1}^{m}\tilde{\mathcal{M}}^i$ is a $\nu-$null set. We define $\nu^i:=\frac{\nu|_{\tilde{\mathcal{M}}^i}}{\tilde{\mathcal{M}}^i}$ as the probability measure induced by $\nu|_{\tilde{\mathcal{M}}^i}$ and $\nu_n^i$ as the empirical measure of $\mathcal{M}_n\cap\tilde{\mathcal{M}}^i$. Since $D^i$ is a bounded domain in $\mathbb{R}^k$, we have known ,with probabilyty one, the existence of the transportation map from $\nu^i\circ\phi^i$ to $\nu_n^i\circ\phi^i$ satisfying \eqref{eq:e2}. We denote such transportation map by $\tilde{T}_{D^i,n}^i$. Then $\tilde{T}_n^i:=\phi^i\circ\tilde{T}_{D^i,n}^i\circ(\phi^i)^{-1}$ is the transportation map from $\nu^i$ to $\nu_n^i$ satisfying \eqref{eq:e2} since we may let $\phi^i$ and $(\phi^i)^{-1}$ are both Lipschitz continuous. For $\mathcal{M}^i$, we define $T_n^i$ similarly.\\
	\indent With the above preparation, we can define the transportation map as follows:   
	\begin{equation*}
		\begin{array}{cccc}
		T_n:&\mathcal{M}&\longrightarrow&\mathcal{M}_n\\
		&y&\longmapsto&\left\{\begin{array}{cc}
		\tilde{T}_n^i(y)&\text{if }y\in\tilde{\mathcal{M}}^i\text{ for some }i,\\
		y&\text{if }y\in\Gamma.
		\end{array}\right.
		\end{array}
	\end{equation*}
	With this definition, it is obvious that $T_n$ satisfies \eqref{eq:e2} since each $\tilde{T}_n^i$ has such property. For any measurable set $A\subset\mathcal{M}$,  we have the following equation:
	\begin{equation*}
		\begin{aligned}
			\frac{1}{\nu(\mathcal{M})}\nu\circ T_n^{-1}(A)&=\frac{1}{\nu(\mathcal{M})}\sum_{i=1}^m\nu\circ T_n^{-1}(A\cap\tilde{\mathcal{M}}^i)\\
			&=\frac{1}{\nu(\mathcal{M})}\sum_{i=1}^m\nu\circ(\tilde{T}_n^i)^{-1}(A\cap\tilde{\mathcal{M}}^i)\\
			&=\frac{1}{\nu(\mathcal{M})}\sum_{i=1}^m\nu|_{\tilde{\mathcal{M}}^i}\circ(\tilde{T}_n^i)^{-1}(A)\\
			&=\frac{1}{\nu(\mathcal{M})}\sum_{i=1}^m\nu(\tilde{\mathcal{M}}^i)\nu_n^i(A)\\
		\end{aligned}
	\end{equation*}
	From the last step, we know that $T_n$ may not be a transportation map from $\frac{\nu}{\nu(\mathcal{M})}$ to $\nu_n$ unless $\mathcal{M}$ has a global coordinate representation so that $m$ can be $1$. But it is still enough for us to get the convergence results for manifolds. 
	\paragraph{Main results}Similarly to the situation of domain, with the transportation map defined as above, we can define the continuum functionals and discrete ones we need as follows:
	\begin{equation}
	\label{eq:discrete-funcm}
	E_n(u)\coloneqq
	\begin{cases}
	\frac{1}{s_n}\mathop{\max}\limits_{x\in\mathcal{M}_n}(\frac{1}{|\mathcal{M}_n|}\sum_{y\in\mathcal{M}_n}\eta_{s_n}(|x-y|)|u(x)-u(y)|^p)^\frac{1}{p}
	& \text{if } u=g \text{ on } \mathcal{O},\\
	\infty
	& \text{otherwise.}
	\end{cases}
	\end{equation}
	\begin{equation}
	E_{n,\mathrm{cons}}(u)\coloneqq
	\begin{cases}
	E_{n}(\tilde{u}) & \text{if } u=\tilde{u}\circ T_n \text{ for some } \tilde{u}:\mathcal{M}_n\longrightarrow\mathbb{R}, \\
	\infty           & \text{otherwise}.
	\end{cases}
	\label{eq:e3m}
	\end{equation}
	\begin{equation*}
	\label{eq:continuum-funcm}
	\mathcal{E}(u)\coloneqq
	\begin{cases}
	\esssup_{x\in\mathcal{M}}|\nabla_{\mathcal{M}} u| & \text{if } u\in W^{1,\infty}(\mathcal{M}), \\
	\infty                         & \text{otherwise}.
	\end{cases}
	\end{equation*}
	\begin{equation}
	\mathcal{E}_{\mathrm{cons}}(u)\coloneqq\left\{\begin{array}{ll}
	\mathcal{E}(u) \ \  & \text{if } u=g\text{ on } \mathcal{O}, \\
	\infty         & \text{otherwise}.
	\end{array}\right.
	\label{eq:e4m}
	\end{equation}
	In these definitions, $\nabla_{\mathcal{M}} u(x)$ can be defined as $(J_\phi(\phi^{-1}(x))^\dagger)^T\nabla(u\circ\phi)(\phi^{-1}(x))$ for any $x\in\mathcal{M}$ where $\phi$ is a local coordinate map of the neighborhood of $x$, $J_\phi$ means the Jacobi matrix and $\dagger$ means the Moore-Penrose generalized inverse for matrices. $u\in W^{1,\infty}(\mathcal{M})$ means that $u\circ\phi$ is a $W^{1,\infty}$ function on the coordinate space for any local coordinate map $\phi$. And $\sigma_\eta$ used in the following theorems define as $(\int_{\mathbb{R}^k}\eta(|z|)|z\cdot e_1|^p dz)^\frac{1}{p}$ which is different from the domain situation by changing dimension $d$ into the dimension of the manifold $k$. $\eta_{s_n}$ also define as $\frac{1}{s_n^k}\eta(\frac{\cdot}{s_n})$ instead of $\frac{1}{s_n^d}\eta(\frac{\cdot}{s_n})$\\
	\indent Our main results are the $\mathrm{\Gamma}$-convergence for these functionals and convergence of their minimizers stated as the following two theorems:  
	
	\begin{theorem}
		[Discrete to continuum $\mathrm{\Gamma}$-convergence]
		\label{thm:1m}
		Let $\mathcal{M}\subset\mathbb{R}^d$ be a bounded, smooth $k$-dimensional manifold,
		and let the Lipschitz continuous kernel $\eta$ fulfil (K1)-(K3).
		Then for any null sequence $\{s_n\}_{n\in\mathbb{N}}\subset\mathbb{R}^+$ which satisfies $\lim_{n\rightarrow\infty}\frac{\delta_n}{s_n}=0$,
		with probability one we have
		\begin{enumerate}
			\item[(1)] for any sequence $\{u_n\}_{n\in\mathbb{N}}\subset L^\infty(\mathcal{M})$ converging to $u\in L^\infty(\mathcal{M})$ in $L^\infty(\mathcal{M})$,
			\[\liminf_{n\rightarrow\infty}E_{n,\mathrm{cons}}(u_n)\geq\sigma_\eta\nu(\mathcal{M})^{-\frac{1}{p}} \mathcal{E}_{\mathrm{cons}}(u)\]
			\item[(2)] for any $u\in L^\infty(\mathcal{M})$, there exists a sequence $\{u_n\}_{n\in\mathbb{N}}\subset L^\infty(\mathcal{M})$ converging to $u$ in $L^\infty(\mathcal{M})$ such that
			\[\limsup_{n\rightarrow\infty}E_{n,\mathrm{cons}}(u_n)\leq\sigma_\eta\nu(\mathcal{M})^{-\frac{1}{p}} \mathcal{E}_{\mathrm{cons}}(u)\]
		\end{enumerate}
	\end{theorem}
	
	\begin{theorem}[Convergence of minimizers]
		\label{thm:2m}
		Let manifold $\mathcal{M}$, kernel $\eta$ and null sequence\\ $\{s_n\}_{n\in\mathbb{N}}\subset\mathbb{R}^+$
		satisfies the assumptions in \autoref{thm:1m}.
		Then with probability one, any bounded sequence $(u_n)_{n\in \mathbb{N}}$ such that
		\begin{equation*}
		\lim_{n\rightarrow \infty} \left( E_{n,\mathrm{cons}}(u_n) - \inf_{u\in L^\infty(\mathcal{M})} E_{n,\mathrm{cons}}(u)
		\right) = 0
		\end{equation*}
		is relatively compact in $ L^\infty(\mathcal{M})$ and
		\begin{equation*}
		\lim_{n\rightarrow \infty} E_{n,\mathrm{cons}}(u_n) = \sigma_{\eta}\nu(\mathcal{M})^{-\frac{1}{p}}\min_{u\in L^\infty(\mathcal{M})} \mathcal{E}_{\mathrm{cons}}(u).
		\end{equation*}
		Furthermore, every cluster point of $(u_n)_{n\in \mathbb{N}}$ is a minimizer of $ \mathcal{E}_{\mathrm{cons}}$.
	\end{theorem}
	The method of proving \autoref{thm:1m} and \autoref{thm:2m} is largely same as \autoref{thm:1} and \autoref{thm:2}, so we will only list some key and different steps instead of all the details.
	\paragraph{Manifolds with global coordinate representations}We will firstly assume that $\mathcal{M}$ has a smooth global coordinate map $\phi:\Omega\rightarrow\mathcal{M}$ where $\Omega$ is a bounded $C^2$ domain in $\mathbb{R}^k$. With such assumption, $T_n$ defined above can be the transportation map from $\frac{\nu}{\nu(\mathcal{M})}$ to $\nu_n$. So there is no difficulty in imitating Lemma \ref{lem:2.4} and \ref{lem:2.5}. What we need to do is proving:
	\[\mathcal{E}_{s_n}\stackrel{\Gamma}{\longrightarrow}\sigma_\eta\mathcal{E},\quad\forall s_n\rightarrow 0\]
	where
	\begin{equation*}
	\begin{aligned}
	\mathcal{E}(u)&=\left\{
	\begin{array}{ll}
	\esssup_{x\in\mathcal{M}}|\nabla_{\mathcal{M}} u| \ \  & \text{if } u\in W^{1,\infty}(\mathcal{M}), \\
	\infty                         & \text{otherwise}.
	\end{array}\right.\\
	&=\left\{\begin{array}{ll}
	\esssup_{x\in\Omega}|(J_\phi(x)^\dagger)^T\nabla(u\circ\phi)(x)| \ \  & \text{if } u\circ\phi\in W^{1,\infty}(\mathcal{M}), \\
	\infty                         & \text{otherwise}.
	\end{array}\right.
	\end{aligned}
	\end{equation*}
	\begin{equation*}
	\begin{aligned}
		\mathcal{E}_{s}(u)&=\frac{1}{s}\esssup_{x\in\mathcal{M}}(\int_{\mathcal{M}}\eta_{s}(|x-y|)|u(x)-u(y)|^p d\nu(y))^\frac{1}{p}\\
		&=\esssup_{x\in\Omega}(\int_{\Omega}\eta_{s}(|\phi(x)-\phi(y)|)(\frac{|u\circ\phi(x)-u\circ\phi(y)|}{s})^p A_\phi(y)dy)^\frac{1}{p}
	\end{aligned}
	\end{equation*} 
	\[A_\phi(y)\coloneqq\sqrt{\det(J_\phi(y)^TJ_\phi(y))}\]
	\indent We firstly prove the liminf inequlity. That is, we want to prove that for any sequence $\{u_n\}_{n\in\mathbb{N}}\subset L^\infty(\mathcal{M})$ converging to $u\in L^\infty(\mathcal{M})$ in $L^\infty(\mathcal{M})$, 
	\[\liminf_{n\rightarrow\infty}\mathcal{E}_{s_n}(u_n)\geq\sigma_\eta\mathcal{E}(u)\]
	To use the methods in Lemma \ref{lem:2.2}, we give two asymptotic approximations of the functional $\mathcal{E}_s$:
	\[\tilde{\mathcal{E}}_s\coloneqq\esssup_{x\in\Omega}(\int_{\Omega}\eta_{s}(|\phi(x)-\phi(y)|)(\frac{|u\circ\phi(x)-u\circ\phi(y)|}{s})^p A_\phi(x)dy)^\frac{1}{p}\] 
	\[\tilde{\tilde{\mathcal{E}}}_s\coloneqq\esssup_{x\in\Omega}(\int_{\Omega}\eta_{s}(|J_\phi(x)(x-y)|)(\frac{|u\circ\phi(x)-u\circ\phi(y)|}{s})^p A_\phi(x)dy)^\frac{1}{p}\]
	And our proof will be divided into three steps:
	\begin{enumerate}
		\item[(1)] $\lim_{n\rightarrow\infty}|\mathcal{E}_{s_n}(u_n)-\tilde{\mathcal{E}}_{s_n}(u_n)|=0$
		\item[(2)] $\lim_{n\rightarrow\infty}|\mathcal{E}_{s_n}(u_n)-\tilde{\tilde{\mathcal{E}}}_{s_n}(u_n)|=0$
		\item[(3)]$\liminf_{n\rightarrow\infty}\tilde{\tilde{\mathcal{E}}}_{s_n}(u_n)\geq\sigma_\eta\mathcal{E}(u)$
	\end{enumerate}
	\indent For (1), w.l.o.g we can assume that $\liminf_{n\rightarrow\infty}\mathcal{E}_{s_n}(u_n)<\infty$ so in the subsequence sense, we have
	\[M\coloneqq \sup_{n}\mathcal{E}_{s_n}(u_n)<\infty\]
	And there also exists positive constant $c\leq C$, such that
	\[c\leq\inf_{x\in\Omega} A_{\phi}(x)\leq\sup_{x\in\Omega} A_{\phi}(x)\leq C\]
	So
	\begin{equation*}
		\begin{aligned}
			&|\mathcal{E}_{s_n}(u_n)^p-\tilde{\mathcal{E}}_{s_n}(u_n)^p|\\
			&\leq \esssup_{x\in\Omega}\int_{\Omega}\eta_{s_n}(|\phi(x)-\phi(y)|)(\frac{|u_n\circ\phi(x)-u_n\circ\phi(y)|}{s_n})^p |A_\phi(x)-A_\phi(y)|dy\\
			&\leq \mathcal{E}_{s_n}(u_n)\sup_{y:|\phi(x)-\phi(y)|\leq s_n}\frac{|A_\phi(x)-A_\phi(y)|}{A_\phi(y)}\\
			&\leq\frac{M}{c}\sup_{y:|\phi(x)-\phi(y)|\leq s_n}|A_\phi(x)-A_\phi(y)|\\
			&\rightarrow 0(n\rightarrow\infty)
		\end{aligned}
	\end{equation*}
	which yields (1).\\
	\indent For (2), w.l.o.g, we assume that $r_\eta=1$ and $\eta(r)>0$ if and only if $r<r_\eta=1$. Since we additionally require that the kernel $\eta$ is Lipschitz continuous. Then, thanks to the Taylor expansion, we have the estimation: 
	\begin{equation*}
		\begin{aligned}
			|\eta_s(|\phi(x)-\phi(y)|)-\eta_s(|J_\phi(x)(x-y)|)|&=\frac{1}{s^k}|\eta(\frac{|\phi(x)-\phi(y)|}{s})-\eta(\frac{|J_\phi(x)(x-y)|}{s})|\\
			&\leq \frac{C_\eta}{s^{k+1}}||\phi(x)-\phi(y)|-|J_\phi(x)(x-y)||\\
			&\leq \frac{C_\eta}{s^{k+1}}|\phi(x)-\phi(y)-J_\phi(x)(x-y)|\\
			&=\frac{o(|x-y|)}{s^{k+1}}
		\end{aligned}
	\end{equation*}
	which yields
	\begin{equation*}
	\begin{aligned}
	&|\mathcal{E}_{s_n}(u_n)^p-\tilde{\tilde{\mathcal{E}}}_{s_n}(u_n)^p|\\
	&\leq\esssup_{x\in\Omega}\int_{\Omega}\mathbf{1}_{\{|\phi(x)-\phi(y)|,|J_\phi(x)(x-y)|\leq s_n\}}\frac{o(|x-y|)}{s_n^{k+1}}(\frac{|u_n\circ\phi(x)-u_n\circ\phi(y)|}{s_n})^pdy
	\end{aligned}
	\end{equation*}
	What we have known from $\sup_{n}\mathcal{E}_{s_n}(u_n)<\infty$ is that
	\[M\coloneqq\sup_n\esssup_{x\in\Omega}\int_{\Omega}\eta_{s_n}(|\phi(x)-\phi(y)|)(\frac{|u_n\circ\phi(x)-u_n\circ\phi(y)|}{s_n})^p dy<\infty\]
	So for any $x,y\in\Omega$, $z$ satisfies that $|\phi(z)-\phi(x)|,|\phi(z)-\phi(y)|\leq\frac{2}{3}s_n$, we have
	\begin{equation*}
		\begin{aligned}
		&\eta_{\frac{7}{6}s_n}(|\phi(x)-\phi(y)|)(\frac{|u_n\circ\phi(x)-u_n\circ\phi(y)|}{s_n})^p\\
		&\leq C_p\eta_{\frac{7}{6}s_n}(|\phi(x)-\phi(y)|)((\frac{|u_n\circ\phi(x)-u_n\circ\phi(z)|}{s_n})^p+(\frac{|u_n\circ\phi(y)-u_n\circ\phi(z)|}{s_n})^p)\\
		&\leq C_{p,\eta}(\eta_{s_n}(|\phi(x)-\phi(z)|)(\frac{|u_n\circ\phi(x)-u_n\circ\phi(z)|}{s_n})^p\\
		&+\eta_{s_n}(|\phi(y)-\phi(z)|)(\frac{|u_n\circ\phi(y)-u_n\circ\phi(z)|}{s_n})^p)
		\end{aligned}
	\end{equation*}
	Integrate with respect to $z$, we get
	\begin{equation*}
	\begin{aligned}
	&\eta_{\frac{7}{6}s_n}(|\phi(x)-\phi(y)|)(\frac{|u_n\circ\phi(x)-u_n\circ\phi(y)|}{s_n})^p\\
	&\leq \frac{C}{s_n^k}\int_{z:|\phi(z)-\phi(x)|,|\phi(z)-\phi(y)|\leq\frac{2}{3}s_n}\eta_{s_n}(|\phi(x)-\phi(z)|)(\frac{|u_n\circ\phi(x)-u_n\circ\phi(z)|}{s_n})^p\\
	&+\eta_{s_n}(|\phi(y)-\phi(z)|)(\frac{|u_n\circ\phi(y)-u_n\circ\phi(z)|}{s_n})^pdz\\
	&\leq \frac{2CM}{s_n^k}
	\end{aligned}
	\end{equation*}
	Integrate with respect to $y$, take supremum with respect to $x$ and $n$, we get
	\[\sup_n\esssup_{x\in\Omega}\int_{\Omega}\eta_{\frac{7}{6}s_n}(|\phi(x)-\phi(y)|)(\frac{|u_n\circ\phi(x)-u_n\circ\phi(y)|}{s_n})^p dy\leq 2\tilde{C}M<\infty\]
	Note that when $n$ is large enough, since $||\phi(x)-\phi(y)|-|J_\phi(x)(x-y)||=o(|x-y|)=o(|\phi(x)-\phi(y)|)$, we have
	\[\{|\phi(x)-\phi(y)|,|J_\phi(x)(x-y)|\leq s_n\}\subset\{|\phi(x)-\phi(y)|\leq\frac{13}{12}s_n\}\]
	So
	\begin{equation*}
		\begin{aligned}
			|\mathcal{E}_{s_n}(u_n)^p-\tilde{\tilde{\mathcal{E}}}_{s_n}(u_n)^p|&\leq 2\tilde{C}M\sup_{x,y\in\Omega}\frac{\mathbf{1}_{\{|\phi(x)-\phi(y)|\leq\frac{13}{12}s_n\}}\frac{o(|x-y|)}{s_n^{k+1}}}{\eta_{\frac{7}{6}s_n}(|\phi(x)-\phi(y)|)}\\
			&\leq C\sup_{x,y\in\Omega}\frac{\mathbf{1}_{\{|\phi(x)-\phi(y)|\leq\frac{13}{12}s_n\}}o(|\phi(x)-\phi(y)|)}{s_n\eta(\frac{6|\phi(x)-\phi(y)|}{7s_n})}\\
			(\eta \text{ is monotonously decreasing})&\leq C\frac{o(\frac{13}{12}s_n)}{s_n\eta(\frac{13}{14})}\\
			&\rightarrow 0(n\rightarrow\infty)
		\end{aligned}
	\end{equation*}
	which yields (2).\\
	\indent For (3), we point out that compared with the proof of Lemma \ref{lem:2.2}, it only needs to add some processes of integral substitution in somewhere. For example, consider the equation \eqref{eq:example}. With the same method, for $\tilde{\tilde{\mathcal{E}_s}}$, we can get
	\[\int_{\mathbb{R}^k}(D^+-D^-)\xi dx=\int_{\mathbb{R}^k}\int_{\mathbb{R}^k}\eta(|J_\phi(x)z|)|\nabla(u\circ\phi)\cdot z|^p\xi(x)A_\phi(x)dzdx,\quad\forall \xi\in\mathcal{C}_c^\infty(\Omega)\]
	Use the singular value decomposition(SVD) of $J_\phi(x)=U^TDV=U^T(\tilde{D}\ O)^TV$, where $U\in\mathbb{R}^{d\times d}, D\in\mathbb{R}^{d\times k},V\in\mathbb{R}^{k\times k}$ and $\tilde{D}\in\mathbb{R}^{k\times k}$ is full rank. Let $\tilde{z}\coloneqq \tilde{D}Vz$, then
	\begin{equation*}
		\begin{aligned}
			\int_{\mathbb{R}^k}(D^+-D^-)\xi dx&=\int_{\mathbb{R}^k}\int_{\mathbb{R}^k}\eta(|\tilde{z}|)|\tilde{D}^{-1}V\nabla(u\circ\phi)(x)\cdot \tilde{z}|^p\xi(x)A_\phi(x)\det\tilde{D}^{-1}d\tilde{z}dx\\
			&=\int_{\mathbb{R}^k}\int_{\mathbb{R}^k}\eta(|\tilde{z}|)|\tilde{D}^{-1}V\nabla(u\circ\phi)(x)\cdot \tilde{z}|^p\xi(x)d\tilde{z}dx
		\end{aligned}
	\end{equation*}
	So we can continue to use the original way to get
	\begin{equation*}
		\begin{aligned}
			\liminf_{n\rightarrow\infty}\tilde{\tilde{\mathcal{E}}}_{s_n}(u_n)&\geq \sigma_\eta\esssup_{x\in\Omega}|\tilde{D}^{-1}V\nabla(u\circ\phi)(x)|\\
			&=\sigma_\eta\esssup_{x\in\Omega}|U^T(\tilde{D}^{-1}\ O)^TV\nabla(u\circ\phi)(x)|\\
			&=\sigma_\eta\esssup_{x\in\Omega}|U^T(\tilde{D}^{-1}\ O)^TV\nabla(u\circ\phi)(x)|\\
			&=\sigma_\eta\esssup_{x\in\Omega}|(J_\phi(x)^\dagger)^T\nabla(u\circ\phi)(x)|\\
			&=\sigma_\eta\mathcal{E}(u)
		\end{aligned}
	\end{equation*}
	\indent To prove the limsup inequality, we need to prove that for any $u\in L^\infty(\mathcal{M})$, there exists a sequence $\{u_n\}_{n\in\mathbb{N}}\subset L^\infty(\mathcal{M})$ converging to $u$ in $L^\infty(\mathcal{M})$ such that
	\[\limsup_{n\rightarrow\infty}\mathcal{E}_{s_n}(u_n)\leq\sigma_\eta \mathcal{E}(u)\]
	Be similar to Lemma \ref{lem:2.3}, we will choose $u_n(y)=\xi_{\frac{1}{n}}*(\hat{u\circ\phi})(\phi^{-1}(y))$, where $\xi$ has the same properties as $\phi$ in the proof of Lemma \ref{lem:2.3}. Like the liminf inequality, our proof is still divided into three steps:
	\begin{enumerate}
		\item[(1)] $\lim_{n\rightarrow\infty}|\mathcal{E}_{s_n}(u_n)-\tilde{\mathcal{E}}_{s_n}(u_n)|=0$
		\item[(2)] $\lim_{n\rightarrow\infty}|\mathcal{E}_{s_n}(u_n)-\tilde{\tilde{\mathcal{E}}}_{s_n}(u_n)|=0$
		\item[(3)]$\limsup_{n\rightarrow\infty}\tilde{\tilde{\mathcal{E}}}_{s_n}(u_n)\leq\sigma_\eta\mathcal{E}(u)$
	\end{enumerate}
	The first and the third steps are the same as the liminf inequality, and the second step is even simpler because $u_n\circ\phi$ now is smooth and $M\coloneqq\sup_n\lVert\nabla(u_n\circ\phi)\rVert_{L^\infty}<\infty$. So
	\begin{equation*}
		\begin{aligned}
			&|\mathcal{E}_{s_n}(u_n)^p-\tilde{\tilde{\mathcal{E}}}_{s_n}(u_n)^p|\\
			&\leq\esssup_{x\in\Omega}\int_{\Omega}\mathbf{1}_{\{|\phi(x)-\phi(y)|,|J_\phi(x)(x-y)|\leq s_n\}}\frac{o(|x-y|)}{s_n^{k+1}}(\frac{|u_n\circ\phi(x)-u_n\circ\phi(y)|}{s_n})^pdy\\
			&\leq\esssup_{x\in\Omega}\int_{\Omega}\mathbf{1}_{\{|x-y|\leq c_\phi s_n\}}\frac{o(|x-y|)}{s_n^{k+1}}(\frac{M|x-y|}{s_n})^pdy\\
			&\leq \esssup_{x\in\Omega}\int_{\Omega}\mathbf{1}_{\{|x-y|\leq c_\phi s_n\}}\frac{o(s_n)}{s_n^{k+1}}dy\\
			&=\frac{o(s_n)}{s_n}\\
			&\rightarrow 0(n\rightarrow\infty)
		\end{aligned}
	\end{equation*}
	which yields (2).\\
	\indent Above all, we have proven \autoref{thm:1m} for $\mathcal{M}$ which has global coordinate representation.
	\paragraph{General situation}Finally, we consider a manifold $\mathcal{M}$ which may not have a global coordinate representation. Our idea mainly based on the finite partition of $\mathcal{M}$ we have done when defining the transportation map $T_n$. That is
	\[\mathcal{M}=\bigcup_{i=1}^m\mathcal{M}^i=\bigcup_{i=1}^{m}\tilde{\mathcal{M}}^i\cup\varGamma\]
	Note that $\mathcal{M}^i,\tilde{\mathcal{M}}^i$ have global coordinate representations so convergence results we have proven hold on them. In the remaining part, we will show that this can lead to the convergence results on $\mathcal{M}$. Firstly, we define the restricted versions of the functionals we need:
	\begin{equation*}
	E_n^i(u)\coloneqq
	\begin{cases}
	\frac{1}{s_n}\mathop{\max}\limits_{x\in\mathcal{M}_n\cap\mathcal{M}^i}(\frac{1}{|\mathcal{M}_n\cap\mathcal{M}^i|}\sum_{y\in\mathcal{M}_n\cap\mathcal{M}^i}\eta_{s_n}(|x-y|)|u(x)-u(y)|^p)^\frac{1}{p}
	& \text{if } u=g \\
	&\text{ on } \mathcal{O}\cap\mathcal{M}^i,\\
	\infty
	& \text{otherwise.}
	\end{cases}
	\end{equation*}
	\begin{equation*}
	E_{n,\mathrm{cons}}^i(u)\coloneqq
	\begin{cases}
	E_{n}^i(\tilde{u}) & \text{if } u=\tilde{u}\circ T_n^i \text{ for some } \tilde{u}:\mathcal{M}_n\cap\mathcal{M}^i\longrightarrow\mathbb{R}, \\
	\infty           & \text{otherwise}.
	\end{cases}
	\end{equation*}
	\begin{equation*}
	\mathcal{E}^i(u)\coloneqq
	\begin{cases}
	\esssup_{x\in\mathcal{M}^i}|\nabla_{\mathcal{M}} u| & \text{if } u\in W^{1,\infty}(\mathcal{M}^i), \\
	\infty                         & \text{otherwise}.
	\end{cases}
	\end{equation*}
	\begin{equation*}
	\mathcal{E}_{\mathrm{cons}}^i(u)\coloneqq\left\{\begin{array}{ll}
	\mathcal{E}^i(u) \ \  & \text{if } u=g\text{ on } \mathcal{O}\cap\mathcal{M}^i, \\
	\infty         & \text{otherwise}.
	\end{array}\right.
	\end{equation*}
	\begin{equation*}
	\tilde{E}_n^i(u)\coloneqq
	\begin{cases}
	\frac{1}{s_n}\mathop{\max}\limits_{x\in\mathcal{M}_n\cap\tilde{\mathcal{M}}^i}(\frac{1}{|\mathcal{M}_n\cap\tilde{\mathcal{M}}^i|}\sum_{y\in\mathcal{M}_n\cap\tilde{\mathcal{M}}^i}\eta_{s_n}(|x-y|)|u(x)-u(y)|^p)^\frac{1}{p}
	& \text{if } u=g \\
	&\text{ on } \mathcal{O}\cap\tilde{\mathcal{M}}^i,\\
	\infty
	& \text{otherwise.}
	\end{cases}
	\end{equation*}
	\begin{equation*}
	\tilde{E}_{n,\mathrm{cons}}^i(u)\coloneqq
	\begin{cases}
	\tilde{E}_{n}^i(\tilde{u}) & \text{if } u=\tilde{u}\circ \tilde{T}_n^i \text{ for some } \tilde{u}:\mathcal{M}_n\cap\tilde{\mathcal{M}}^i\longrightarrow\mathbb{R}, \\
	\infty           & \text{otherwise}.
	\end{cases}
	\end{equation*}
	\begin{equation*}
	\tilde{\mathcal{E}}^i(u)\coloneqq
	\begin{cases}
	\esssup_{x\in\tilde{\mathcal{M}}^i}|\nabla_{\mathcal{M}} u| & \text{if } u\in W^{1,\infty}(\tilde{\mathcal{M}^i}), \\
	\infty                         & \text{otherwise}.
	\end{cases}
	\end{equation*}
	
	\begin{equation*}
	\tilde{\mathcal{E}}_{\mathrm{cons}}^i(u)\coloneqq\left\{\begin{array}{ll}
	\tilde{\mathcal{E}}^i(u) \ \ \  & \text{if } u=g\text{ on } \mathcal{O}\cap\tilde{\mathcal{M}}^i, \\
	\infty         & \text{otherwise}.
	\end{array}\right.
	\end{equation*}
	
	What we have known is that for all $i\in\{1,2,\dots,m\}$
	\[E_{n,\mathrm{cons}}^i\stackrel{\Gamma}{\rightarrow}\sigma_{\eta}\nu(\mathcal{M}^i)^{-\frac{1}{p}}\mathcal{E}_{\mathrm{cons}}^i\]
	\[\tilde{E}_{n,\mathrm{cons}}^i\stackrel{\Gamma}{\rightarrow}\sigma_{\eta}\nu(\tilde{\mathcal{M}}^i)^{-\frac{1}{p}}\tilde{\mathcal{E}}_{\mathrm{cons}}^i\]
	And we hope to prove
	\[E_{n,\mathrm{cons}}\stackrel{\Gamma}{\rightarrow}\sigma_{\eta}\nu(\mathcal{M})^{-\frac{1}{p}}\mathcal{E}_{\mathrm{cons}}\]
	\indent For the liminf inequality, for all $u_n\stackrel{L^\infty}{\rightarrow} u$, we have
	\begin{equation*}
		\begin{aligned}
			\sigma_{\eta}^p\nu(\mathcal{M})^{-1}\mathcal{E}_{\mathrm{cons}}^p(u)&=\sigma_{\eta}^p\nu(\mathcal{M})^{-1}\max_{i}(\tilde{\mathcal{E}}_{\mathrm{cons}}^i)^p(u)\\
			&=\sigma_{\eta}^p\nu(\mathcal{M})^{-1}\max_{i}(\tilde{\mathcal{E}}_{\mathrm{cons}}^i)^p(u|_{\tilde{\mathcal{M}}^i})\\
			&\leq\nu(\mathcal{M})^{-1}\max_{i}\nu(\tilde{\mathcal{M}}^i)\liminf_{n\rightarrow\infty}(\tilde{E}_{n,\mathrm{cons}}^i)^p(u_n|_{\tilde{\mathcal{M}}^i})\\
			&=\liminf_{n\rightarrow\infty}\max_{i}\frac{\nu(\tilde{\mathcal{M}}^i)}{\nu(\mathcal{M})}(\tilde{E}_{n,\mathrm{cons}}^i)^p(u_n)
		\end{aligned}
	\end{equation*}
	W.l.o.g, we can assume that $\sup_nE_{n,\mathrm{cons}}(u_n)<\infty$ as same as before. With such assumption,
	\begin{equation*}
		\begin{aligned}
			\liminf_{n\rightarrow\infty}E_{n,\mathrm{cons}}^p(u_n)&=\liminf_{n\rightarrow\infty}E_{n}^p(u_n)\\
			&\geq\liminf_{n\rightarrow\infty}\max_i\frac{\nu(\tilde{\mathcal{M}}^i)}{\nu(\mathcal{M})}(\tilde{E}_{n}^i)^p(u_n)\\
			&=\liminf_{n\rightarrow\infty}\frac{\nu(\tilde{\mathcal{M}}^i)}{\nu(\mathcal{M})}(\tilde{E}_{n,\mathrm{cons}}^i)^p(u_n)\\
			&\geq\sigma_{\eta}^p\nu(\mathcal{M})^{-1}\mathcal{E}_{\mathrm{cons}}^p(u)
		\end{aligned}
	\end{equation*}
	which yields the liminf inequality. To get the second inequality, we use the Strong Law of Large Number(SLLN), from which we have $\lim_{n\rightarrow\infty}\frac{|\mathcal{M}_n\cap\tilde{\mathcal{M}}^i|}{|\mathcal{M}_n|}=\frac{\nu(\tilde{\mathcal{M}}^i)}{\nu(\mathcal{M})}$ with probability one.\\
	\indent For the limsup inequality, for all $u\in L^\infty(\mathcal{M})$, we will choose $u_n\coloneqq u|_{\mathcal{M}_n}\circ T_n$. Then $u_n=u$ on $\mathcal{M}_n$ and $\lVert u_n-u\rVert_{L^\infty}\rightarrow 0$ because $T_n$ satisfies \eqref{eq:e2}. We hope to prove
	\[\limsup_{n\rightarrow\infty}E_{n,\mathrm{cons}}(u_n)\leq\sigma_{\eta}\nu(\mathcal{M})^{-\frac{1}{p}}\mathcal{E}_{\mathrm{cons}}(u)\]
	Actually,
	\begin{equation*}
		\begin{aligned}
			\limsup_{n\rightarrow\infty}E_{n,\mathrm{cons}}^p(u_n)&=\limsup_{n\rightarrow\infty}E_{n}^p(u_n)\\
			&=\limsup_{n\rightarrow\infty}\max_i\frac{\nu(\mathcal{M}^i)}{\nu(\mathcal{M})}(E_{n}^i)^p(u_n)\\
			&=\limsup_{n\rightarrow\infty}\max_i\frac{\nu(\mathcal{M}^i)}{\nu(\mathcal{M})}(E_{n,\mathrm{cons}}^i)^p(u|_{\mathcal{M}\cap\mathcal{M}^i}\circ T_n^i)\\
			&\leq\max_i\frac{\nu(\mathcal{M}^i)}{\nu(\mathcal{M})}\limsup_{n\rightarrow\infty}(E_{n,\mathrm{cons}}^i)^p(u|_{\mathcal{M}\cap\mathcal{M}^i}\circ T_n^i)\\
			&\leq\sigma_{\eta}^p\nu(\mathcal{M})^{-1}\max_i \mathcal{E}_{\mathrm{cons}}^i(u)\\
			&=\sigma_{\eta}^p\nu(\mathcal{M})^{-1}\mathcal{E}_{\mathrm{cons}}(u)
		\end{aligned}
	\end{equation*}
	So we prove \autoref{thm:1m} for general manifolds.\\
	\indent For \autoref{thm:2m}, we only need to give the relatively compactness of minimizers. Suppose that $\{u_n\}$ satisfies $\sup_n E_{n,\mathrm{cons}}(u_n)<\infty$, then for any $i\in\{1,2\dots,m\}$, $\sup_nE_n^i(u_n)<\infty$. Note that $\mathcal{M}^i$ has a smooth global coordinate map $\phi:\Omega\rightarrow\mathcal{M}^i$. We have
	\begin{equation*}
		\begin{aligned}
			&E_n^i(u_n)\\
			&=\frac{1}{s_n}\mathop{\max}\limits_{x\in\mathcal{M}_n\cap\mathcal{M}^i}(\frac{1}{|\mathcal{M}_n\cap\mathcal{M}^i|}\sum_{y\in\mathcal{M}_n\cap\mathcal{M}^i}\eta_{s_n}(|x-y|)|u_n(x)-u_n(y)|^p)^\frac{1}{p}\\
			&=\frac{1}{s_n}\mathop{\max}\limits_{x\in\phi^{-1}(\mathcal{M}_n\cap\mathcal{M}^i)}(\frac{1}{|\mathcal{M}_n\cap\mathcal{M}^i|}\sum_{y\in\phi^{-1}(\mathcal{M}_n\cap\mathcal{M}^i)}\eta_{s_n}(|\phi(x)-\phi(y)|)|u_n\circ\phi(x)-u_n\circ\phi(y)|^p)^\frac{1}{p}\\
			&\geq \frac{C_\phi^k}{s_n}\mathop{\max}\limits_{x\in\phi^{-1}(\mathcal{M}_n\cap\mathcal{M}^i)}(\frac{1}{|\mathcal{M}_n\cap\mathcal{M}^i|}\sum_{y\in\phi^{-1}(\mathcal{M}_n\cap\mathcal{M}^i)}\eta_{C_\phi s_n}(|x-y|)|u_n\circ\phi(x)-u_n\circ\phi(y)|^p)^\frac{1}{p}
		\end{aligned}
	\end{equation*}
	where the last inequality holds because $\phi$ is smooth so that it is also Lipschitz continuous with Lipschitz constant $C_\phi$ and $\eta$ is monotonously decreasing. So with Lemma \ref{lem:compactness}, we have $\{u_n\circ\phi\}$ has a convergent subsequence in $L^\infty(\Omega)$, i.e. $\{u_n\}$ has a convergent subsequence in $L^\infty(\mathcal{M}^i)$. This holds for any $i\in\{1,2,\dots,m\}$ which yields $\{u_n\}$ has a convergent subsequence in $L^\infty(\mathcal{M}^i),\forall i\in\{1,2,\dots,m\}$, also in $L^\infty(\mathcal{M})$. The relatively compactness is proven.

\section{Algorithm}
\label{sec:alg}

In this section, we propose an algorithm to solve the discrete infinity Laplacian problem \eqref{prob_cons} for the special case of $p = 2$. Let us denote by $f$ the objective function as follows
\begin{equation*}
f(\mathbf{u}) = \mathop{\max}_{i\in I} \big\{\sum_{j \in I} w_{ij}(u_i-u_j)^2\big\}  + \alpha \sum_{i\in I} \sum_{j \in I} w_{ij}(u_i-u_j)^2
\end{equation*}
in which $I=\{1,\dots,m+n\}$, $\{u_i\}_{i=1}^m$ is the set of label values, $\mathbf{u} = (u_{m+1},\dots,u_{m+n})^T$ is the optimization variable. $w_{ij} \geq 0$ is the weight of the points $x_i$ and $x_j$. The minimizer of $f$ exists because it is a convex function defined on $\mathbb{R}^n$ with a lower bound $0$. Let $G$ be the graph constructed with the point cloud $P$. We assign an edge to points $x_i$ and $x_j$ if and only if $w_{ij} > 0$. In the following we assume that the graph $G$ is connected. When $\alpha > 0$, the minimizer is unique because $f$ is strictly convex. When $\alpha = 0$, a counterexample is provided in appendix \ref{appendix:non-uniqueness} to show that the minimizer may not be unique.

To solve the discrete problem \eqref{prob_cons} with $p=2$, we write it equivalently as follows.
\begin{equation*}
\begin{aligned}
\mathop{\min}_{\mathbf{u}, \mathbf{D}} & \Big\{\mathop{\max}_{i \in I} \big\{\sum_{j \in I} D_{ij}^2\big\}  + \alpha \sum_{i \in I} \sum_{j \in I} D_{ij}^2\Big\}\\
s.t.\quad & D_{ij}= \sqrt{w_{ij}}(u_i-u_j), \quad i,j \in I
\end{aligned}
\end{equation*}
in which $I = \{1,\dots,m+n\}$ and $\mathbf{D} = (D_{ij})$ is a matrix of size $(m+n) \times (m+n)$ . We use the split Bregman method \cite{SBI} to enforce these constraints, which leads to the following algorithm.

\begin{itemize}
    \item Let $k=0$, $D_{ij}^0 = 0$, $q_{ij}^0 = 0$ with $i,j = 1,2,\dots, m+n$.
    \item Update $u_{m+j}$ with $j=1,2,\dots,n$ by solving
\end{itemize}   

\begin{equation}
\label{update_u}
    \{u^{k+1}_{m+j}\}_{j=1}^n = \mathop{\arg\min}_{\{u_{m+j}\}_{j=1}^n} \sum_{i \in I} \nu_i \sum_{j \in I} (D_{ij}^{k} - \sqrt{w_{ij}}(u_i-u_j)+q_{ij}^k)^2
\end{equation}
in which $\{\nu_i\}_{i=1}^{m+n}$ are fixed parameters.

\begin{itemize}
\item Update $D_{ij}$ with $i,j = 1,2,\dots, m+n$ by solving
\end{itemize}

\begin{equation}
\label{update_D}
\begin{aligned}
    \{D_{ij}^{k+1}\}_{i,j \in I} = \mathop{\arg\min}_ { \{D_{ij}\}_{i,j\in I}} \Big\{ &\mathop{\max}_{i\in I} \big\{\sum_{j\in I}  D_{ij}^2 \big\} + \alpha \sum_{i\in I} \sum_{j \in I} D_{ij}^2\\
    & + \sum_{i \in I} \nu_i \sum_{j \in I} (D_{ij} - \sqrt{w_{ij}}(u_i^{k+1}-u_j^{k+1})+q_{ij}^k)^2 \Bigg\}
\end{aligned}
\end{equation}
\begin{itemize}
    \item Update $q_{ij}$ with $i,j=1,2,\dots,m+n$ as
\end{itemize}
\begin{equation*}
    q_{ij}^{k+1} = q_{ij}^k + D_{ij}^{k+1} - \sqrt{w_{ij}}(u_i^{k+1} - u_j^{k+1})
\end{equation*}

The sub-problem \eqref{update_u} is a least squares problem, whose solution can be obtained by solving the linear system
\begin{equation}
\label{lin_solve_u}
    \sum_{ \substack{j \in I\\j \neq i}}(\nu_i w_{ij} + \nu_j w_{ji}) (u_i - u_j) = \sum_{\substack{j \in I\\j \neq i}} (\nu_i \sqrt{w_{ij}}s_{ij}^k - \nu_j \sqrt{w_{ji}}s_{ji}^k)
\end{equation}
in which $i = m+1, \dots, m+n$ and $s_{ij} = D_{ij}^k + q_{ij}^k$.

To solve the sub-problem of updating $\mathbf{D}$, we rewrite the objective function of  \eqref{update_D} as follows.
\begin{equation*}
\begin{aligned}
    \varphi(\mathbf{D}) = &\mathop{\max}_{i\in I} \big\{\sum_{j \in I}  D_{ij}^2 \big\} + \alpha \sum_{i\in I} \sum_{j \in I} D_{ij}^2\\
    & + \sum_{i \in I} \nu_i \sum_{j \in I} (D_{ij} - \sqrt{w_{ij}}(u_i^{k+1}-u_j^{k+1})+q_{ij}^k)^2\\
    =& \mathop{\max}_{i\in I} \big\{\sum_{j \in I}  D_{ij}^2 \big\} + \sum_{i \in I} \sum_{j \in I} a_i ( D_{ij} - C_{ij}^k )^2
\end{aligned}
\end{equation*}
in which $a_i = \alpha + \nu_i > 0$ and $C_{ij}^{k+1} = \frac{\nu_i}{\alpha + \nu_i}(\sqrt{w_{ij}}(u_i^{k+1} - u_j^{k+1}) - q_{ij}^k)$. For simplicity, we omit the superscript $k$ or $k+1$ in the following. By denoting $\mathbf{D}_i = (D_{i1},\dots,D_{i,m+n})^T$ and  $\mathbf{C}_i = (C_{i1},\dots,C_{i,m+n})^T$ , we have
\begin{equation}
\label{solve_D_1}
    \varphi(\mathbf{D}) = \max_{i \in I} \Vert \mathbf{D}_i \Vert^2 + \sum_{i=1}^{m+n} a_i \Vert \mathbf{D}_i - \mathbf{C}_i \Vert^2
\end{equation}
One easily checks that the desired solution is given by $\mathbf{D}_i^* = r_i \mathbf{C}_i$, in which $r_i$ is the solution of the following optimization problem
\begin{equation}
\label{solve_D_2}
\mathop{\min}_{\{x_i\}_{i \in I}} \Big\{  \mathop{\max}_{i \in I } \{x_i^2\} + \sum_{i \in I} a_i (x_i - \Vert \mathbf{C}_i \Vert)^2 \Big\}
\end{equation}
In fact, it suffices to consider the problem of the following type, which can be solved efficiently by Algorithm~\ref{alg:max_L2_alg}. A detailed discussion is provided in appendix \ref{subproblem_update_D}.

\begin{equation}
\label{general_solve_D_2}
\mathop{\min}_{\{x_i\}_{i=1}^{N}} \Big\{  \mathop{\max}_{i \in \{1,\dots,N\}} \{x_i^2\} + \sum_{i=1}^{N} a_i (x_i - c_i)^2 \Big\}, \ c_1 > c_2 > \dots > c_N \geq 0.
\end{equation}

\begin{algorithm}[t]
\caption{Solving the problem \eqref{general_solve_D_2}.}
\label{alg:max_L2_alg}

\begin{algorithmic}
\REQUIRE $N$, $\{a_{i}\}_{i=1}^N$, $\{c_i\}_{i=1}^N$ with $a_i \geq 0$, $c_1 > c_2 > \dots > c_N \geq 0$.
\ENSURE $\{x_{i}^*\}_{i=1}^N$.

\STATE{Let $\phi_i=\frac{a_i}{a_i+1}c_i$, $i=1,\dots,N$. $\phi = max_j\{\phi_j\}_{j=1}^N$.}
\STATE{Let $K_1 = \{k:c_k \leq \phi\}$, $K_2 = \{k:c_k>\phi\}$. $\forall k \in K_1$, $x_k^* = c_k$.}
\IF{$|K_2| = 1$}
    \STATE $\forall k \in K_2$, $x_{k^*}^* =\phi$.
\ELSE
    \STATE Let $t = |K_2|$. For $k=1,2,\dots,t-1$, compute $A_k=\sum_{j=1}^ka_{j}$,  $\delta_k = c_{k} - c_{k+1}$. 
    \STATE Let $T = max_k\{k:2\leq k\leq t, \sum_{i=1}^{k-1}A_i\delta_i\leq c_{k}\}$ and $\phi^* = \frac{\sum_{j=1}^{T}a_{j}c_{j}}{\sum_{j=1}^{T}a_{j}+1}$.
    \STATE For $1\leq p\leq T$, $x_p^* = \phi^*$. For $T+1 \leq p \leq t$, $x_{p}^* = c_{p}$.
\ENDIF
\end{algorithmic}
\end{algorithm}

\sloppy The choice of parameters $\nu_i$ plays an important role for fast convergence and a heuristic procedure is used to set these parameters adaptively. We denote by $\mathbf{T}^1 = (t_{ij}^1)_{\vert P \vert \times \vert P \vert} = (\sqrt{w_{ij}}(u_i^1 - u_j^1))_{\vert P \vert \times \vert P \vert}$ the matrix consisting of the values of the non-local gradient and choose $\nu_i = c$ for all $i \in I$, in which $c$ is a constant to be determined. In the first iteration, $\mathbf{u}^1$ can be updated with \eqref{lin_solve_u}, independent of $c$. To obtain $\mathbf{D}^1$, we solve the problem \eqref{update_D} with $q_{ij}^0 = 0$. In this problem, the solution $\mathbf{D}^1$ can be viewed as a kind of threshold with respect to $\mathbf{T}^1$. In Algorithm~\ref{alg:choose_nu} we provide an approach to determine $c$ such that $\Vert \mathbf{D}^1 - \mathbf{T}^1  \Vert_\text{F} / \Vert \mathbf{T}^1 \Vert_\text{F} = 1/2$. The preset threshold $\epsilon$ is chosen as $10^{-4}$.

Finally, the algorithm for solving the discrete problem \eqref{prob_cons} with $p = 2$ is summarized in Algorithm~\ref{alg:main}. When the sub-problems \eqref{update_u} and \eqref{update_D} are solved accurately, convergence to the minimizer is guaranteed under the ADMM framework with the Slater condition \cite{Opt}. In the proposed algorithm, $\mathbf{D}$ can be updated accurately with \eqref{solve_D_1}, \eqref{solve_D_2} and Algorithm~\ref{alg:max_L2_alg}. To update $\mathbf{u}$, the problem \eqref{lin_solve_u} is solved efficiently using the minimal residual method implemented in the Matlab software \cite{minres1,minres2}. In all experiments, the iteration stops when the objective value decreases by no more than $0.0001\%$ in two adjacent steps, or the maximum number of iterations is reached.

\begin{algorithm}[t]
\caption{Choosing the parameter $c$ adaptively.}
\label{alg:choose_nu}
\begin{algorithmic}
\REQUIRE $\alpha$, $\epsilon$, $\mathbf{T}^1$.
\ENSURE $c$.
\IF{$\alpha > 0$} 
\STATE $c = \alpha$.
\ELSE
\STATE $c = 1$.
\ENDIF
\STATE {Update $\mathbf{D}^1$ by \eqref{solve_D_1}, \eqref{solve_D_2}, and Algorithm  \ref{alg:max_L2_alg} with $\nu_i = c$.}
\WHILE{$\Big| \frac{\Vert \mathbf{D}^1 - \mathbf{T}^1  \Vert ^2}{\Vert \mathbf{T}^1 \Vert^2} - \frac{1}{4}\Big| > \epsilon$}
\STATE $c = 4c \cdot \frac{\Vert \mathbf{D}^1 - \mathbf{T}^1 \Vert ^2}{\Vert \mathbf{T}^1 \Vert^2}$.
\STATE {Update $\mathbf{D}^1$ by \eqref{solve_D_1}, \eqref{solve_D_2}, and Algorithm~\ref{alg:max_L2_alg} with $\nu_i = c$.}
\ENDWHILE
\end{algorithmic}
\end{algorithm}

\begin{algorithm}[t]
\caption{Solving the discrete infinity Laplacian problem.}
\label{alg:main}
\begin{algorithmic}

\REQUIRE $m$, $n$, $\{w_{ij}\}_{i,j=1}^{m+n}$, $\{u_i\}_{i=1}^m$.
\ENSURE $\{u_{m+j}\}_{j=1}^n$.
\STATE{$k=0$, $D_{ij}^0 = 0$, $q_{ij}^0 = 0$, $i,j = 1,2,\dots, m+n$.}

\STATE{Update $\{u_{m+j}^{1}\}_{j=1}^n$ with \eqref{lin_solve_u}.}
\STATE{Compute $\mathbf{T}^1 = (t_{ij}^1) = (\sqrt{w_{ij}}(u_i^1 - u_j^1))$}.
\STATE{Choose $c$ with Algorithm~\ref{alg:choose_nu}. Set $\nu_i = c$, $i=1,\dots,m+n$}.
\STATE{Update $\{D_{ij}^{1}\}_{i,j=1}^{m+n}$ with \eqref{solve_D_1}, \eqref{solve_D_2} and Algorithm~\ref{alg:max_L2_alg}. $k=1$.}
\WHILE{not converge}
\STATE{Update $\{u_{m+j}^{k+1}\}_{j=1}^n$ with \eqref{lin_solve_u}.}
\STATE{Update $\{D_{ij}^{k+1}\}_{i,j=1}^{m+n}$ with \eqref{solve_D_1}, \eqref{solve_D_2} and Algorithm~\ref{alg:max_L2_alg}.}
\STATE{$q_{ij}^{k+1} = q_{ij}^k + D_{ij}^{k+1} - \sqrt{w_{ij}}(u_i^{k+1} - u_j^{k+1})$.}
\ENDWHILE
\end{algorithmic}
\end{algorithm}

We remark that solutions of the Graph Laplacian (GL) \cite{GL} and Weighted Non-local Laplacian (WNLL) \cite{WNLL} methods can also be obtained with the proposed algorithm. In GL, the objective function writes
\begin{equation*}
\mathop{\min}_{\{u_{m+j}\}_{j=1}^n} \sum_{i\in I} \sum_{j \in I} w_{ij}(u_i-u_j)^2
\end{equation*}
which is the regularization term in problem \eqref{prob_cons}. The objective function of WNLL is
\begin{equation*}
\mathop{\min}_{\{u_{m+j}\}_{j=1}^n} \Big\{\frac{\vert P \vert}{\vert S \vert}\sum_{i\in I_1} \sum_{j \in I} w_{ij}(u_i-u_j)^2 + \sum_{i\in I_2} \sum_{j \in I} w_{ij}(u_i-u_j)^2\Big\}
\end{equation*}
in which $I_1 = \{1,\dots,m\}$, $I_2 = \{m+1,\dots,m+n\}$, $I = I_1 \cup I_2$. $S$ is the label set, $T$ is the unlabel set, $P = S\cup T$ is the point cloud. In fact, if we take $\nu_i = 1$ for all $i \in I$, the result of the first iteration of Algorithm  \ref{alg:main} is exactly the solution of GL. If we take $\nu_i = \vert P \vert / \vert S \vert$ for $i \in I_1$ and $\nu_j = 1$ for $j \in I_2$, the result of the first iteration turns out to be the solution of WNLL. In this sense, GL or WNLL can be viewed as a first step to solve the IL model.

\section{Experimental results}
\label{sec:experiments}

In the following experiments we compare the interpolation results of the proposed IL method with GL and WNLL. The algorithm is implemented with Matlab on a laptop equipped with CPU Intel i5-1135G7 2.4GHz.

\subsection{A toy example in two dimensions}

\begin{figure}[t]
  \centering
   \subfigure[$f(x,y) = \sin(x)\cos(y) $]{\label{fig:1a}\includegraphics[width=0.45\textwidth]{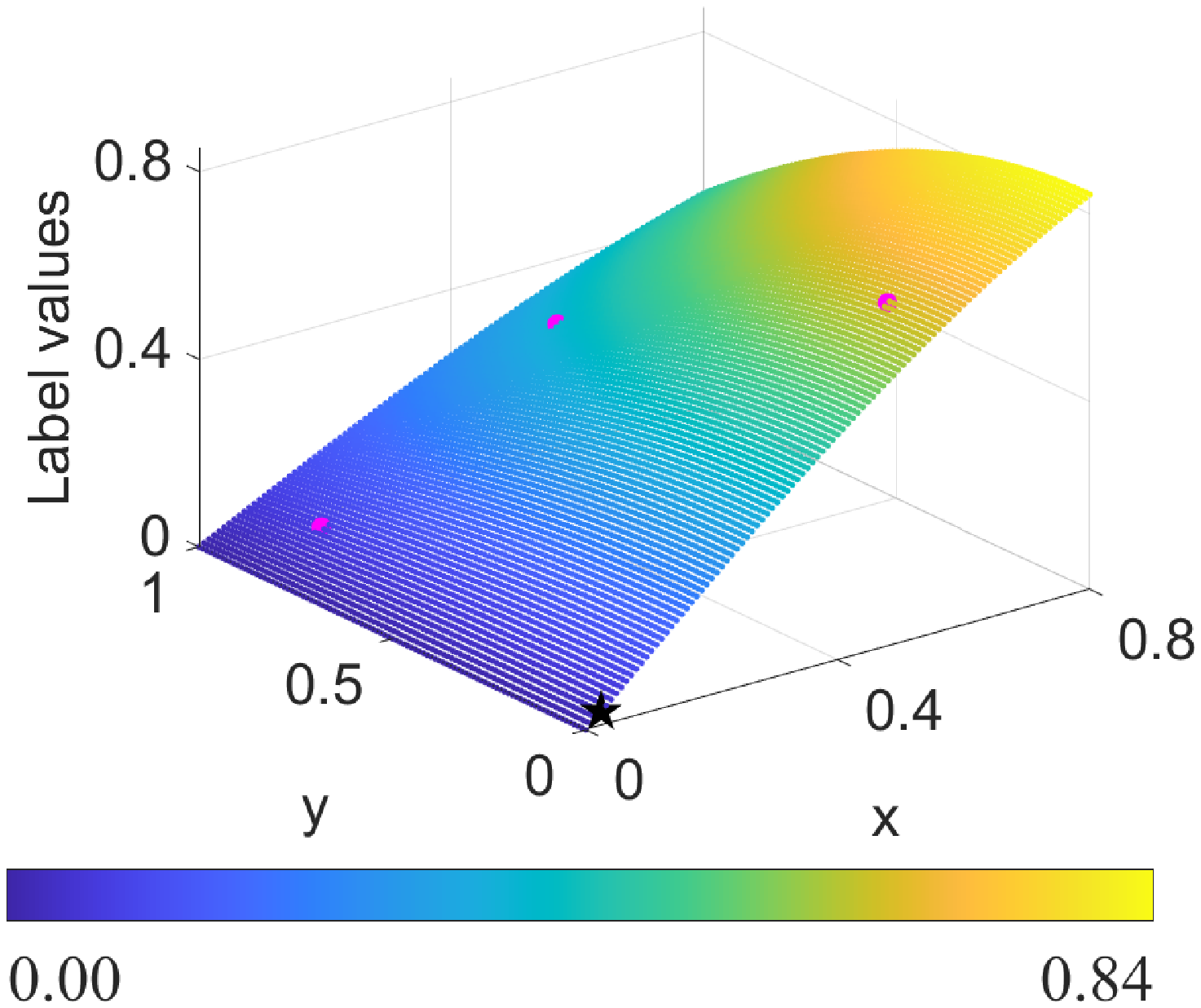}}
   \subfigure[GL interpolation]{\label{fig:1b}\includegraphics[width=0.45\textwidth]{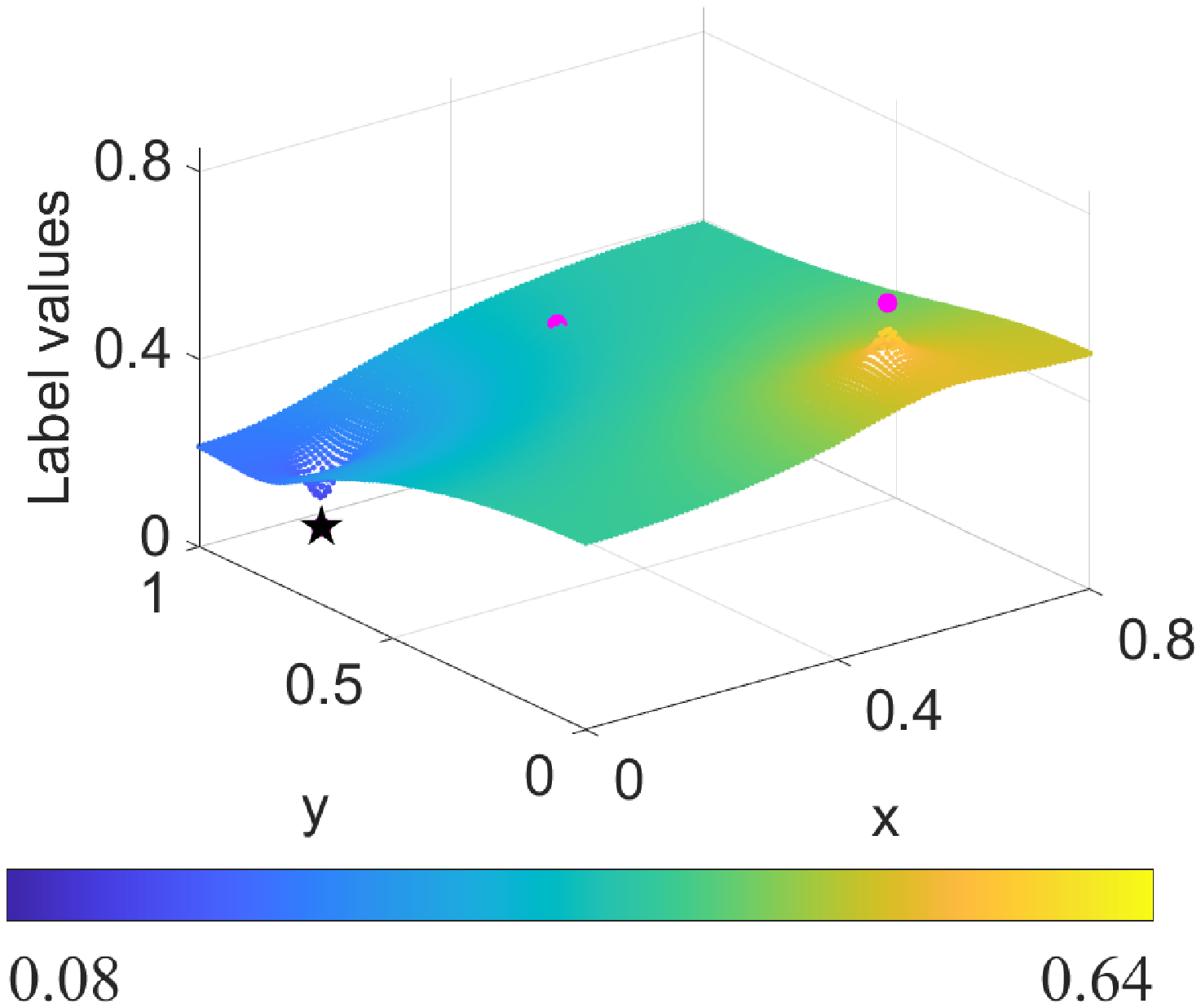}}
   \subfigure[WNLL interpolation]{\label{fig:1c}\includegraphics[width=0.45\textwidth]{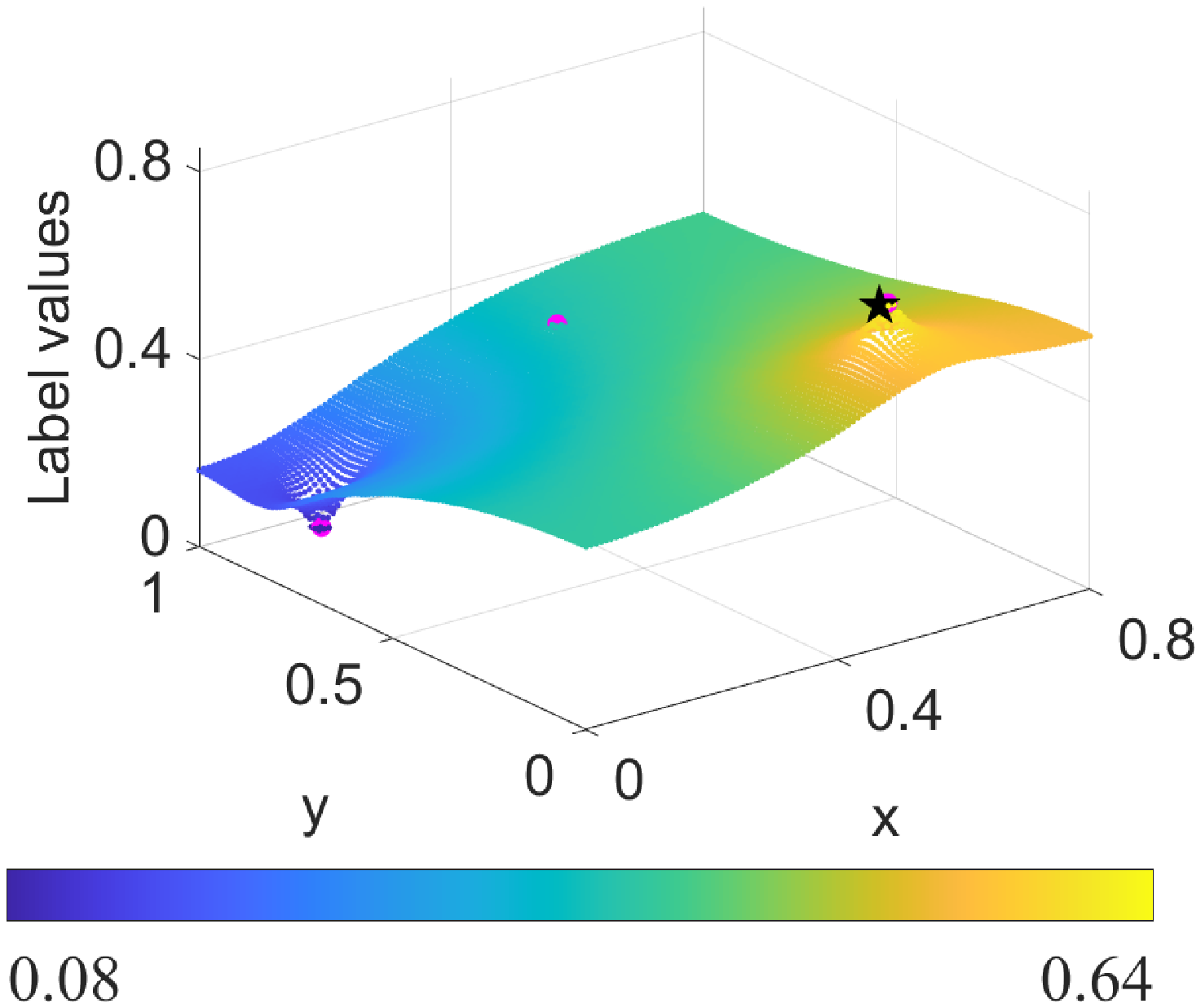}}
   \subfigure[IL interpolation]{\label{fig:1d}\includegraphics[width=0.45\textwidth]{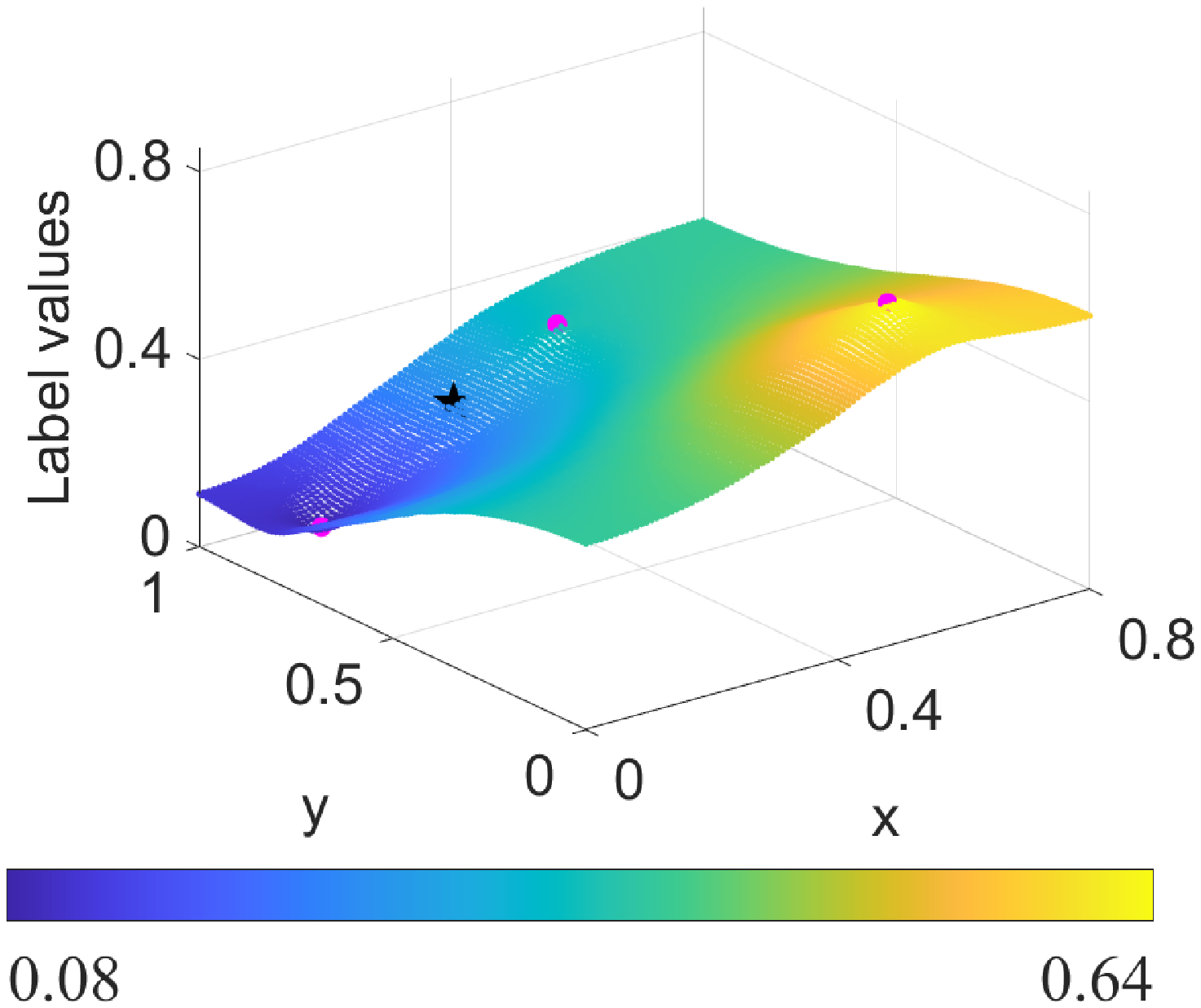}}
  \caption{The function that generates the label values is shown in (a). Interpolation results of Graph Laplacian, Weighted Non-local Laplacian and Infinity Laplacian are shown in (b), (c) and (d) respectively. The label values are shown in magenta. The black stars show points at which the non-local gradient reaches the maximum.}
  \label{fig:fig1}
\end{figure}

We use a toy example in two dimensions to test the proposed algorithm. The set of labeled points is chosen as $S = \{ (\frac{\sqrt{2}}{2},1-\frac{\sqrt{3}}{10}), (\frac{\sqrt{2}}{10},\frac{\sqrt{5}}{20}), (\frac{\sqrt{3}}{3},\frac{\sqrt{11}}{4})\}$ and the set of unlabeled points is $T = \{(0.01i,0.01j)\}_{i,j=0}^{100}$. We use the function $f(x,y) = \sin(x)\cos(y)$ to generate the label values. The weights are constructed by $w_{ij} = \operatorname{exp}(-\Vert x_i-x_j \Vert^2 / \sigma^2)$ with $\sigma = 0.02$. In order to obtain a sparse weight matrix, the weight is truncated to 10 nearest neighbors, which are searched using the approximate nearest neighbor (ANN) algorithm in the vl-feat package \cite{vlfeat}. For the IL method, we choose $\alpha = 0$. Interpolation results of graph Laplacian (GL), weighted non-local Laplacian (WNLL) and the proposed infinity Laplacian (IL) models are illustrated in \autoref{fig:fig1}. It is shown that GL fails to interpolate the points smoothly around the labeled points. The WNLL algorithm also yields large difference of the label values in the neighborhoods of the labeled points. In contrast, the proposed infinity Laplacian model provides an interpolation which is globally smooth. The infinity norm of the non-local gradient of the function $f(x,y) = \sin(x)\cos(y)$, which generates the label values, is $6.85 \times 10^{-4}$. The infinity norms of the non-local gradient of GL, WNLL and IL are $3.52\times 10^{-2}$, $4.36\times 10^{-3}$ and $1.92\times 10^{-4}$ respectively. We run each algorithm for five times independently and the average computation times of GL, WNLL and IL are 0.08 s, 0.13 s and 9.04 s respectively. 

\begin{figure}[t]
  \centering
  \includegraphics[width=0.5\textwidth]{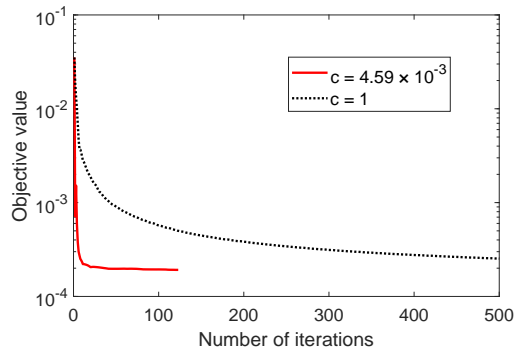}
  \caption{Objective value of the IL model with different choices of $c$. The parameter $c$ adaptively chosen by Algorithm~\ref{alg:choose_nu} yields fast converegence (red line) and meets the stopping criterion at iteration 123.}
  \label{fig:fig2}
\end{figure}

In the IL algorithm, the parameter $c = \nu_i$ adaptively chosen by Algorithm  \ref{alg:choose_nu} is $c^* = 4.59 \times 10^{-3}$ and the stopping criterion is fulfilled at iteration 123 with $ (f(\mathbf{u}^{123}) - f(\mathbf{u}^{122})) / f(\mathbf{u}^{122}) = -1.11\times 10^{-7}$. For a comparison, we fix $\nu_i = c_0 = 1$ and run the code again. The stopping criterion is not reached in the first 500 iterations, and the objective value at the last iteration is $2.53\times 10^{-4}$. As is shown in \autoref{fig:fig2}, the parameter chosen by Algorithm~\ref{alg:choose_nu} leads to fast convergence of the proposed algorithm.

\subsection{A toy example of image inpainting}

\begin{figure}[t]
  \centering
  \subfigure[clear image]{\label{fig:3a}\includegraphics[width=0.3\textwidth]{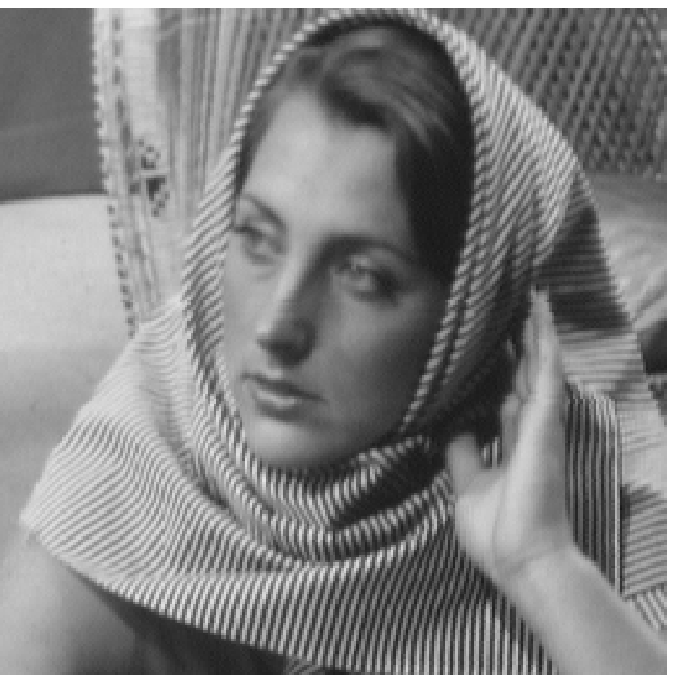}}
  \subfigure[0.05\% sampled pixels]{\label{fig:3b}\includegraphics[width=0.3\textwidth]{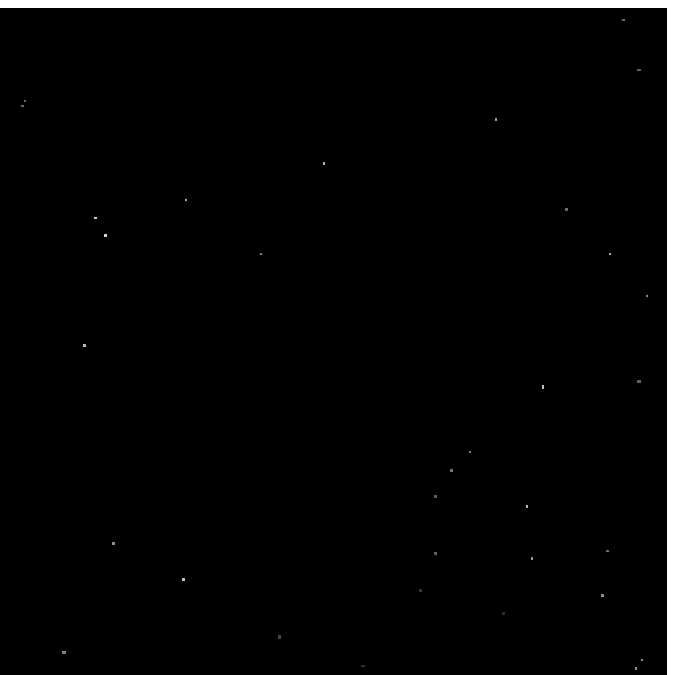}}
  \subfigure[GL \  14.58 dB ]{\label{fig:3c}\includegraphics[width=0.3\textwidth]{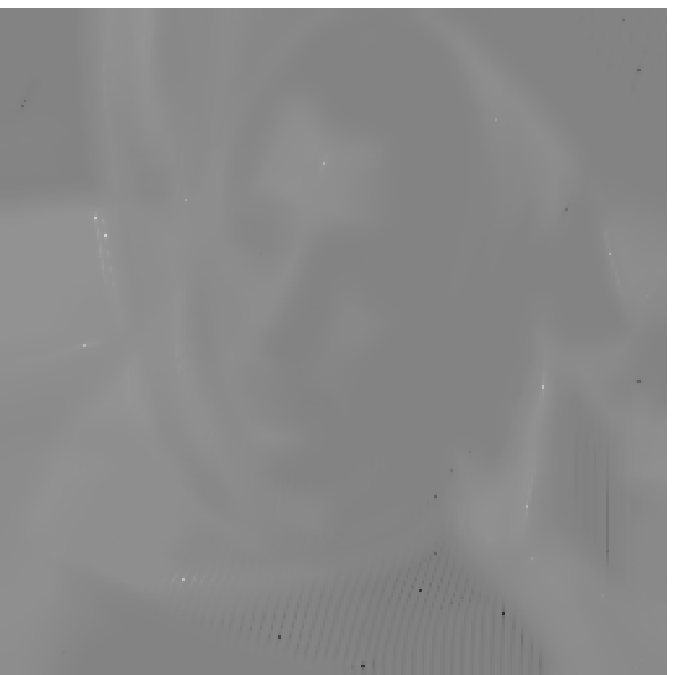}}
  \subfigure[WNLL \  18.52 dB]{\label{fig:3d}\includegraphics[width=0.3\textwidth]{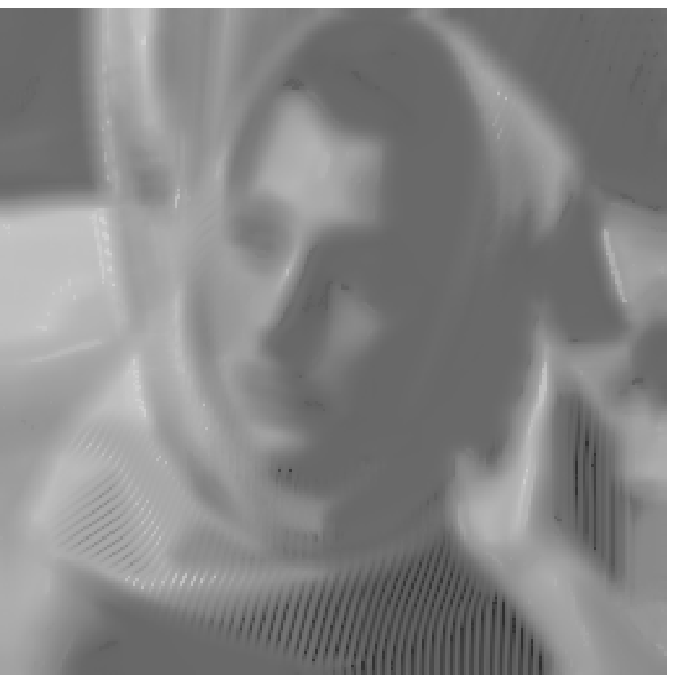}}
  \subfigure[IL ($\alpha = 0$) \  19.62 dB
  ]{\label{fig:3e}\includegraphics[width=0.3\textwidth]{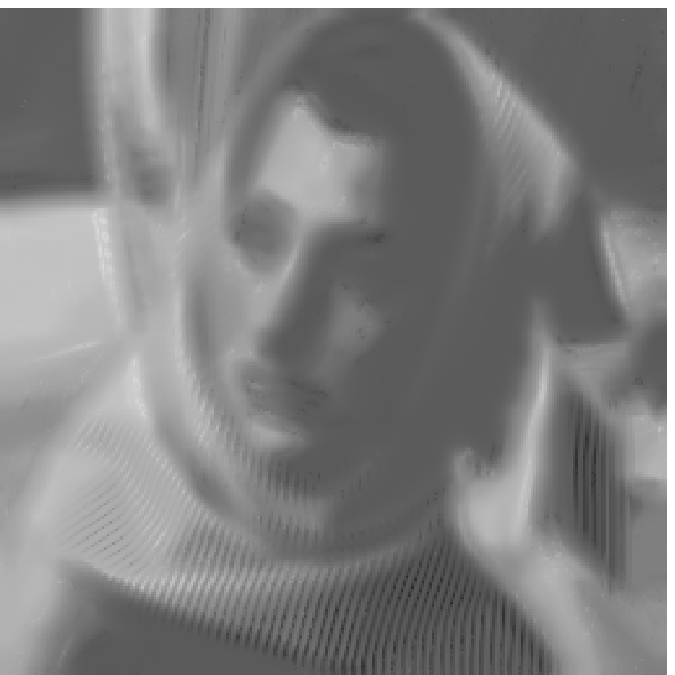}}
  \subfigure[IL ($\alpha = 10^{-5}$) \  19.38 dB ]{\label{fig:3f}\includegraphics[width=0.3\textwidth]{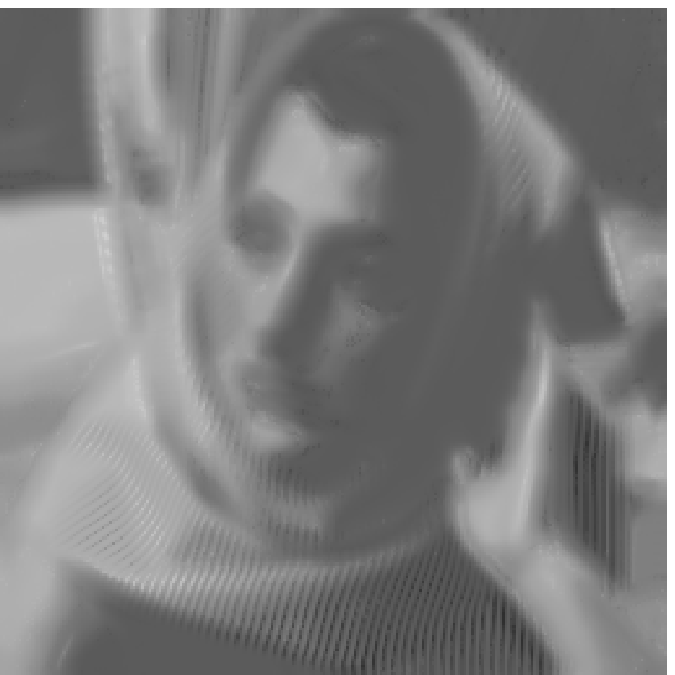}}
\caption{Inpainting results of the image Barbara. The weight is constructed using the clear image. (a) Clear image (b) 0.05\% randomly sampled pixels (c) The inpainting result of GL (d) The inpainting result of WNLL (e) The inpainting result of IL with $\alpha = 0$ (f) The inpainting result of IL with $\alpha = 10^{-5}$}
  \label{fig:fig3}
\end{figure}

A gray scale discrete image can be viewed as a function $f: \Omega =\{1,\dots,m\}\times\{1,\dots,n\} \rightarrow{[0,255]}$ or a matrix with the element $f_{ij}$ representing the intensity of the pixel in the $i^{th}$ row and $j^{th}$ column. In the task of image inpainting, values of $f$ on a subset $\Omega_0\subset \Omega$ are known, and the goal is to predict the function $f$ on $\Omega-\Omega_0$. This problem is ill-posed due to non-uniqueness of the solution and priors of the image are needed to complete the inpainting task. One of the most successful priors ever proposed is the non-local self-similarity prior \cite{NLM,BM3D,BLMSVDTV,NLH}, which means that local structures repeat many times in the whole image. In recent works, it was also observed that the inpainting problem can be solved efficiently by studying the patch manifold of the image \cite{mnfdimg,LDMM}. For each point $(i,j) \in \Omega$, we denote by $p_{ij}(f)$ the local patch at the location $(i,j)$, which is a sub-matrix of $f$ with the pixel $(i,j)$ in the center. We call the size of this sub-matrix the patch size and denote it by $p_x \times p_y$. For pixels near the edges, mirror reflection can be used to deal with the lack of data \cite{WNLL}. The set of all patches is called the patch dataset of the image.
\begin{equation*}
    \mathcal{P}(f) = \{ p_{ij}(f):(i,j)\in \Omega\}\subseteq \mathbb{R}^d, \quad d = p_x \times p_y
\end{equation*}

To inpaint the image, a commonly used assumption is that the patch dataset samples a manifold $\mathcal{M}\subset\mathbb{R}^d$, which is called the patch manifold of the image \cite{LDMM}. A coordinate function on this manifold $u:\mathcal{M}\rightarrow\mathbb{R}$ can be used to map each patch to the intensity value in the center of this patch, that is $u(p_{ij}(f)) = f(i,j)$. The image inpainting task is then formulated as an interpolation problem in $\mathbb{R}^{d}$, with the label set $\{p_{ij}(f),(i,j)\in{\Omega_0}\}$, the label values $\{u_{ij} = f(i,j),(i,j)\in{\Omega_0}\}$ and the unlabel set $\{p_{ij}(f),(i,j)\in{\Omega-\Omega_0}\}$. The interpolated values $\{u_{ij} = f(i,j),(i,j)\in{\Omega - \Omega_0}\}$, together with the sampled pixels, are then used to reconstruct the whole image. In the beginning, the label set, unlabel set, and weights are unknown. An efficient method to deal with this issue is to fill the image with random values, and update the image and weights iteratively \cite{WNLL}.

To test the efficiency of the proposed algorithm under extreme low sampling rate, we use a toy image inpainting example where the weights are constructed using the clear image. In Figures \ref{fig:3a} and \ref{fig:3b} we show the widely used image Barbara and its 0.05\% randomly sampled pixels. We choose the patch size as $p_x = p_y = 11$. For each patch, we truncate the weight to its 50 nearest neighbors searched with the approximate nearest neighbor (ANN) algorithm in the vl-feat package \cite{vlfeat}. The weight of the patches $x$ and $y$ can be computed using the clear image with $w_{xy} = \operatorname{exp}(-\Vert x - y \Vert^2 / \sigma^2)$. While in this experiment we choose $w_{xy} =({\operatorname{exp}}(-\Vert x - y \Vert^4 / \sigma(x)^4) )^2$ because it generates better results. $\sigma(x)$ is chosen to be the distance between $x$ and its $20^{th}$ nearest neighbors. 

\begin{figure}[t]
  \centering
  \subfigure[the enlarged area]{\label{fig:4a}\includegraphics[width=0.3\textwidth]{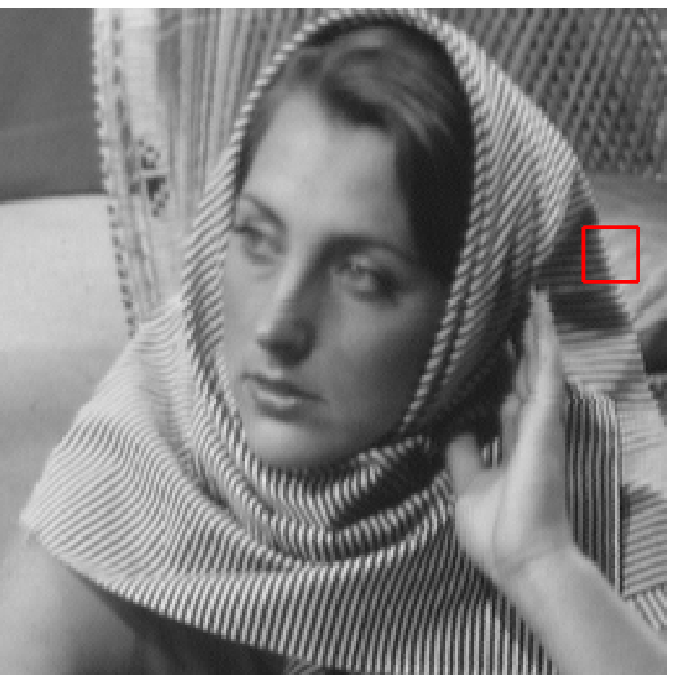}}
  \subfigure[ground truth]{\label{fig:4b}\includegraphics[width=0.3\textwidth]{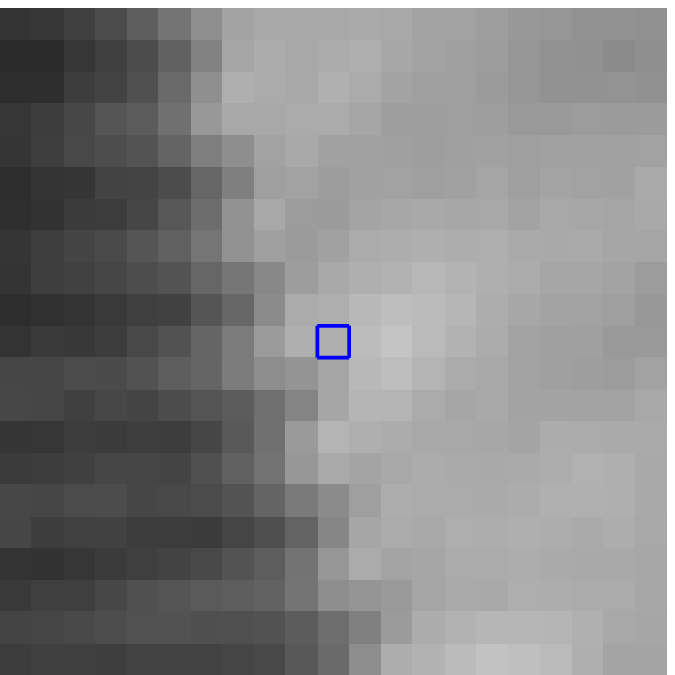}}
  \subfigure[GL \ 15.42 dB ]{\label{fig:4c}\includegraphics[width=0.3\textwidth]{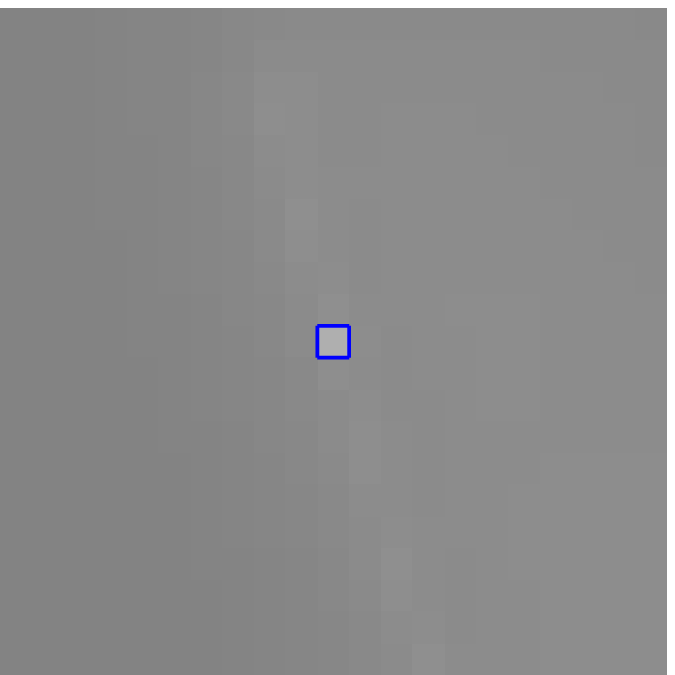}}
  \subfigure[WNLL \ 19.93 dB ]{\label{fig:4d}\includegraphics[width=0.3\textwidth]{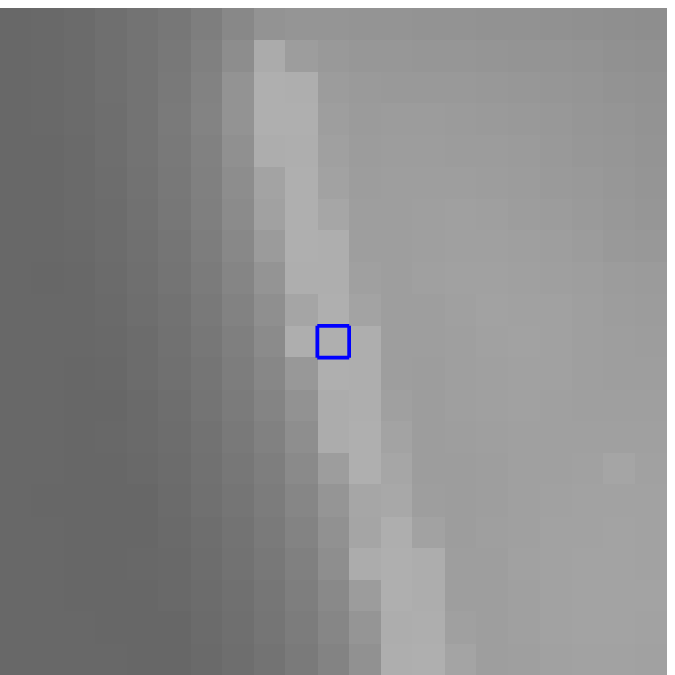}}
  \subfigure[IL ($\alpha = 0$) \  22.59 dB
  ]{\label{fig:4e}\includegraphics[width=0.3\textwidth]{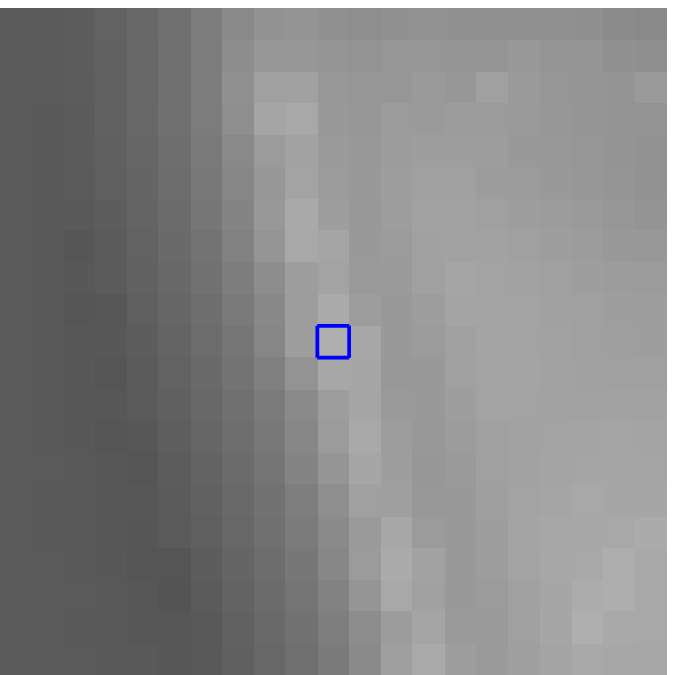}}
  \subfigure[IL ($\alpha = 10^{-5}$) \  22.62 dB ]{\label{fig:4f}\includegraphics[width=0.3\textwidth]{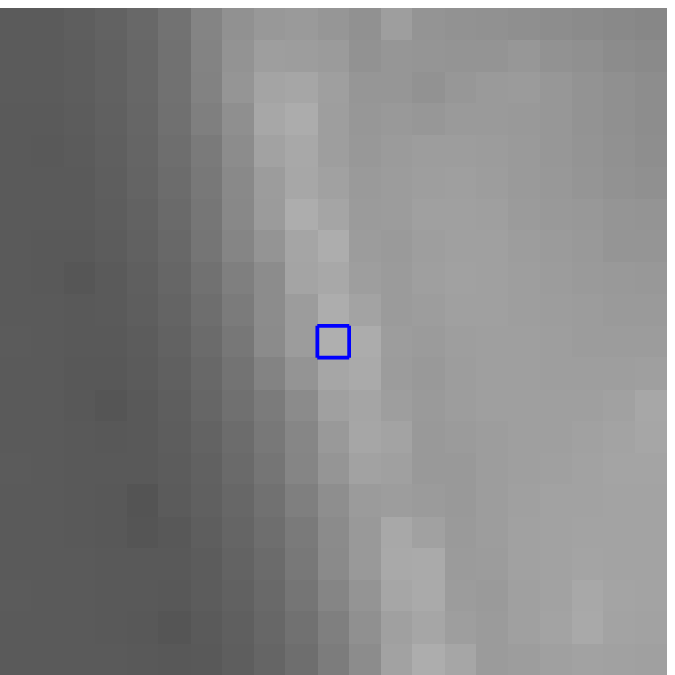}}
\caption{Zoom in of a patch of the image Barbara. (a) The enlarged area of the clear image is shown in the red box. (b) Zoom in of the clear image. There is only one sampled pixel in this region, as is shown in the blue square. (c) Zoom in of the inpainting result of GL (d) Zoom in of the inpainting result of WNLL (e) Zoom in of the inpainting result of IL with $\alpha = 0$ (f) Zoom in of the inpainting result of IL with $\alpha = 10^{-5}$}
  \label{fig:fig4}
\end{figure}

In \autoref{fig:fig3} we compare the inpainting result of the proposed infinity Laplacian model with GL and WNLL. The interpolated values of GL are not consistent with the label values. For the IL model, we choose $\alpha = 0$ (Figure \ref{fig:3e}) and $\alpha = 10^{-5}$ (Figure \ref{fig:3f}) respectively. It is shown that the IL inpainting has higher contrast at the forehead and hair than the result obtained with WNLL. We use the Peak Signal to Noise Ratio (PSNR) as a quantitive criteria to compare the results. The PSNR of the inpainted image $f$ is defined as
\begin{equation*}
    \operatorname{PSNR}(f,f^*) = 20 \log_{10}\frac{255}{\sqrt{\Vert f - f^* \Vert^2/(mn)}}
\end{equation*}
in which $f^*$ is the clear image, 255 is the possible maximum value of the intensity, $m$ and $n$ are the height and width of the image respectively. The GL model yields lower PSNR value (14.58 dB) than WNLL (18.52 dB) and IL (19.62 dB). The highest PSNR is obtained by using the IL method with $\alpha = 0$. In \autoref{fig:fig4} we compare zoom in of a $21 \times 21$ patch of these reconstruction results. The region enlarged is shown in the red box. There is only one sampled pixel in this region, which is shown in the blue square. The GL inpainting fails to capture the boundary and the interpolated values are inconsistent with the intensity of the sampled pixel. Both WNLL and IL capture the boundary, while IL reconstructs the structure better and yields higher PSNR value. 

\section{Conclusion and future work}
\label{sec:conclusion}

In this paper, we propose a novel infinity Laplacian method for the interpolation task. The infinity norm of the labeling function is minimized by introducing the non-local gradient and convergence of the discrete minimizer to the optimal continuous function is proved. The efficiency of the IL method is verified with two toy examples. In the future, we will consider other regularizations other than graph Laplacian in the IL model. We will also report the performance of the proposed method in real-world applications, such as inpainting, classification, super resolution, colorization, and denoising.

\section{Appendix}
\subsection{Auxiliary Properties of Bounded $C^2$ Domain}
	\label{subsec:auxiliary-properties}
	In this part we will prove the properties of Bounded Domains of Class $C^2$ used in our proof above.
	
	\begin{proposition}
		\label{prop:a1}
		Let $\Omega \subset \mathbb{R}^d$ be a bounded $C^2$ domain.
		For any $n\in\mathbb{N}^+$, define
		\[\tilde{\Omega}_n=\{x\in\Omega\big|B(x,\frac{1}{n})\subset\Omega\}\]
		where $B(x,r)$ represents the ball in $\mathbb{R}^d$ with center x and radius $r$.
		Then the following equation holds
		\begin{align*}
		d_n\coloneqq\sup_{x\in\Omega\backslash\tilde{\Omega}_n}dist(x,\tilde{\Omega}_n)\longrightarrow 0(n\longrightarrow\infty)
		\end{align*}
	\end{proposition}
	
\begin{proposition}
\label{prop:a2}
Let $\Omega \subset \mathbb{R}^d$ be a bounded $C^2$ domain, then there exists positive constants $C,\ r_0 > 0$ such that for all $0 < r < r_0$,
\begin{equation*}
		\nu(B(x,r)\cap \Omega) \geq C r^d,\; \forall x \in \Omega
\end{equation*}
here $B(x,r) \subset \Omega$ is an open ball centered at $x$ with radius $r$.
\end{proposition}

\begin{proof}

For all $x\in\partial\Omega$, since $\partial\Omega$ is $C^2$, we can denote $B(c(x,r_x),r_x)$ an inscribed sphere tangent to $\partial\Omega$ at $x$, where $c(x,r_x)-x=r_x\mathbf{n}_x$ and $\mathbf{n}_x$ is the unit inward normal vector of $\partial\Omega$ at $x$. We also denote $U_x$ a small open neighborhood of $x$ on $\partial\Omega$ such that $B(c(x,r_x),r_x)\cap U_x=\{x\}$. Since $\partial\Omega\backslash U_x$ is a closed set in $\mathbb{R}^d$, there exists a constant $\varepsilon>0$ such that $B(x,\varepsilon)\cap \partial\Omega\backslash U_x=\emptyset$. So if we choose $r_x<\frac{\varepsilon}{2}$, we have $B(c(x,r_x),r_x)\cap\partial\Omega=\{x\}$ which yields $B(c(x,r_x),r_x)\subset\Omega\cup\{x\}$.

Next, we hope to prove that there exists $r_0>0$ such that $B(c(x,r_0),r_0)\subset\Omega\cup\{x\}\ ,\forall x\in\partial\Omega$. We prove this with apagoge, assuming that there exists a sequence $\{x_n\}\subset\partial\Omega$ and $\{r_n\}\subset\mathbb{R}^+$ satisfied $r_n\rightarrow 0$ and $B(c(x_n,r_n),r_n)\cap\partial\Omega\varsupsetneqq\{x_n\}$. Since $\partial\Omega$ is compact in $\mathbb{R}^d$, we may assume $x_n\rightarrow x_0$ for some $x_0\in\partial\Omega$. Be same as above, choose $\varepsilon>0$ such that $B(x_0,\varepsilon)\cap \partial\Omega\backslash U_{x_0}=\emptyset$ and let $n$ be large enough so that $r_n<\frac{\varepsilon}{2}$ and $|x_n-x_0|<\frac{\varepsilon}{2}$. Then $B(c(x_n,r_n),r_n)\subset B(x_0,\varepsilon)$ which yields $B(c(x_n,r_n),r_n)\cap\partial\Omega=\{x_n\}$ and derive a contradiction.

Now we consider a large enough $n$ such that $\frac{1}{n}<r_0$. Note that
\begin{equation*}
    \tilde{\Omega}_n=\{x\big|dist(x,\Omega^c)\geq\frac{1}{n}\}
\end{equation*}
so we have $\Omega\backslash\tilde{\Omega}_n=\{x\big|dist(x,\Omega^c)\in(0,\frac{1}{n})\}$. Because $\Omega^c$ is a closed set, for all $x\in\Omega\backslash\tilde{\Omega}_n$, we can choose $y\in\Omega^c$ such that $|x-y|=dist(x,\Omega^c)<\frac{1}{n}$. Moreover, we have $y\in\partial\Omega$ and $x-y$ has the same direction with $\mathbf{n}_y$ thanks to $\partial\Omega$ is $C^2$. Let $\delta=\min\{r_0-\frac{1}{n},dist(x,\Omega^c)\}$, then $B(c(y,\frac{1}{n}+\delta),\frac{1}{n}+\delta)\subset B(c(y,r_0),r_0)\subset\Omega\cup\{y\}$. So $B(c(y,\frac{1}{n}+\delta),\frac{1}{n})\subset\Omega$, $c(y,\frac{1}{n}+\delta)\in\tilde{\Omega}_n$ and $dist(x,\tilde{\Omega}_n)\leq|x-c(y,\frac{1}{n}+\delta)|=\frac{1}{n}+\delta-dist(x,\Omega^c)\leq\frac{1}{n}$ which yields $d_n\leq\frac{1}{n}\rightarrow 0(n\rightarrow \infty)$ and we get Proposition \ref*{prop:a1}.

For Proposition \ref*{prop:a2}, define $\tilde{\Omega}_r=\{x\in\Omega\big|B(x,r)\subset\Omega\}$ for $r>0$. Be similar with the above part, we can prove for all $r<r_0$ and all $x\in\Omega\backslash\tilde{\Omega}_r$, there exists $y\in\tilde{\Omega}_{r}$ such that $|x-y|\leq r$. So there exists a constant $C\in(0,1)$ only with respect to the dimension $d$ such that
\[\nu(B(x,r)\cap\Omega)\geq\nu(B(x,r)\cap B(y,r))\geq  Cr^d\]
\end{proof}
	
\begin{proposition}
\label{prop:a3}
Let $\Omega \subset \mathbb{R}^d$ be a bounded $C^2$ domain. Then the following holds:
	\begin{equation}
	\limsup_{\delta \downarrow 0} \left\{\frac{d_\Omega(x,y)}{|x-y|}\;  : \; x,y\in \Omega, |x-y| < \delta \right\} \leq 1
	\label{eq:convex}
	\end{equation}
\end{proposition}
\begin{proof}
This proposition is introduced in \cite{ref1} Example 3.6.
\end{proof}

\subsection{Non-uniqueness of the minimizer when $\alpha = 0$}
\label{appendix:non-uniqueness}
Recall the discrete infinity Laplacian problem
\begin{equation*}
\mathop{\min}_{\{u_{m+j}\}_{j=1}^n} \Big\{\mathop{\max}_{i\in I} \big\{\sum_{j \in I} w_{ij}(u_i-u_j)^2\big\}  + \alpha \sum_{i\in I} \sum_{j \in I} w_{ij}(u_i-u_j)^2\Big\}\\
\end{equation*}
in which $I = \{1,\dots, m+n\}$. $u_{m+j}$ is the label value of $x_{m+j}$. $w_{ij}$ is the weight of the points $x_i$ and $x_j$ and $\alpha$ is a constant. When $\alpha = 0$, we use an example to show that the minimizer may not be unique. Let $S = \{x_1, x_2\}$ be the set of labeled points with label values $u(x_1) = 2$ and $u(x_2) = 0$. The set of unlabeled points is $T = \{x_3, x_4\}$. The weight matrix $\mathbf{W} = (w_{ij})_{4\times4}$ is given as
\begin{equation*}
\mathbf{W} = 
\begin{bmatrix}
1 & 1/3 & 1/2 & 0\\ 
1/3 & 1 & 1/2 & 1/2\\
1/2 & 1/2 & 1 & 0\\
0 & 1/2 & 0 & 1\\
\end{bmatrix}   
\end{equation*}
The set of minimizers of this example is given by $\{(x_3,x_4) | x_3 = 1, -1 \leq x_4 \leq 1\}$.

\subsection{Solving the sub-problem for updating $\mathbf{D}$}
\label{subproblem_update_D}
To update $\mathbf{D}$, we need to solve the problem of the following type.
\begin{equation*}
\mathop{\min}_{\mathbf{x}} \Big\{ \mathop{\max}_{i \in \{1,\dots,N\}} \{x_i^2\} + \sum_{i=1}^{N} a_i (x_i - c_i)^2 \Big\}
\end{equation*}
in which $\mathbf{x}=(x_1,\dots,x_N) \in \mathbb{R}^N$, $a_i > 0$, $c_i \geq 0$. We denote this objective function by $\phi(\mathbf{x})$. Note that $\phi$ is strictly convex on $\mathbb{R}^N$ with a lower bound 0, so there exists a unique minimizer $\mathbf{x}^* = (x_1^*,\dots,x_n^*)$. We denote $x^* = \max_i\{x_i^*\}$. Note that $\mathbf{x}^* \in \prod _{i=1}^n [0,c_i]$, it has to be the case that $x_i^* = x_j^*$ whenever $c_i = c_j$. Without loss of generality, it suffice to consider the case of $c_1 > c_2 > \dots > c_N$. In the following we prove that the output of Algorithm~\ref{alg:max_L2_alg} is the minimizer of this problem. We begin with a trivial fact.
\begin{lemma}
\label{trivial_fact}
Let $f_1(x) = a_1(x - y_1)^2 + c_1$ and $f_2(x) = a_2(x - y_2)^2 + c_2$ be two quadratic functions with $a_1 > 0$, $a_2 > 0$, and $y_1 \neq y_2$. Let $x^*$ be the minimizer of $f_1(x) + f_2(x)$, then $x^* \in (\mathop{\min}\{y_1,y_2\},\mathop{\max}\{y_1,y_2\})$.
\end{lemma}

We define quadratic functions $f_i(x) = x^2 + a_i(x-c_i)^2$ and denote the corresponding minimizers as $x_{f_i}^* = a_i c_i / (1 + a_i)$. We also denote by
\begin{equation*}
    x_{F_t}^* = \frac{\sum_{i=1}^t a_i c_i}{ 1 + \sum_{i=1}^t a_i}
\end{equation*}
the minimizer of $F_t(x) = x^2 + \sum_{j=1}^t a_j(x - c_j)^2$ for $t = 1,\dots,N$.

Let $x_f^* = \mathop{\max}_i \{ x_{f_i}^* \}$. We claim that $x^* \geq x_f^*$. If not, choose $k$ such that $x_{f_k}^* = \mathop{\max}_i \{ x_{f_i}^* \}$ and we have
\begin{equation*}
\begin{aligned}
&\phi(x_1^*,\dots,x_{f_k}^*,\dots,x_n^*) - \phi(\mathbf{x}^*)\\
= &x_{f_k}^* + a_k(x_{f_k}^* - c_k)^2 - (x^*)^2 - a_k(x_k^* - c_k)^2\\
\leq & x_{f_k}^* + a_k(x_{f_k}^* - c_k)^2 - (x_k^*)^2 - a_k(x_k^* - c_k)^2\\
= & f_k(x_{f_k}^*) - f_k(x_k^*)\\
< &0
\end{aligned}
\end{equation*}
which is a contradiction. Therefore $c_i \leq x_f^*$ implies $x_i^* = c_i$.

For the special case of $x_{f_1}^* \geq c_2 > \dots > c_N$, we have $x_1^* = x^*$ and the solution is $\mathbf{x}^* = (x_{f_1}^*,c_2,c_3,\dots,c_N)$. For the general case of $c_1 > c_2 > \dots > c_t > x_f^* \geq c_{t+1} > \dots > c_N $, clearly $x_j^* = c_j$ for $j > t$ and it suffices to find $\{x_1^*,\dots, x_t^*\}$. This problem breaks down into several cases.

If $x^* \in [x_f^*, c_t]$, we have $x_1^* = x_2^* \dots = x_t^* = x^*$ and the objective function reduces to a quadratic function
\begin{equation*}
\phi_{[x_f^*,c_t]}(x) = x^2 + \sum_{i=1}^t a_i(x_i - c_i)^2
\end{equation*}
in which the subscript $[x_f^*,c_t]$ emphasizes the domain of $x$. The abscissa of the vertex of this quadratic function is $x_{F_t}^*$, which is bigger than $x_f^*$. In fact, for each $1\leq j \leq t$, we have
\begin{equation*}
\begin{aligned}
\frac{\sum_{i=1}^t a_i c_i}{ 1 + \sum_{i=1}^t a_i} &= \frac{a_j c_j}{1 + \sum_{i=1}^t a_i} + \frac{\sum_{i \neq j} a_i c_i}{1 + \sum_{i=1}^t a_i}\\
&> \frac{a_j c_j}{1 + \sum_{i=1}^t a_i} + \frac{\sum_{i \neq j} a_i}{1 + \sum_{i=1}^t a_i} x_f^*\\
&\geq \frac{a_j c_j}{1 + \sum_{i=1}^t a_i} + \frac{\sum_{i \neq j} a_i}{1 + \sum_{i=1}^t a_i}\frac{a_j c_j}{1+a_j}\\
&= \frac{a_j c_j}{1+a_j}\\
& = x_{f_j}^* 
\end{aligned}
\end{equation*}
On the other hand, we have
\begin{equation*}
\label{c_t}
\frac{\sum_{i=1}^t a_i c_i}{ 1 + \sum_{i=1}^t a_i} - c_t =
\frac{\sum_{i=1}^{t-1} a_i(c_i - c_t) - c_t}{1+\sum_{i=1}^t a_i}
\end{equation*}
When $\sum_{i=1}^{t-1} a_i(c_i - c_t) \leq c_t$, the abscissa of the vertex of $\phi_{[x_f^*,c_t]}(x)$ lies in the interval $[x_f^*, c_t]$ and the optimal solution is $x_1^* = x_2^* = \dots = x_t^* = x_{F_t}^*$. Now we claim that $\sum_{i=1}^{t-1} a_i(c_i - c_t) \leq c_t$ implies $x^* \in [x_f^*, c_t]$ and prove it with reduction to absurdity. Suppose for some $1 \leq j \leq t-2$, $x^* \in [c_{j+1},c_{j}]$, then it has to be the case that $x_k^* = c_k$, for $j < k \leq t$ and $x_1^* = x_2^* = \dots = x_j^*$. Now the objective function turns out to be
\begin{equation*}
\phi_{[c_{j+1},c_j]}(x) = x^2 + \sum_{i=1}^j a_i(x_i - c_i)^2
\end{equation*}
in which the subscript $[c_{j+1},c_j]$ emphasizes the domain of $x$. The abscissa of the vertex of this quadratic function is less than or equal to $c_t$. (Suppose not, adding the terms $a_k(x-c_k)^2$ with $k = j+1,\dots,t$ one by one to $\phi_{[c_{j+1},c_j]}(x)$ and we obtain $\phi_{[x_f^*,c_t]}(x)$. By Lemma \ref{trivial_fact}, the abscissa of the vertex of $\phi_{[x_f^*,c_t]}(x)$ should be greater than $c_t$, which is a contradiction with \eqref{c_t} and our assumption. Based on this observation, we see that $\phi_{[c_{j+1},c_j]}(x)$ increases strictly on $[x_{j+1}, x_{j}]$ and the optimal solution is $x_1^* = \dots = x_{j+1}^* = x^* = c_{j+1}$. Now it is also the case that $x^* \in [x_{j+2},x_{j+1}]$, and the same argument shows that $x_1^* = \dots = x_{j+1}^* = x_{j+2}^* = x^* = c_{j+2}$, a contradiction. Suppose $x^* \in [c_{t},c_{t-1}]$, the same argument shows that $ x_1^* = \dots = x_t^* = x^* = c_t$. It is also the case that $x^* \in [x_f^*,x_{t}]$. By our previous discussion, it should be the case that $x_1^* = x_2^* = \dots = x_t^* = x_{F_t}^*$, which is a contradiction. To sum up, when $\sum_{i=1}^{t-1} a_i(c_i - c_t) \leq c_t$, the unique minimizer is $x_1^* = \dots = x_t^* = x_{F_t}^*$ and $x_{j}^* = c_j$ for $j > t$.

When $\sum_{i=1}^{t-1} a_i(c_i - c_t) > c_t$, $\phi_{[x_f^*,c_t]}(x)$ decreases strictly on $[x_f^*,c_t]$. If $x^* \in [x_f^*,c_t]$, then it has to be the case that $x_1^* = \dots = x_t^* = x^* = c_t$. Now it is also the case that $x^* \in [c_t, c_{t-1}]$, so $x^* = c_t$ should be the minimizer of $\phi_{[c_t,c_{t-1}]}^* = \sum_{i=1}^{t-1}a_i(x_i-c_i)^2 + x^2$. However, by Lemma \ref{trivial_fact}, the abscissa of its vertex $x_{F_{t-1}}^*$ is greater than $c_t$, which is a contradiction. Thus $x^* \notin [x_f^*, c_t]$. If $x^* \in [c_t, c_{t-1}]$, we have $x_{F{t-1}}^* < c_{t-1}$ when $\sum_{i=1}^{t-2} a_i(c_i - c_{t-1}) \leq c_{t-1}$. In this case, $x_1^* = \dots = x_{t-1}^* = x_{F_{t-1}}^*$, $x_t^* = c_t$. Similar argument leads to the fact that $\sum_{i=1}^{t-2} a_i(c_i - c_{t-1}) \leq c_{t-1}$ if and only if $x^* \in [c_t, c_{t-1}]$ on the premise that $\sum_{i=1}^{t-1} a_i(c_i - c_t) > c_t$.

Continuing this process, we find the greatest positive integer $S \leq t$, such that $\sum_{i=1}^{s-1} a_i(c_i - c_s) > c_s$ for all $S < s \leq t$ and $\sum_{i=1}^{S-1} a_i(c_i - c_S) \leq c_S$. Then the optimal solution is given by $x_1^* = \dots = x_S^* = x_{F_S}^*$ and $x_p^* = c_p$ for all $S < p \leq N$. The integer $S$ exists because there is some positive integer $k$ such that $x^* \in [c_{k+1},c_k]$. When $S = 1$ and $s = 2$, the inequality $\sum_{i=1}^{s-1} a_i(c_i - c_s) > c_s$ reduces to $x_{f_1}^* > c_2$, which is exactly the special case we discuss in the beginning. For fast implementations, we use $\sum_{i=1}^{s-1} a_i(c_i - c_s) = \sum_{j=1}^{s-1} \sum_{i=1}^j a_i(c_j - c_{j+1})$ in Algorithm~\ref{alg:max_L2_alg}.

\begin{acknowledgements}
The authors would like to thank Peng Zhang and Tangjun Wang for valuable discussion.
\end{acknowledgements}

\section*{Conflict of interest}
The authors declare that they have no conflict of interest.

\section*{Code availability}
The open source code of our algorithm is available at \url{https://cloud.tsinghua.edu.cn/d/4d3d9e828be44426b5b8/}

\bibliographystyle{spmpsci}      
\bibliography{references}     

\end{document}